\documentclass[a4,reqno,11pt]{amsart}  

\usepackage{gensymb}
\usepackage{amsmath}
\usepackage{amssymb} 
\usepackage{mathtools}   % loads »amsmath«
\usepackage{graphicx}
\usepackage{graphics}  
\usepackage{color}
\usepackage{verbatim} 
\usepackage[normalem]{ulem} 
\usepackage[mathscr]{eucal}
\usepackage{upgreek}
\usepackage{enumerate} 
\usepackage{cite}
\usepackage{amsmath}
\usepackage{amsfonts}
\usepackage{multirow}
\usepackage{amssymb}
\usepackage{bbm}
\usepackage{cancel}
\usepackage{graphicx}
\usepackage{multicol}
\usepackage[utf8]{inputenc}
\usepackage{framed}
\usepackage{setspace}
\usepackage{amsthm}
\usepackage{hyperref}
\usepackage{caption}
\usepackage{subcaption}
\usepackage{bbm}
\usepackage{graphicx}
\usepackage{dsfont}
\usepackage{tikz}
\usepackage{mwe}
\usepackage[titletoc]{appendix}
\numberwithin{equation}{section}  
\usepackage[refpage,noprefix]{nomencl}
\usepackage{nomencl} 
\newtheorem{theorem}{Theorem}[section] 
\newtheorem{lemma}[theorem]{Lemma} 
\newtheorem{proposition}[theorem] {Proposition} 
\newtheorem{cor}[theorem]  {Corollary} 
\newtheorem{remark}[theorem]  {Remark} 
\newtheorem{definition}[theorem] {Definition} 
 
\newtheorem{assump}[theorem]{Assumption}

\theoremstyle{definition}
\newtheorem{example}[theorem] {Example}

 \usepackage{float}
\restylefloat{figure}

\DeclareMathAlphabet{\mathpzc}{OT1}{pzc}{m}{it}

\DeclarePairedDelimiter{\abs}{\lvert}{\rvert}
\DeclarePairedDelimiter{\norm}{\lVert}{\rVert}

\renewcommand{\L} {\Lambda} %

\def\d{\delta} 
\newcommand{\e} {\varepsilon}

\newfam\Bbbfam 
\font\tenBbb=msbm10 
\font\sevenBbb=msbm7 
\font\fiveBbb=msbm5 
\textfont\Bbbfam=\tenBbb 
\scriptfont\Bbbfam=\sevenBbb 
\scriptscriptfont\Bbbfam=\fiveBbb 
 
\newcommand{\B}     {\mathbb{B}} 
  
\newcommand{\R}     {\mathbb{R}} 
\newcommand{\Z}     {\mathbb{Z}} 
\newcommand{\N}     {\mathbb{N}} 
\renewcommand{\P}   {\mathbb{P}} 
 
\newcommand{\E}     {\mathbb{E}}

 \newcommand{\floor}[1]{\left\lfloor #1 \right\rfloor}

\def\1{{\mathchoice {1\mskip-4mu\mathrm l}      % Blackboard bold 1 
{1\mskip-4mu\mathrm l} 
{1\mskip-4.5mu\mathrm l} {1\mskip-5mu\mathrm l}}} 
 
\def\comment#1{} 
\newtheoremstyle{thm}{2ex}{2ex}{\itshape\rmfamily}{} 
{\bfseries\rmfamily}{}{1.7ex}{} 
 
\newtheoremstyle{rem}{1.3ex}{1.3ex}{\rmfamily}{} 
{\itshape\rmfamily}{}{1.5ex}{} 
 
\newcommand{\bB} {\boldsymbol{B}}

 \newcommand{\bmu}{\boldsymbol{\mu}}

\newcommand{\Ocal}   {{\mathcal O }}

\newcommand{\Pfrak}{\mathfrak{P}}
\newcommand{\Zfrak}{\mathfrak{Z}}

\newcommand{\hy}{{\textsf{HY}}}

 \newcommand{\ex}{{\rm e}} 
 
\renewcommand{\d}{{\rm d}}

\newcommand{\supp}{{\operatorname {supp}}}

\newcommand{\Exp}{\mathscr{E}\kern-0.2mm{\operatorname{xp}}}
\newcommand{\Log}{\mathscr{L}\kern-0.2mm{\operatorname{og}}}

\renewcommand{\emptyset} {\varnothing} 
 
\newcommand{\p} {\partial}

\newcommand\NoBlackBoxes{\global\overfullrule0pt}
\NoBlackBoxes

\newcommand{\redline}{\raisebox{2pt}{\tikz{\draw[very thick, red](0,0) -- (1,0);}}}
\newcommand{\bluedottedline}{\raisebox{2pt}{\tikz{\draw[very thick, dotted, blue](0,0) -- (1,0);}}}

\renewcommand{\p}{\mathfrak{p}}

\renewcommand{\c}{{\mathfrak{c}}}

\newcommand{\pmf}{{\mathrm{PMF}}}

\renewcommand{\hy}{{\mathrm{HYL}}}
\newcommand{\gmf}{{\mathrm{GMF}}}
\newcommand{\Phy}{\P^\hy}
\newcommand{\bN}{\overline{N}}

\newcommand{\rhob}{\bar{\rho}}

\newcommand{\brho}{\boldsymbol{\rho}}
\renewcommand{\bmu}{\boldsymbol{\mu}}
\newcommand{\Hpmf}{\mathsf{H}^\pmf}

\newcommand{\Hgmf}{\mathsf{H}^\gmf}

\newcommand{\Ham}{\mathsf{H}}

\newcommand{\Ns}{N^{\mathrm{short}}}
\newcommand{\Nl}{N^{\mathrm{long}}}
\newcommand{\Mfrak}{\mathfrak{M}}

\renewcommand{\min}{\mathrm{min}}
\newenvironment{proofsect}[1] 
{\vskip0.1cm\noindent{\scshape #1.}\hskip0.5cm} 

\newcommand{\GCanonFree}[2]{\P_{#1,\beta,#2}}
\newcommand{\GCanonHY}[2]{\P^\hy_{#1,\beta,#2}}
\newcommand{\CanonFree}[2]{\P^{({\rm Can})}_{#1,\beta,#2}}
\newcommand{\CanonHY}[2]{\P^{({\rm Can},\hy)}_{#1,\beta,#2}}

\makeindex

\begin{document}

\title[\hfill Formation of infinite loops for an interacting bosonic loop soup\hfill]
{Formation of infinite loops for an interacting bosonic loop soup}

\author{Matthew Dickson \and Quirin  Vogel}
\address[Matthew Dickson]{Ludwig--Maximilians--Universität München, Mathematisches Institut, Theresienstr.\ 39, 80333 München, Germany; Email: dickson@math.lmu.de}
\address[Quirin  Vogel]{NYU-ECNU Institute of Mathematical Sciences at NYU Shanghai, 3663 Zhongshan Road North, Shanghai, 200062, China; Email: qtv203@nyu.edu}

% \email{qtv203@nyu.edu}
%\thanks{}

\thanks{}
 \maketitle  
%    \date is required; it is the date received by the editor.
%\date{September 25, 2012}

%    The 2010 edition of the Mathematics Subject Classification is
%    now available.  If you are citing a classification from the
%    new scheme, use the following input coding instead.
%\subjclass[2010]{Primary }
 
\begin{abstract}
We compute the limiting measure for the Feynman loop representation of the Bose gas for a non mean-field energy. As predicted in previous works, for high densities the limiting measure gives positive weight to random interlacements, indicating the quantum Bose--Einstein condensation. We prove that in many cases there is a shift in the critical density compared to the free/mean-field case, and that in these cases the density of the random interlacements has a jump-discontinuity at the critical point.
\end{abstract}
%\dedication{Dedication text (use \\[2pt] for line break if necessary)} 

% \tableofcontents

\section{Introduction}
Ever since Feynman introduced his representation of the Bose gas as a soup of interacting loops in \cite{feynman1948space}, there has been a continued interest in its properties. Amongst the wide variety of interesting questions, our investigation focuses on the macroscopic behaviour of the Bose gas as the particle density varies. Bose and Einstein (in \cite{bose1924plancks} and \cite{einstein2005quantentheorie}) predicted that above a certain density, a macroscopic fraction of the bosons aggregate into a single quantum state, commonly referred to as the \textit{Bose--Einstein condensate}. Feynman in \cite{feynman1953atomic} gave arguments linking the formation of the condensate to that of macroscopic loops. Much research has been undertaken in that direction, especially focusing on the induced measure on permutations. See, for example, \cite{suto1993percolation,0305-4470-35-33-303,ueltschi2006feynman,betz2009spatial,adams2011free,adams2021large,AV19} amongst many others.

The specific question we are interested in can be informally stated as:\begin{center}
    \textit{What is the limiting\footnote{in the thermodynamic sense} measure governing the Feynman loop representation of the Bose gas?}
\end{center} 
There has been some progress in this direction (especially in the aforementioned papers), however the stochastic process of \textit{random interlacements} used to describe the limiting state has only been introduced a few years ago, see \cite{sznitman2010vacant}. The work \cite{armendariz2019gaussian} was the first to draw a rigorous connection between random interlacements and the Bose gas. In that publication the authors showed that the superposition of the bosonic loop process on the whole space and the random interlacements gives the same distribution as the random permutations described in \cite{macchi1975coincidence}. Note that in \cite{armendariz2019gaussian} the interactions between the different loops were neglected. In \cite{vogel2021emergence}, the author proved that by taking the thermodynamic limit along boxes of diverging volume, the limiting process is indeed given by the superposition of the random interlacements and the bosonic loop soup. That paper considered the case of the free Bose gas as well as the mean-field case, where the interaction energy between loops is given by the square of the total particle number.

\begin{figure}[ht]
    \centering
    \begin{subfigure}[t]{0.45\textwidth}
    \centering
    \begin{tikzpicture}[scale=1.2]
    \draw[->] (0,0) -- (4,0) node[below]{Total density $\rho$};
    \draw[->] (0,0) -- (0,4)node[above]{Limiting density};
    \draw[loosely dashed] (0,0) -- (4,4);
    \draw[very thick, dotted, blue] (0,0) -- (1.5,1.5);
    \draw[very thick, dotted, blue] (1.5,1) to [out=290,in=180] (4,0.1);
    \draw[very thick, red] (0,0) -- (1.5,0);
    \draw[very thick, red] (1.5,0.5) to [out=70,in=225] (4,3.9);
    \draw[loosely dashed] (1.5,0) node[below]{$\rho^\hy_c$} -- (1.5,4);
    \draw[loosely dashed] (2.5,0) node[below]{$\rho_{\mathrm{c}}$} -- (2.5,4);
    \end{tikzpicture}
    \caption{Partial HYL model, $d=3,4$.}
    \end{subfigure}
    \hfill
    \begin{subfigure}[t]{0.45\textwidth}   
    %\hspace{-10mm}
    \centering
    \begin{tikzpicture}[scale=1.2]
    \draw[->] (0,0) -- (4,0) node[below]{Total density $\rho$};
    \draw[->] (0,0) -- (0,4)node[above]{Limiting density};
    \draw[loosely dashed] (0,0) -- (4,4);
    \draw[very thick, dotted, blue] (0,0) -- (2.5,2.5);
    \draw[very thick, dotted, blue] (2.5,2.5) -- (4,2.5);
    \draw[very thick, red] (0,0) -- (2.5,0);
    \draw[very thick, red] (2.5,0) -- (4,1.5);
    \draw[loosely dashed] (1.5,0) node[below]{$\rho^\hy_c$} -- (1.5,4);
    \draw[loosely dashed] (2.5,0) node[below]{$\rho_{\mathrm{c}}$} -- (2.5,4);
    \end{tikzpicture}
    \caption{Mean-field and Free model, $d\geq 3$.}    
    \end{subfigure}
    \caption{The density of the random interlacements/condensate (red, \protect\redline) and the finite loops/bulk (blue, \protect\bluedottedline) for various models in certain dimensions.}
    \label{figure1}
\end{figure}
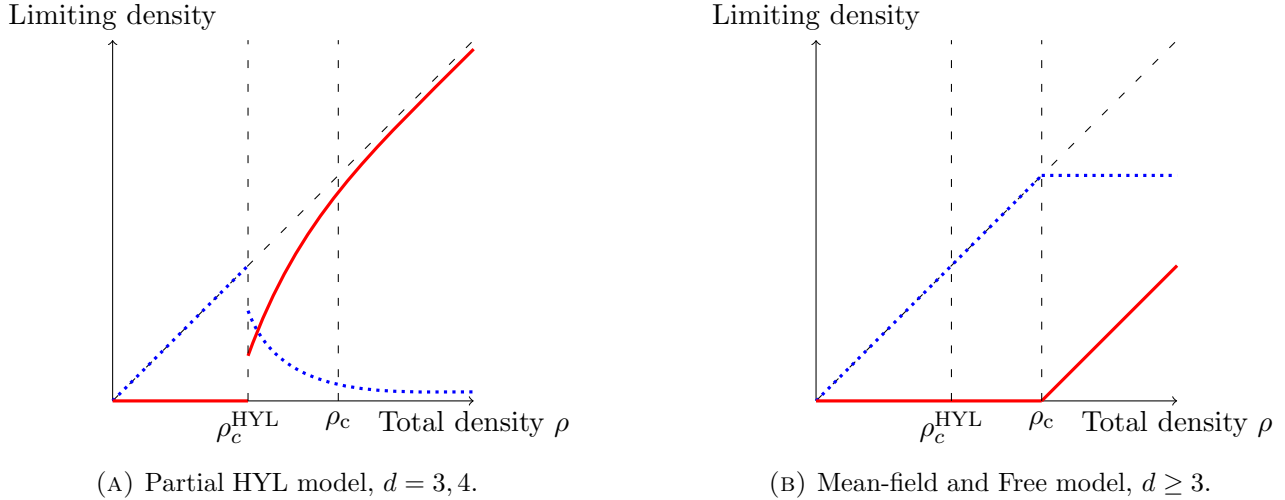

In this work, we consider the (partial) HYL energy (named after Huang, Yang and Luttinger), a non mean-field interaction between loops, inspired by the publications \cite{huang1957imperfect,lewis1986bosons,adams2021explicit}, mimicking the repulsion between particles. Such Bose soups behave qualitatively differently to the free and mean-field cases. While in the latter, there is a \textit{continuous} transition in the density of random interlacements as we cross the critical density $\rho_{\mathrm{c}}>0$ of the Bose gas, the former has a \textit{jump discontinuity} in low dimensions, see Figure \ref{figure1}. We comment in greater length on this important result in Section \ref{subsectioncondensatediscont}. Furthermore, the value of the critical density changes, compared to the free and mean-field cases.
Our proof is based on the approach of combining large deviation theory in the heavy-tail regime, based on recent estimates such as \cite{berger2019notes}, with order-one large deviation results from \cite{bahadur1960deviations,martin1982laplace} for the contribution from the exponential part. Both the heavy-tail and the exponential parts contribute to the condensate. Following \cite{bratteli2003operator}, the first part can be interpreted as being due to the quantum statistics of the bosons, while the second contribution is due to the repulsive force between particles at high densities.

The truncation and asymptotic expansion of the free large deviation rate function near the critical point is crucial in the analysis of the HYL-energy. Once we have shown that macroscopic loops exist with the correct density, we can refer back to \cite{vogel2021emergence} where the convergence of the macroscopic loops to the random interlacements is proven. As a by-product of our proofs, we get the order-one asymptotics of the partition functions.

Whilst \cite{vogel2021emergence} proved the emergence of macroscopic loops for a model on $\Z^d$, the arguments and proofs are sufficiently robust to easily apply to a model in $\R^d$. Working in the continuum has the additional advantage of aligning with previous work on random interlacements in \cite{armendariz2019gaussian} and on the interaction energy in \cite{adams2021explicit}.

\subsection*{The Models}

We first describe the (free) \emph{grand canonical ensemble} with inverse temperature $\beta>0$ and chemical potential $\mu\leq 0$. The configuration space of each individual finite loop is
\begin{equation}
    \Gamma_\mathrm{F} = \bigcup^{\infty}_{j=1}\left\{\omega:\left[0,\beta j\right]\to \R^d, \omega\text{ is continuous}, \text{ and } \omega\left(0\right)=\omega\left(\beta j\right)\right\}.
\end{equation}
Given $\omega\in  \Gamma_\mathrm{F}$, we write $\ell\left(\omega\right)=j$ if $\omega:\left[0,\beta j\right]\to \R^d$. Such a loop can be thought to represent $\ell\left(\omega\right)\in\N$ particles at inverse temperature $\beta$. Given a loop with length $\ell\left(\omega\right)=j$, it is distributed according to a Brownian bridge with a time horizon $\beta j$ and co-incident start and end points.

The grand canonical ensemble on Borel measurable $\Lambda\subset\R^d$ is then given as a Poisson point process on $\Gamma_\mathrm{F}$ (with a given intensity measure). We let $\P^t_{x,y}$ denote the un-normalised Brownian bridge measure with time horizon $t>0$ and start and end points $x,y\in\R^d$. Then the intensity measure is the bosonic loop measure, given by
\begin{equation}
    M^B_{\Lambda,\beta,\mu}\left(A\right) =\int_\L \sum_{j\geq 1}\frac{1}{j}\ex^{\beta\mu j} \P^{\beta j}_{x,x}\left(A\right)\d x
\end{equation}
for each measurable $A\subset  \Gamma_\mathrm{F}$. We denote the law of this (free) grand canonical ensemble on $\Lambda$ as $\GCanonFree{\L}{\mu}$.

Recall $\ell\left(\omega\right)$ denotes the number of particles associated with the loop $\omega$. Let $\eta_\Lambda$ be a locally finite counting measure on $ \Gamma_\mathrm{F}$ restricted to loops with $\omega(0)\in\Lambda$. Then we define $N_\Lambda$ to be the total particle number on $\Lambda$, that is
\begin{equation}
    N_\Lambda(\eta) = \sum_{\omega\in \eta_{\Lambda}}\ell\left(\omega\right).
\end{equation}
For density $\rho>0$, we define the (free) canonical measure $\CanonFree{\Lambda}{\rho}$ as the grand canonical measure $\GCanonFree{\Lambda}{0}$ conditioned on $N_\Lambda = \lfloor\rho\abs*{\Lambda}\rfloor$.

To define the Hamiltonian with which we are primarily concerned, we introduce two interaction strengths $a>b>0$ and a loop length scale $q_{\Lambda}\in\N$ such that $q_{\Lambda} = o\left(\abs*{\Lambda}\right)$ as $\abs*{\Lambda}\to\infty$. Then
\begin{equation}
    \label{eqn:partialHYLinteraction}
    \Ham_{\Lambda}\left(\eta_{\Lambda}\right) = \Ham_{\Lambda,a,b}\left(\eta_{\Lambda}\right) = \frac{a}{2\abs*{\Lambda}}N_\Lambda^2 - \frac{b}{2\abs*{\Lambda}}\sum_{k\geq q_\Lambda}k^2\#\left\{\omega\in\eta_\Lambda:\ell\left(\omega\right)=k\right\}^2.
\end{equation}
We call this the partial HYL Hamiltonian (`partial' because the counter-term only sees cycles longer than the scale $q_\Lambda$). A similar looking ``sum-of-the-squares" counter-term was introduced by \cite{huang1957imperfect} in the process of approximating the hard-sphere interaction for a gas of bosons. Their counter-term picked out all the momentum eigenstates of a quantum Bose gas rather than the cycle lengths greater than some diverging cut-off in a loop soup as ours does. Therefore there is no a priori reason to suppose that the two models are related. Nevertheless it has been proven in \cite{adams2021explicit} that this partial HYL Hamiltonian on loops produces precisely the same thermodynamic pressure as the original momentum eigenstate based model (whose pressure was derived in \cite{van1990pressure}). We also prove in Proposition~\ref{propositionfullvspartial} that if there is no cut-off, then this equality no longer holds. More details and discussion about the interaction can be found in Section~\ref{HYLDiscussion}.

The Hamiltonian influences the loop soup's behaviour via a tilting of the non-interacting measure. For finite volume $\Lambda$, we define the grand-canonical measure with HYL interaction, $\Phy_{\Lambda,\beta,\mu}$, by its Radon--Nikodym derivative with respect to the non-interacting measure. Specifically, given a chemical potential $\mu\in\R$ we have
\begin{equation}
\label{eqn:GCHYLradonnikodymDerivative}
    \frac{\d \Phy_{\Lambda,\beta,\mu}}{\d \GCanonFree{\Lambda}{0}}\left(\eta\right) = \frac{1}{Z^\hy_{\Lambda,\beta,\mu}}\exp\left(-\beta \Ham_{\Lambda}\left(\eta\right)+\mu\frac{1}{\abs*{\Lambda}}N_\Lambda\left(\eta\right)\right),
\end{equation}
where $Z^\hy_{\Lambda,\beta,\mu}$ is the partition function that normalises $\GCanonHY{\Lambda}{\mu}$. Under the measure $\GCanonFree{\Lambda}{0}$ the particle number $N_\Lambda$ is almost surely finite, and to avoid ambiguity when $\mu>0$ we set $\frac{\d \Phy_{\Lambda,\beta,\mu}}{\d \GCanonFree{\Lambda}{0}}=0$ when $N_\Lambda = +\infty$. It is worth noting that whilst the counter-term in $\Ham_\Lambda$ is itself attractive, the condition $a>b$ ensures that the Hamiltonian as a whole is bounded from below by some quadratic function of the particle number $N_\Lambda$ and is repulsive. It is also this quadratic feature which means that it is possible to include $\mu>0$ in the parameter space for this model whilst it was not possible for the free grand canonical ensemble.

Our main result is concerned with a canonical version of this interacting measure. Given $\GCanonHY{\Lambda}{\mu}$, we define the canonical measure with HYL interaction, $\CanonHY{\Lambda}{\rho}$
\begin{equation}
\label{eqn:canonHYLRadonNykodym}
    \frac{\d \CanonHY{\Lambda}{\rho}}{\d \CanonFree{\Lambda}{\rho}}\left(\eta\right) = \frac{1}{Z^{(\mathrm{Can, HYL})}_{\Lambda,\beta,\rho}}\ex^{-\beta \Ham_{\Lambda}\left(\eta\right)}= \frac{\ex^{-\beta \Ham_{\Lambda}\left(\eta\right)}}{\E^{(\mathrm{Can})}_{\L,\beta,\rho}\left[\ex^{-\beta \Ham_{\Lambda}}\right]}\, .
\end{equation}
Note that for this canonical model the chemical potential term and the first term in the Hamiltonian \eqref{eqn:partialHYLinteraction} are both constant, and therefore do not affect the model. In the context of the canonical ensemble we can therefore treat the partial HYL Hamiltonian as if it were just the last counter-term with parameter $b>0$.

The case of this Hamiltonian with $b=0$ is known as the Particle Mean-field Hamiltonian and is well understood. The previous work \cite{vogel2021emergence} contains a proof that random interlacements emerge in the thermodynamic limit, for example. In addition to the inclusion of the counter-term to get the HYL Hamiltonian above, we also consider a generalisation of the Particle Mean-field Hamiltonian to show that the precise choice of a quadratic function of the particle number is not important and our techniques are sufficiently resilient to be applied to other functions of the particle number.

We are going to be interested in the large-volume behaviour of these models. We will set $\Lambda=\left[-n/2,n/2\right)^d$ and aim to describe their $n\to\infty$ limits. The fundamental claim of this paper is that in this limit \textit{Brownian random interlacements} emerge from our systems of interacting finite loops. Loosely speaking these interlacements are doubly infinite continuous paths where we quotient out reparametrizations. We let $\P^\iota_u$ denote the Poisson point process on this space with an intensity that ensures they have two main properties. Firstly that the interlacements have Brownian finite dimensional distributions, and secondly that the expected local time is given by the parameter $u$. For background on \textit{Brownian} interlacements we refer the reader to \cite{sznitman2013scaling}, and to \cite{drewitz2014introduction} for a general introduction. The novelty in this paper is that we are able to show that Brownian interlacements emerge from our \emph{interacting} models, and at different intensities to those emerging from the models considered in \cite{vogel2021emergence}.

\subsection*{Organisation of the paper} In Section \ref{section results}, we give the main results and state a corollary regarding the free energy. In Section \ref{sectiondiscussion} we discuss the choice of models and assumptions, relate our work to other results and state open questions. Section \ref{SectionProofs} gives the proofs of the results from Section \ref{section results}. Given the multitude of parameters used throughout Section \ref{SectionProofs}, we give a table containing frequently used notation in the Appendix, Table \ref{table1}.
\section{Results}\label{section results}

We will now be more precise on the \emph{thermodynamic limit} we will be taking. By Gnedenko's Local Limit Theorem (see \cite[Theorem 8.4.1]{bingham1989regular}), the particle number $N_\Lambda$ under the non-interacting grand-canonical measure $\P_{\Lambda,\beta,0}$ satisfies a central limit theorem with scale $a_\Lambda$, where
\begin{equation}
    a_\Lambda = \begin{cases}
    \abs*{\Lambda}^{2/3} &\text{if }d=3,\\
    \abs*{\Lambda}^{1/2}\left(\log \abs*{\Lambda} \right)^{1/2} &\text{if }d=4,\\
    \abs*{\Lambda}^{1/2} &\text{if }d\geq 5.\\
    \end{cases}
\end{equation}
This scale will affect our result in two ways. First we will require that the scale $q_\Lambda$ appearing in the HYL Hamiltonian \eqref{eqn:partialHYLinteraction} not only satisfies $q_{\Lambda} = o\left(\abs*{\Lambda}\right)$, but also $a_{\Lambda} = o\left(q_{\Lambda}\right)$ as $\Lambda\to\R^d$. We explain the reasoning for these conditions on $q_\Lambda$ in Section~\ref{HYLDiscussion}.

The second way in which the scale $a_\Lambda$ appears in our result is in the sense in which the limit is taken. First let $n\in\N$ and set $\Lambda = \left[-n/2,n/2\right)^d$. It is then natural to interpret $a_\Lambda$ and $q_\Lambda$ as sequences $(a_n)_n$ and $(q_n)_n$. We then use $\Lambda$ to tessellate a larger box in $\R^d$. Let $(r_n)_n$ be some positive increasing sequence of real numbers that diverges to infinity at most logarithmically in $n$. Then set $C_n=[-r_n n^{d/2-1}/2,r_n n^{d/2-1}/2)^d\cap \Z^d$ and 
\begin{equation}
    \L_n=\bigcup_{x\in C_n}\left(xn+\L\right) =\left[-r_n n^{d/2}/2,r_n n^{d/2}/2\right)^d
\end{equation}
so that $\L_n$ is a finite (hyper-)cubic tessellation of boxes $\L$. We then define  $\GCanonHY{n}{\mu} := \bigotimes_{x\in C_n}\GCanonHY{xn+\L}{\mu}$ to get a measure that describes an independent superposition of loop soups distributed according to $\GCanonHY{xn+\L}{\mu}$ for each $x\in C_n$. The same is done for the canonical version to produce $\CanonHY{n}{\rho}$ from the $\CanonHY{xn+\L}{\rho}$.

Our main result is then a convergence result for the measure $\CanonHY{n}{\rho}$. Like in \cite{vogel2021emergence}, this convergence is with respect to the topology of local convergence, denoted $\xrightarrow{\mathsf{loc}}$. The definition of this topology is given precisely in Definition~\ref{DefinitionLocConv}, but can be thought of as the topology describing local, parametrization-invariant events.

We introduce some notation here. For $d\geq 3$ define the function $\brho\colon \left(-\infty,0\right]\to\left(0,\rho_{\mathrm{c}}\right]$ with $\brho(\mu) = \left(2\pi\beta\right)^{-d/2}\textsf{Li}_{d/2}(\ex^{\beta \mu})$ and $\rho_{\mathrm{c}} = \left(2\pi\beta\right)^{-d/2}\textsf{Li}_{d/2}(1)$. Here $\textsf{Li}_{s}(z)$ is the polylogarithm of order $s$ with argument $z$, and $\textsf{Li}_{d/2}(1) = \zeta\left(d/2\right)$ where $\zeta$ is the Riemann zeta function. The value $\brho(\mu)$ can be thought of as the mean density of the non-interacting grand canonical ensemble with chemical potential $\mu$ (see for example, \cite{AV19}). This is a strictly increasing convex function, and we let $\bmu$ be the inverse of $\brho$ defined on the range of $\brho$. Thus $\bmu(\rho)$ can be thought of as the chemical potential that produces a given mean density $\rho$ in the non-interacting grand canonical ensemble. In the following theorem, the notation $\otimes$ means that the resulting combined measure describes an independent superposition of the two processes described by the two separate measures.

\begin{theorem}\label{TheoremMain}
Fix $b>0$ and $d\ge 3$. Then there exists a constant $\rho_{\mathrm{c}}^\hy\le\rho_{\mathrm{c}} $ and a function $\rhob\colon\rho\mapsto \rhob(\rho)$, such that
\begin{equation}
    \CanonHY{n}{\rho}\xrightarrow{\mathsf{loc}}\begin{cases}
    \GCanonFree{\R^d}{\bmu(\rho)}&\text{ for }\rho<\rho_{\mathrm{c}}^\hy\, ,\\
    \GCanonFree{\R^d}{\bmu(\rho-\rhob)}\otimes\P_{\rhob}^\iota&\text{ for }\rho>\rho_{\mathrm{c}}^\hy\, .
    \end{cases}
\end{equation}
Furthermore
\begin{enumerate}
    \item For $d=3,4$, we have $\rho_{\mathrm{c}}^\hy<\rho_{\mathrm{c}} $ for any value of $b>0$.
    \item For $d\ge 5$, we have $\rho_{\mathrm{c}}^\hy<\rho_{\mathrm{c}} $ if $b\in\left(b_{\mathrm{c}},\infty\right)$ for $b_{\mathrm{c}}=\tfrac{1}{\brho'(0)}$. Otherwise, $\rho_{\mathrm{c}}^\hy=\rho_{\mathrm{c}} $.
    \item The function $\rhob\colon (0,\infty)\to(0,\infty)$ is continuous unless $\rho_{\mathrm{c}}^\hy<\rho_{\mathrm{c}} $, in which case it has a jump discontinuity at $\rho_{\mathrm{c}}^\hy$.
\end{enumerate}
\end{theorem}
Note that the red line in Figure \ref{figure1} gives the plot of $\rho\mapsto \rhob$. Its discontinuity is qualitatively different to the continuous transition established in \cite{vogel2021emergence}. Also, the case $\rho=\rho_{\mathrm{c}}^\hy$ is not treated in this work, we comment on that in Section \ref{sectioncrit}.

We also give a similar result for the grand-canonical case, relating it to the canonical model.
\begin{theorem}\label{TheoremGrandCanonical}
Fix $a>b>0$ and $d\ge 3$. Then there exists a function $\mu\mapsto\rho^{\mathrm{GC}}(\mu)$, mapping \textnormal{chemical potential} to \textnormal{density}, such that
\begin{equation}
     \lim_{n\to\infty}\GCanonHY{n}{\mu}= \lim_{n\to\infty}\CanonHY{n}{\rho^{\mathrm{GC}}(\mu)}\, ,
\end{equation}
in the topology of local convergence. Furthermore, $\rho^{\mathrm{GC}}(\mu)$ is monotone and satisfies $\lim_{\mu\to\infty}\rho^{\mathrm{GC}}=\infty$ and $\lim_{\mu\to-\infty}\rho^{\mathrm{GC}}=0$.
\end{theorem}
\begin{remark}
The result in Theorem \ref{TheoremGrandCanonical} can be interpreted as an equivalence of ensembles. For more discussion, we refer the reader to Section~\ref{subsectionEquiv}.

The functions $\rho^{\mathrm{GC}}$ and $\rhob$ involve polynomial term, polylogarithms (or Riemann zeta functions), their inverses, and their derivatives. Their precise form is rather lengthy, so we defer their definition to Section~\ref{SectionProofs}.
\end{remark}

The following proposition helps to illustrate why having the scale $q_\Lambda\to\infty$ is an important feature of the model. Let $\widetilde{\Ham}_\Lambda$ be the HYL Hamiltonian appearing in \eqref{eqn:partialHYLinteraction} modified so that the counter-term's sum runs over $k\geq 1$ rather than $k\geq q_\Lambda$. If we let $\E_{\Lambda,\beta,0}$ denote the expectation with respect to $\GCanonFree{\Lambda}{0}$, then we can define the thermodynamic pressures for $\Ham_\Lambda$ and $\widetilde{\Ham}_\Lambda$ as
\begin{equation}
\label{eqn:thermodynamicPressures}
    P^{\hy}\left(\beta,\mu\right) = \lim_{\Lambda \to \R^d}\frac{1}{\beta \abs*{\Lambda}}\E_{\Lambda,\beta,0}\left[\ex^{\beta\mu N_\Lambda-\beta \Ham_\Lambda}\right],\qquad \widetilde{P}^{\hy}\left(\beta,\mu\right) = \lim_{\Lambda \to \R^d}\frac{1}{\beta \abs*{\Lambda}}\E_{\Lambda,\beta,0}\left[\ex^{\beta\mu N_\Lambda-\beta \widetilde{\Ham}_\Lambda}\right].
\end{equation}
To be clear, by $\lim_{\Lambda \to \R^d}$ we mean setting $\Lambda = \left[-n/2,n/2\right)^d$ and taking $n\to\infty$. The recent results of \cite{adams2021explicit,adams2021large} prove that both these limits exist and describe them using variational expressions.
\begin{proposition}\label{propositionfullvspartial}
For all $\beta>0$ and $\mu\in\R$,
\begin{equation}
    P^{\hy}\left(\beta,\mu\right) < \widetilde{P}^{\hy}\left(\beta,\mu\right).
\end{equation}
\end{proposition}
A heuristic explanation of why the short cycles make a difference can be found in Section~\ref{HYLDiscussion}, and the detailed proof of Proposition~\ref{propositionfullvspartial} in Section~\ref{SubsectionPressureDifference}. The fact that the Hamiltonian $\Ham_\Lambda$ produces the same thermodynamic pressure as the model described by \cite{huang1957imperfect} and derived by \cite{van1990pressure}, whilst the Hamiltonian $\widetilde{\Ham}_\Lambda$ does not, suggests that the model we consider here is in a sense the `right' model and is more closely related to physical behaviour.

Next, we give a corollary to our results, which may be of independent interest.
\begin{definition}
We define the \textnormal{free energy}
\begin{equation}
   f^\hy(\beta,\rho):=-\lim_{\Lambda\to\R^d}\frac{1}{\beta\abs{\L}}\log Z^{(\mathrm{Can, HYL})}_{\Lambda,\beta,\rho} = -\lim_{\Lambda\to\R^d}\frac{1}{\beta\abs{\L}}\log\E^{(\mathrm{Can})}_{\L,\beta,\rho}\left[\ex^{-\beta \Ham_{\Lambda}}\right],
\end{equation}
where $Z^{(\mathrm{Can, HYL})}_{\Lambda,\beta,\rho}=\E^{(\mathrm{Can})}_{\L,\beta,\rho}\left[\ex^{-\beta \Ham_{\Lambda}}\right]$ is the partition function appearing in \eqref{eqn:canonHYLRadonNykodym}. 

\end{definition}
\begin{cor}\label{freeenergycorolary}
Given the conditions of Theorem \ref{TheoremMain}, we have
\begin{equation}
\begin{split}
    f^\hy(\beta,\rho)=(\rho-\rhob)\bmu(\rho-\rhob) -\int_0^{\bmu(\rho-\rhob)}\brho(s)\d s -\frac{b\rhob^2}{2}\, .
\end{split}
\end{equation}
Furthermore
\begin{equation}
\begin{split}
    P^\hy(\beta,\mu)&=\sup_{\rho>0}\left\{\mu\rho-a\rho^2/2-\beta^{-1}I(\rho)-f^\hy(\beta,\rho)\right\}\\
    &=\mu \rho^{\mathrm{GC}}-a\left(\rho^{\mathrm{GC}}\right)^2/2-\beta^{-1}I\left(\rho^{\mathrm{GC}}\right)- f^\hy\left(\beta,\rho^{\mathrm{GC}}\right)\, ,
\end{split}
\end{equation}
where the function $I$ is defined in Equation \eqref{EquationEcplicitI}.
\end{cor}

The final result is a generalisation of \cite[Theorem~2.3]{vogel2021emergence}.
Let $G\colon\left[0,\infty\right)\to \R\cup\{+\infty\}$ be measurable and bounded below. Then we define the Generalised Mean-field Hamiltonian as
\begin{equation}
    \Hgmf_\Lambda\left(\eta\right) = \abs*{\Lambda}G\left(\frac{N_\Lambda\left(\eta\right)}{\abs*{\Lambda}}\right).
\end{equation}
In turn we can define the grand-canonical measure with GMF interaction, $\P^\gmf_{\Lambda,\beta}$, by its Radon--Nikodym derivative with respect to the non-interacting measure. Specifically,
\begin{equation}
\label{eqn:GMFmeasure}
    \frac{\d \P^\gmf_{\Lambda,\beta}}{\d \GCanonFree{\Lambda}{0}}\left(\eta\right) = \frac{1}{Z^\gmf_{\Lambda,\beta}}\ex^{-\beta \Ham^\gmf_{\Lambda}\left(\eta\right)},
\end{equation}
where $Z^\gmf_{\Lambda,\beta}$ is the partition function that normalises $\P^\gmf_{\Lambda,\beta}$. Note that we are omitting the supplementary chemical potential from the GMF model. This is because we are free to replace the function $x\mapsto G\left(x\right)$ with the function $x\mapsto G\left(x\right) - \beta \mu x$ as long as the latter is bounded below.
% $x\mapsto G\left(x\right) - \beta \mu x - \inf_{y\geq 0}\left\{G\left(y\right) - \beta \mu y\right\}$ when that infimum is finite. This infimum term ensures that the function remains non-negative and doesn't affect the behaviour of the model because it is a finite constant and will get cancelled out by the partition function. 
Note that the usual Particle Mean-field measure is a special case of the Generalised Mean-field measure where $G$ is set to be a quadratic function.
% and applies to Generalised Mean-Field models (recall \eqref{eqn:GMFmeasure}). 
% That is, it concerns Hamiltonians that are solely functions of the total particle number, i.e.  $\Hgmf_\Lambda\left(\eta_\L\right) =  \abs{\L}G( N_\L/\abs{\L})$ for some measurable function $G\colon [0,\infty)\to[0,\infty]$. 
As we did for the HYL models, we define $\P^\gmf_n := \bigotimes_{x\in C_n}\P^\gmf_{xn+\L}$ to get random loop soups that are independent superpositions of loop soups distributed according to $\P^\gmf_{xn+\L}$. We first state the theorem before giving the conditions required of $G$. Let $I$ be the rate function for the particle number with no interaction, given in Equation \eqref{EquationEcplicitI}. Note that this can be also written as the sum of the free energy and the pressure of the system (from a simpler application of the techniques in \cite{adams2021explicit}).
\begin{theorem}
\label{thm:GMFlimit}
Let $d\geq 3$ and suppose that $G$ satisfies Assumption \ref{conditionG}. If $I+G$ attains its unique minimum at $\rho>0$, then as $n\to\infty$
\begin{equation}
    \P^\gmf_n\xrightarrow{\mathrm{loc}}\begin{cases}
    \P_{\R^d,\beta,\bmu(\rho)}&\text{ for }\rho<\rho_{\mathrm{c}}\, ,\\
    \P_{\R^d,\beta,0}\otimes\P_{\rho-\rho_{\mathrm{c}}}^\iota&\text{ for }\rho>\rho_{\mathrm{c}}\, .
    \end{cases}\,
\end{equation}
\end{theorem}
For a discussion on the case of multiple minimizers, see Remark \ref{remarkmin}.

We now give the (not very strict) conditions on $G$.
\begin{assump}\label{conditionG}
We assume that $I+G$ has a unique minimizer $x_\min$ with $x_\min>0$. Set $G(x_\min)=K$. Furthermore, we require that one of the following holds:
\begin{enumerate}
    \item $x_\min<\rho_{\mathrm{c}}$, and $G$ is twice differentiable at $x_\min$ and continuous in a neighbourhood of $x_\min$. We furthermore require that for any $\e>0$, there exists a $\delta>0$ such that $G^{-1}\left[[K,K+\delta)\right]\subset (x_\min-\e,x_\min+\e)$.
    \item $x_\min>\rho_{\mathrm{c}}$, and for any $\e>0$, there exists a $\delta>0$ such that $G^{-1}\left[[K,K+\delta)\right]\subset(x_\min-\e,x_\min+\e)$. We also require that
    \begin{enumerate}
        \item either $G$ is twice differentiable at $x_\min$ at has a local minimum there,
        \item or $G$ has a jump-discontinuity from the left (or the right), is differentiable in a right (resp. left) neighbourhood of $x_\min$ and the first derivative\footnote{at $x_\min$, we take the right (resp. left) derivative} is bounded uniformly from below.
    \end{enumerate}
\end{enumerate}
\end{assump}

\begin{example}
There are various models that satisfy the conditions of Assumption \ref{conditionG}. For example:
\begin{enumerate}
    \item \label{Example:PMF} $G\left(x\right) = \tfrac{a}{2}x^2+\infty\1\{x>\rho\}$ for $a>0$ and $\rho>\rho_{\mathrm{c}}$ was considered in \cite[Theorem~2.3]{vogel2021emergence}. Therefore Theorem~\ref{thm:GMFlimit} is indeed a generalisation of \cite{vogel2021emergence}.
    \item $G$ is strictly convex, $G(x)\to+\infty$ as $x\to\infty$, and twice differentiable at $x_\min$. 
    \item $G''(x)\geq 0$ for all $x\geq 0$, and $G'(\rho_{\mathrm c}) >0$.
    % \item \label{Example:NonConvex} $G(x) = x^3 - 3x^2 + 2x$ (or any non-convex, increasing function with strictly positive derivative).
\end{enumerate}
% In examples (2)--(3), the differentiability criterion in Assumption~\ref{conditionG} must be checked separately. A
In general, as we have a rather explicit representation (see Equation \eqref{EquationEcplicitI}) of $I$ in terms of polylogarithms, checking Assumption \ref{conditionG} for a specific $G$ is usually rather straight-forward.
\end{example}

%
%%%%
%
\section{Discussion}\label{sectiondiscussion}
\subsection{Momentum HYL and full HYL}\label{HYLDiscussion}
In \cite{huang1957imperfect}, Huang, Yang, and Luttinger considered a gas of bosons experiencing a hard-sphere interaction. To study the virial coefficients, they expanded the thermodynamic pressure of the interacting gas to second order in the dimensionless parameter $a/\lambda$. Here $a$ is the diameter of the hard-sphere interaction and $\lambda$ is the so-called 'thermal wavelength' --- a length scale corresponding to the de Broglie wavelength of a massive particle with energy $1/\beta$. Of particular interest in this was the first order perturbation, which they described using the single-particle momentum eigenstates of the non-interacting model. In fact, they were able to express the first order perturbation solely in terms of the expectation of these eigenstates' occupation numbers. Inspired by this, Huang, Yang, and Luttinger `invented' a fictitious Bose gas whose energy levels were given by taking the first order perturbation and replacing the expectations with the raw variables. If we let $\alpha$ label the countably many single-particle momentum eigenstates and $n_\alpha$ denote their occupation numbers, this `momentum HYL' interaction energy can be written as
\begin{equation}
\label{eqn:momentumHYLenergy}
    H^{\left(\mathrm{mHYL}\right)}_\L\left(n\right) = \frac{a}{2\abs{\L}}\left(\left(\sum_{\alpha\geq 1}n_\alpha\right)^2-\frac{1}{2}\sum_{\alpha\geq 1} n^2_\alpha\right).
\end{equation}

As described in \cite{huang1957quantum}, the first 'square of the sum' term can be expected classically --- on the basis of an ``index of refraction approximation" --- whereas the second `sum of squares' term is purely quantum mechanical. The Heisenburg uncertainty principle applied to the relative distance of two particles and their relative momentum indicates that particles prefer to be in the same momentum state in order to minimise the spatial repulsion from the hard spheres.

In fact, when discussing their fictitious model in \cite{huang1957imperfect}, Huang, Yang, and Luttinger replaced $\sum_\alpha n^2_\alpha$ with $n^2_0$, omitting all terms other than that arising from the single-particle ground state. At the time, they justified this simplification by noting that near condensation the average occupation number $n_0$ would be much higher than that of the other `excited' states. However, it was not until 1990 that van den Berg, Dorlas, Lewis, and Pul{\'e} proved in \cite{van1990pressure} that the thermodynamic pressures given by the full momentum HYL energy and the ground-state-only version truly are equal. They did this by using large deviation techniques --- applying Varadhan's Lemma with two different topologies allowed them to tightly bound the pressure from above and below. However, the use of different topologies meant that this did not prove a large deviation principle for the model.

At first glance, one may mistakenly expect the loop Hamiltonian we consider in this paper to have only superficial similarity to the original momentum HYL model. Not only does it replace momentum eigenstate occupation numbers with loop type occupation numbers, but the second term omits more and more types of `short' loop as the thermodynamic limit is taken. However the similarity in fact runs much deeper:
\begin{center}
It has been proven in \cite[Theorem~2.2]{adams2021explicit} that the thermodynamic pressure of this loop model is \emph{precisely} equal to the thermodynamic pressure of the HYL Bose model derived in \cite{van1990pressure}.
\end{center}
In Proposition~\ref{propositionfullvspartial} below we show that omitting the short loops is crucial: including the short loops in the Hamiltonian destroys this equality.

It is a conjecture of Feynman that the emergence of ``long loops" in bosonic loop soup models should correspond to Bose--Einstein condensation, and that the fraction of `particles' on such loops will equal the fraction in the condensate (see for example \cite{feynman1953atomic,feynman1972statistical}). The Heisenberg uncertainly principle also suggests that low momenta states should relate to long cycles. This relation motivates why we want the Hamiltonian to keep influencing long loops: \cite{van1990pressure} showed that it was only the effect on the `condensate state' that mattered, and for the loop model this corresponds to our `long loops.'

It is important that the loop model Hamiltonian does not have the `sum of squares' term include all types of cycle. That is, we do not have the parameter $q_\Lambda$ in \eqref{eqn:partialHYLinteraction} remain bounded as $\Lambda\to\R^d$. For the momentum model there are countably many discrete states and in the thermodynamic limit these get closer and closer and approach a continuum of states. The occupation density of each of the individual discrete \emph{excited} states vanishes in the thermodynamic limit and so the energy contribution to the `sum of squares' term from these is not significant. For the loop model, the states stay separated and each maintains a positive density in the thermodynamic limit. Hence the thermodynamic pressure of the full cycle HYL model differs from that of the partial cycle HYL model. This reasoning is made rigorous in the proof of Proposition~\ref{propositionfullvspartial} contained in Section~\ref{SubsectionPressureDifference}.

\begin{figure}[ht]
    \centering
    \begin{subfigure}[t]{0.45\textwidth}
    \centering
    \begin{tikzpicture}[scale=1.2]
    \draw[->] (0,0) -- (4,0) node[below]{Total density $\rho$};
    \draw[->] (0,0) -- (0,4)node[above]{Limiting density};
    \draw[loosely dashed] (0,0) -- (4,4);
    \draw[very thick, blue, dotted] (0,0) -- (1.5,1.5);
    \draw[very thick, blue, dotted] (1.5,1) to [out=290,in=180] (4,0.1);
    \draw[very thick, red] (0,0) -- (1.5,0);
    \draw[very thick, red] (1.5,0.5) to [out=70,in=225] (4,3.9);
    \draw[loosely dashed] (1.5,0) node[below]{$\rho^\hy_c$} -- (1.5,4);
    \draw[loosely dashed] (2.5,0) node[below]{$\rho_{\mathrm{c}}$} -- (2.5,4);
    \end{tikzpicture}
    \caption{Partial HYL model for $d=3,4$, and $d\geq5$ with certain parameters.}
    \end{subfigure}
    \hfill
    \begin{subfigure}[t]{0.45\textwidth}
    \centering
    \begin{tikzpicture}[scale=1.2]
    \draw[->] (0,0) -- (4,0) node[below]{Total density $\rho$};
    \draw[->] (0,0) -- (0,4)node[above]{Limiting density};
    \draw[loosely dashed] (0,0) -- (4,4);
    \draw[very thick, blue, dotted] (0,0) -- (2,2);
    \draw[very thick, blue, dotted] (2,2) to [out=280,in=180] (4,0.1);
    \draw[very thick, red] (0,0) -- (2,0);
    \draw[very thick, red] (2,0) to [out=80,in=225] (4,3.9);
    \draw[loosely dashed] (2,0) node[below]{$\rho_{\mathrm{c}}$} -- (2,4);
    \end{tikzpicture}
    \caption{Partial HYL model for $d\geq5$ with certain parameters.}
    \end{subfigure}
    \caption{The density of the random interlacements/condensate (\protect\redline) and the finite loops/bulk (\protect\bluedottedline) for the partial HYL model in different dimensions and different parameters.}
    \label{figure2}
\end{figure}
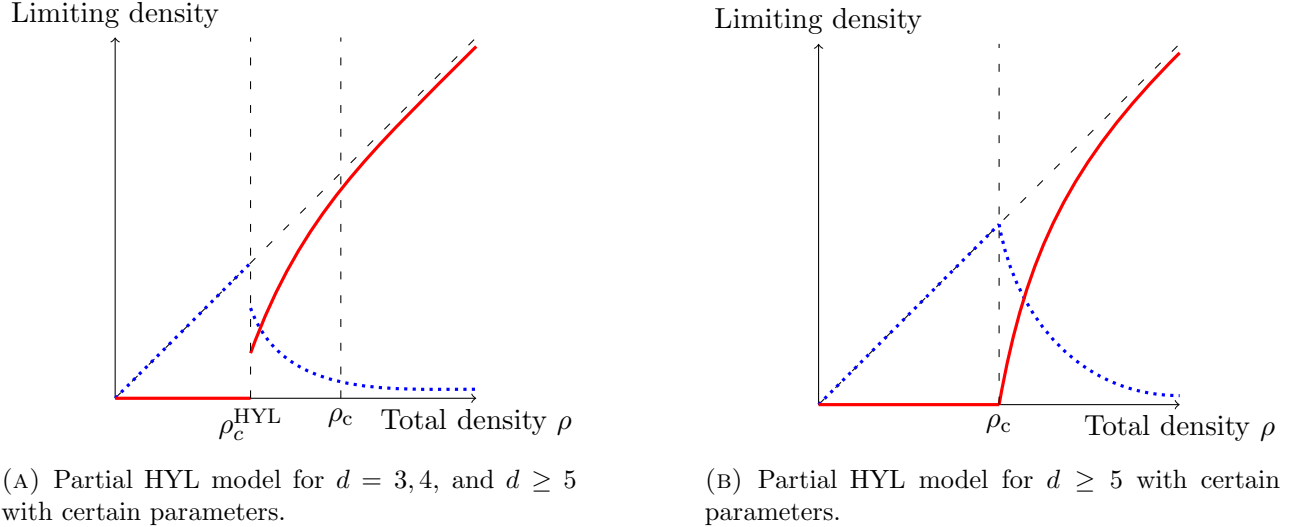

\subsection{Equivalence of Ensembles}\label{subsectionEquiv}
The question of equivalence of ensembles can be viewed from different standpoints. On one level, we can ask if there is a relation between various thermodynamic functions of the ensembles. For example, \cite{ruelle1969statistical} proves that the pressure (associated with the grand-canonical ensemble), the free energy (associated with the canonical ensemble), and the entropy (associated with the microcanonical ensemble) in the thermodynamic limit for various interacting particle models can be related to each other by Legendre--Fenchel transforms over appropriate parameter spaces. An alternative approach is to study the measures more directly. Under an appropriate topology, is it possible to relate the accumulation points of the finite-volume measures in the various ensembles? In \cite{georgii1995equivalence} it is proven that such an equivalence at the level of measures for microcanonical and grand-canonical ensembles of particle models (not bosons) with suitable pair-wise interactions is indeed possible.

Theorem \ref{TheoremGrandCanonical} implies an equivalence of ensembles on the level of measures. It furthermore shows that the limiting measures are universal, in a limited sense: not only do the HYL models converge to the superposition of the free Bose gas and the random interlacements, but also the mean-field and the free Bose gas do as well (with different parameters, of course). This can be seen as a confirmation of Feynman's prediction that
\begin{center}
``\textit{...the strong interactions between particles do not prevent these particles from behaving very much as though they move freely among each other.}",
\end{center} see \cite{feynman1953atomic}. We predict that similar results will hold for more involved interactions between particles.

\subsection{Condensate discontinuity}\label{subsectioncondensatediscont}
The study of the grand-canonical ensemble for the partial loop HYL model in \cite{adams2021explicit} derived a large deviation principle --- a lower resolution result than that derived here. Nevertheless, it was sufficient to derive the `condensate' density of that loop soup model. Much like \cite{huang1957imperfect,lewis1986bosons} did for the momentum HYL Bose gas, it was shown that for certain parameters a discontinuity in the condensate density can occur as the chemical potential $\mu$ crosses the transition point. Here we show that a similar discontinuity in the density of the random interlacement occurs for the canonical ensemble as the density crosses the new critical density. Figure~\ref{figure2} shows how this can occur.

\subsection{On the choice of the intermediate scale}\label{subsectionchoicescale}
In the statement of our results, we gave two requirements the sequence $(q_n)_n$ has to satisfy, in order for our results to remain valid. We justify our choice as follows:

The requirement that $q_n=o\left(\abs{\L}\right)$ is due to the ``unphysical" results larger $q_n$'s would give. Indeed, no changes to the proof are necessary to treat this case. However, a choice of $q_n\ge \e\abs{\L}$ results in a built-in prevention of the existence of interlacements with densities lower than $\e>0$. Therefore, our requirements on $(q_n)_n$ ensure that macroscopic (and mesoscopic) loops are not being artificially excluded. A study of the large deviations of a system that does exclude such loops and the consequences for the condensate  in such a system can be found in \cite{adams2021explicit}.

Requiring $a_n=o(q_n)$ can be seen as the bigger restriction of our model. The condition can be motivated by the fact that $q_n$'s which grow slower lead to interaction between the small loops, encouraging a clumping. By clumping, we mean that loops share the same lengths more often than we would expect in the free case. Applying the Heisenberg uncertainty principle to the conjugate pair of position and momentum suggests that particles residing on short loops (and therefore having a smaller length scale) would somehow be related to particles having a higher momentum. Since these higher momentum states are more sparsely occupied, one may guess that the contribution from this clumping of short loops would be small --- perhaps negligible. It may be that at the large deviation scale it indeed doesn't matter, but that it does for our higher-resolution study here. Certainly, as it can be seen from Proposition \ref{propositionfullvspartial}, the model for fixing $q_n\equiv 0$ is different to both the models studied here and in \cite{adams2021explicit}. We plan to address the physicality of the full HYL (loop) model in a future publication.
\subsection{The critical case}\label{sectioncrit}
We do not prove anything about the case $\rho=\rho_{\mathrm{c}}$ (resp. $\rho=\rho_{\mathrm{c}}^\hy$). It is standard to require for the canonical ensemble that $N_\L/\abs{\L}$ converges to $\rho$. However, depending on the sequence we choose, different global phenomena emerge: if we choose $N_\L=\rho\abs{\L}+\abs{\L}^{\alpha}$ with $\alpha\in(0,1)$ \textit{greater than} the CLT coefficient, macroscopic loops form (as the estimates from \cite{berger2019notes} are valid in that regime). The density of these macroscopic loops is vanishing, so that they cannot be detected from a local perspective. For $N_\L=\rho\abs{\L}+b_\L$ with $b_\L=\Ocal\left(a_\L\right)$, we do not expect any formation of infinite loops. Instead, there is tilting in the distribution of $N_\L$. We believe that these phenomena warrant an independent investigation.

\subsection{Interlacements in low dimensions}\label{sectionLowDimensions}
In studying their versions of an HYL interaction, both \cite{lewis1986bosons} and \cite{adams2021explicit} found that their versions of condensate behaviour occurred in every dimension $d\geq 1$. Contrast this with the non-interacting models, in which no condensation occurs for $d=1,2$. Naturally then, \cite{armendariz2019gaussian,vogel2021emergence} studied the emergence of interlacements for $d\geq 3$ only. Whilst interlacements in $d=2$ can make sense, their construction is quite different (see \cite{comets2016two}), and for $d=1$ there is currently no framework for them. For these reasons we have restricted our attention in this paper to $d\geq 3$. Nevertheless the question of whether the emergence of `long loops' for the partial loop HYL model can be understood via random interlacements is an interesting avenue of future study.

%%%%%%%%%%%%%%%%%%%
%%%%%%%%%%%%%%%%%%%%
%%%%%%%%%%%%%%%%%%%%
\section{Proofs}\label{SectionProofs}
As there are multiple phenomena contributing the formation of infinite loops, we split the proof in different sections. Overall, the structure is as follows:
\begin{enumerate}
    \item In Section \ref{SubsectionNotation}, we introduce further notation.
    \item In Section \ref{SubsectionAnalysisRF} we perform a careful analysis of the large deviation rate function for the free Bose gas. The most important result is the analysis of its behaviour near the critical point $\rho_{\mathrm{c}}$.
    \item In Section \ref{SubsectionCalculationPartitionFUnction}, the partition function is calculated up to order $(1+o(1))$. We identify the two sources making up the density of the random interlacements: the contribution from the free loops and the one induced by the Hamiltonian. This presents a substantial difference to the mean-field model considered in \cite{vogel2021emergence}.
    \item The computation of the limiting measure is performed in Section \ref{subsectionlimiting}. Given the previous section, this follows after some approximation arguments, using a quicker strategy compared to \cite{vogel2021emergence}.
    \item Section \ref{subsectionbelowctricial} does the analysis for the case $\rho<\rho_{\mathrm{c}}$, which was previously excluded. We also prove the discontinuity of the density of infinite loops as $\rho$ varies.
    \item The grand-canonical case is solved in Section \ref{subsectionGC}. We use the results from the canonical case together with a large-deviation principle for distribution of particle number under mean-field interaction.
    \item The generalisation of the mean-field results in \cite{vogel2021emergence} is given in Section \ref{SubsectionGMF}.
\end{enumerate}
Recall that a table containing frequently used notation is given in the Appendix, Table \ref{table1}.
\subsection{Further notation}\label{SubsectionNotation}
Let $ \Gamma_\mathrm{F}$ be the space of \textit{finite} loops, i.e.
\begin{equation}
     \Gamma_\mathrm{F}=\bigcup_{j\ge 1}\big\{\omega\colon [0,\beta j]\to\R^d,\, \omega(0)=\omega(\beta j)\text{ and }\omega\text{ continuous}\big\}\, .
\end{equation}
We also set the space of \textit{random interlacements}
\begin{equation}
    \Gamma_\mathrm{I}=\Big\{\omega\colon(-\infty,\infty)\to\R^d\colon \lim_{\abs{t}\to\infty}\abs{\omega(t)}=+\infty\text{ and }\omega\text{ continuous}\Big\}\, .
\end{equation}
Denote $\Gamma=\Gamma_\mathrm{F}\cup\Gamma_\mathrm{I}$. Given $t\in \R$, we define the shift $\theta_t$ as follows:
\begin{enumerate}
    \item $\omega\circ\theta_t(s)=\omega(t+s)$, if $\omega\in \Gamma_\mathrm{I}$.
    \item $\omega\circ\theta_t(s)=\omega(t+s \!\!\mod\beta j)$, if $\Gamma_\mathrm{F}\ni\omega\colon [0,\beta j]\to\R^d$.
\end{enumerate}
We define an equivalence on $\Gamma$ as follows: $\omega_1$ is equivalent to $\omega_2$ if there exists a $t\in \R$ such that $\omega_1=\omega_2\circ \theta_t$. Let $\Gamma^*$ (and $\Gamma_\mathrm{F}^*,\Gamma_\mathrm{I}^*$) be the space of equivalence classes on $\Gamma$ (resp. $\Gamma_\mathrm{F},\Gamma_\mathrm{I}$). Let $\Pi$ denote the projection from $\Gamma$ to $\Gamma^*$ and let $\amalg$ be the preimage of $\Pi$, i.e. $\amalg(A)=\Pi^{-1}[A]$ for any set $A\in \Gamma^*$. For $u>0$, let $\nu_u$ be the intensity measure of the Brownian random interlacements on $\Gamma^*$, as defined in \cite{sznitman2013scaling}.

Let $\p_t(x,y)=\p_t(x-y)$ be the transition kernel of a standard Brownian motion in $\R^d$, for $x$ and $y$ two points in $\R^d$ and $t>0$. For $x\in\R^d$, we set $\B_{x,x}^t$ the measure of a standard Brownian bridge, conditioned to return to $x$ at time $t>0$. We set
\begin{equation}\label{EqDefBridgeMe}
    \P_{x,x}^t=\p_{t}(x,x)\B_{x,x}^t\, .
\end{equation}
For $\L\subset\R^d$, we write $M_\L$ for the loop measure:
\begin{equation}
    M_\L=M_{\L,\beta,\mu}=\int_{\L}\d x\sum_{j\ge 1}\frac{\ex^{\beta \mu j}}{j}\P_{x,x}^{\beta j}\, .
\end{equation}
Let $\P_{\L}=\P_{\L,\beta,\mu}$ be the Poisson point process (PPP) with intensity measure $M_{\L,\beta,\mu}$. A sample of $\P_\L$ will be denoted by $\eta$ and can be written as
\begin{equation}
    \eta=\sum_{k}\delta_{\omega_k}\, ,
\end{equation}
with $\omega_k\in\Gamma_\mathrm{F}$. We write $\omega\in\eta$ whenever $\omega\in\supp(\eta)$. For $\omega\in\Gamma_\mathrm{F}$, we set $\ell(\omega)=j$, if $\omega\colon[0,j\beta]\to\R^d$. For $\Delta\subset\R^d$, we set
\begin{equation}
    \eta_\Delta=\sum_{\omega\in\eta}\delta_{\omega}\1\{\omega(0)\in\Delta\}\, ,
\end{equation}
and
\begin{equation}
    N_\Delta(\eta)=\sum_{\omega\in\eta_\Delta}\ell(\omega)\, .
\end{equation}
We also set $\bN_\Delta(\eta)=\abs{\Delta}^{-1}N_\Delta(\eta)$. Finally, set $N_o=N_{[0,1]^d}$. 
\begin{definition}\label{DefinitionLocConv}
The topology of local convergence is generated by functions $F$ of the type $F(\eta)=\ex^{-\eta[f]}$ for $f\colon \Gamma\to [0,\infty)$ which satisfy the following properties
\begin{enumerate}
    \item $f(\omega)=f(\omega\circ\theta_t)$, for any $t\in\R$,
    \item $f$ depends only on the values of $\omega$ on some compact set.
    \item $f$ is continuous in the Skorokhod topology on $\Gamma$.
\end{enumerate}
\end{definition}
It is known (see \cite[Theorem 24.7]{klenke2013probability}) that such $F$ generate the topology of continuous bounded functions $\eta\mapsto F(\eta)$ which are invariant under reparametrization of the loops, local and invariant under permutation of the loops. For more details (for example of the Skorokhod topology), see also \cite[Section 3]{vogel2021emergence}.

We introduce some functions, relating the macroscopic behaviour of the gas: for $x\le 0$, set the pressure
\begin{equation}
    P(x)=\sum_{j\ge 1}\frac{\ex^{\beta x j}}{j}\p_{\beta j}(0)\, ,
\end{equation}
and the density
\begin{equation}
    \brho(x)=\sum_{j\ge 1}{\ex^{\beta x j}}\p_{\beta j}(0)=\frac{1}{\beta}\frac{\d}{\d x}P(x)\, .
\end{equation}
Since for $d\geq 1$,
\begin{equation}
    \p_{\beta j}(0)= \frac{\c_d}{(\beta j)^{d/2}}, \qquad \text{where }\c_d=\frac{1}{(2\pi)^{d/2}},
\end{equation}
the $n$-th left derivative of $\brho$ at the origin exists only for $d>2n+2$.

\begin{figure}[ht]
    \centering
    \begin{subfigure}[t]{\textwidth}
    \centering
    \begin{tikzpicture}[scale=1.2]
        \draw[->] (4,-2) -- (4,2)node[above]{$\brho(\mu)$};
        \draw[->] (0,-1.5) -- (4.2,-1.5)node[below]{$\mu$};
        \draw[thick] (0,-1.45) to [out=0,in=268] (4,1.5);
        \filldraw[black] (4,1.5) circle (1pt) node[right]{$\rho_{\mathrm{c}}$};
        
        \draw[->] (6,-2) -- (6,2)node[above]{$\bmu(\rho)$};
        \draw[->] (6,1.5) -- (10.2,1.5)node[below]{$\rho$};
        %\draw (9,1.5)node[below]{$\rho_{\mathrm{c}}$};
         \draw[thick] (6.05,-2) to [out=90,in=180] (9,1.5);
         \draw[thick]  (9,1.5) to (10.2,1.5);
        \filldraw[black] (9,1.5) circle (1pt) node[below]{$\rho_{\mathrm{c}}$};
    \end{tikzpicture}
    \caption{Sketches of $\brho$ and $\bmu$, for $d\ge 3$}
    \end{subfigure}
    \hfill
    \begin{subfigure}[t]{\textwidth}
    \centering
    \begin{tikzpicture}[scale=1.2]
        \draw[->] (4,-2) -- (4,2)node[above]{$\brho'(\mu)$};
        \draw[->] (0,-1.5) -- (4.2,-1.5)node[below]{$\mu$};
        \draw[thick] (0,-1.45) to [out=0,in=268] (4,0.5);
        \filldraw[black] (4,0.5) circle (1pt) node[right]{$\brho'(0)$};
        \draw[very thick, blue, dotted] (0,-1.45) to [out=0,in=270] (3.95,2);
        
        \draw[->] (6,-2) -- (6,2)node[above]{$\bmu'(\rho)$};
        \draw[->] (6,-1.5) -- (10.2,-1.5)node[below]{$\rho$};
        %\draw (9,1.5)node[below]{$\rho_{\mathrm{c}}$};
         \draw[very thick, blue, dotted] (6.05,2) to [out=270,in=180] (9,-1.5);
         \draw[thick] (6.05,2) to [out=270,in=180] (9,-1);
        \filldraw[black] (9,-1.5) circle (1pt) node[below]{$\rho_{\mathrm{c}}$};
        \draw[dashed] (9,-1) -- (9,-1.5);
        \draw[thick] (9,-1.5) -- (10.2,-1.5);
        \filldraw[black] (9,-1) circle (1pt) node[above right]{\small{$\bmu'(0)=\frac{\beta}{\brho'(0)}$}};
    \end{tikzpicture}
    \caption{Sketches of $\brho'$ and $\bmu'$, for $d\ge 3$}
    \end{subfigure}
    \caption{Sketch of the behaviour the thermodynamic functions relation density and chemical potential. The behaviour for $d\geq5$ is drawn as a solid line. For the first derivatives, the behaviour is qualitatively different for $d=3,4$ and is drawn as a blue dotted line.}
    \label{figurethermodynamic functions}
\end{figure}
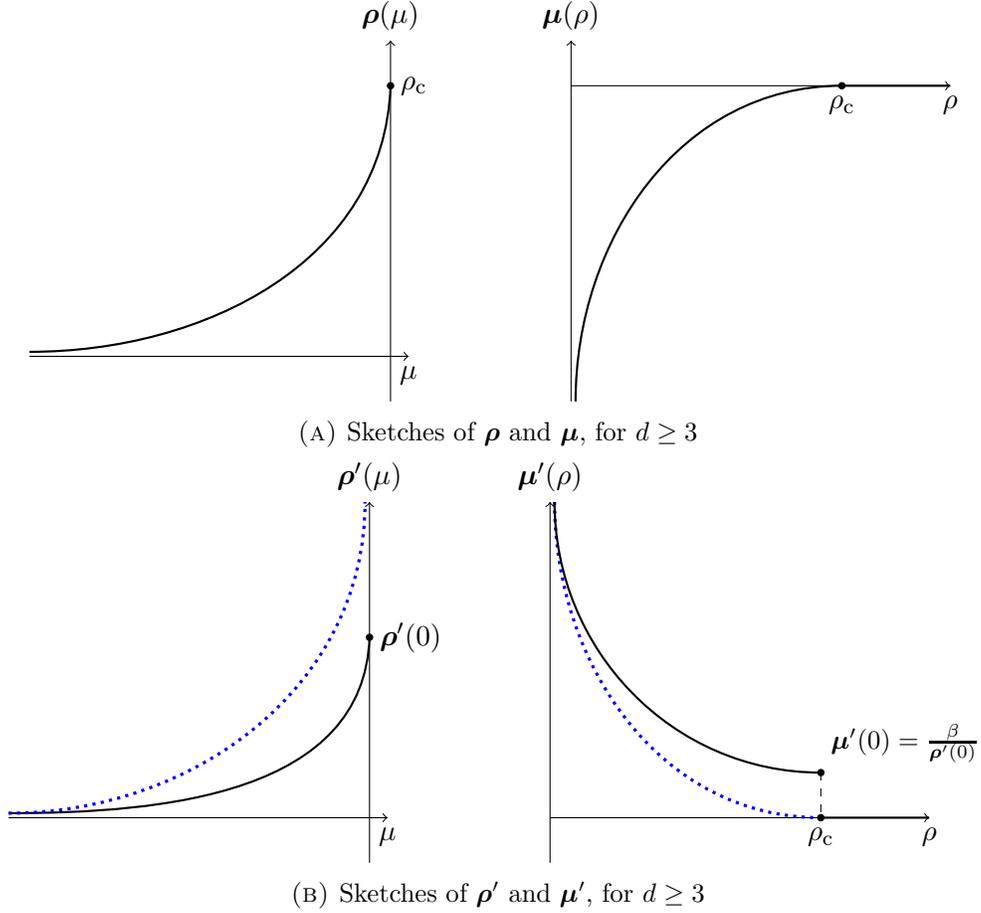

For $x\in (0,\brho(0)]$, we let $\bmu(x)$ be the unique number such that $\brho(\bmu(x))=x$. For larger $x$, we extend $\bmu$ by setting it to zero. As $\rho$ and $\mu$ are also parameters of the model, we have chosen to use the \textbf{boldsymbol}, to stress the difference. See Figure \ref{figurethermodynamic functions} for a sketch of the functions.

%  From now on, we incorporate the $\beta$-scaling $\Ham$ and $\Hpmf$ into the constants $a,b>0$, so that $a$ (and $b$) are understood as $a\beta$ (resp. $\beta b$). 
 For two sequences $(a_n)_n,(b_n)_n$, we write $a_n\sim b_n$ whenever as $a_n=b_n(1+o(1))$, as $n\to\infty$. We also set $\bB_\e(x)$ as the open ball of radius $\e>0$ around $x\in\R^d$.
\subsection{Analysis of the free rate function}\label{SubsectionAnalysisRF}
Let $\phi(t)$ be the log moment generating function of $N_o$ with respect to $\P_{\L,\beta,\mu}$. We then have that,
\begin{equation}\label{MGF}
    \phi(t)=\begin{cases}
    P(\mu+t/\beta)-P(\mu)&\text{ if }t\le -\mu\, ,\\
    +\infty &\text{ otherwise}    \, .\end{cases}
\end{equation}
Indeed, by the Campbell formula, for $t\le \mu$
\begin{equation}
    \E_{\L,\beta,\mu}\left[\ex^{tN_o}\right]=\exp\left(M_{\L,\beta,\mu}[\ex^{tN_o}-1]\right)=\exp\left(\sum_{j\ge 1}\frac{\ex^{\beta\mu j}}{j}\E_{0,0}^{\beta j}[\ex^{tj}-1]\right)=\exp\left(P(\mu+t/\beta)-P(\mu)\right),
\end{equation}
where $\E_{0,0}^{\beta j}$ is the expectation with respect to the unnormalised Brownian bridge measure $\P_{0,0}^{\beta j}$ defined in Equation \eqref{EqDefBridgeMe}. We use this in the next result.
\begin{lemma}\label{lemmafreeratefimctopm}
The large deviation rate function $I_\mu(x)$ associated to $\bN_\L$ is given by
\begin{equation}\label{EquationEcplicitI}
   I_\mu(x)=I_{\beta,\mu}(x)=\begin{cases}
   +\infty &\text{ if }x<0\, ,\\
   P(\mu)&\text{ if }x=0\, ,\\
   \beta x\left(\bmu(x)-\mu\right)-P\left(\bmu(x)\right)+P(\mu)&\text{ if }0<x\le \rho_{\mathrm{c}}\, ,\\
  -x\beta \mu -P(0)+P(\mu)&\text{ otherwise}
   \, . 
   \end{cases}
\end{equation}
For $\mu=0$, the rate function is not good. We write $I(x)$ instead of $I_0(x)$.
\end{lemma}
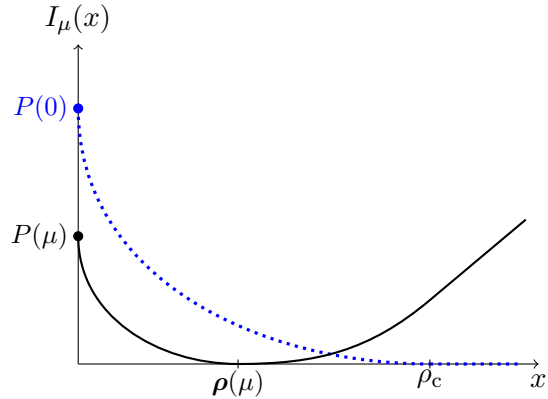
\begin{figure}[ht]
    \centering
% \begin{tikzpicture}[scale=0.85]
%     \draw[->] (0,0) -- (8.2,0) node[below]{$x$};
%     \draw[->] (1,0) -- (1,5)node[above]{$I_\mu(x)$};
%     %\draw(0.5,4.5) node[above]{\small{$+\infty$}} ;
%     % \draw(0,0) circle (2pt)
%     \draw[very thick, black] (0,4.5) -- (1,4.5);
%     \filldraw[blue] (1,4) circle (1pt) node[left]{\small{$P(0)$}};
%     \draw[very thick, blue, dotted] (1,4) to [out=290,in=175] (6.5,0);
%     \draw[very thick, blue, dotted] (6.5,0) to (7.9,0);
%       \filldraw[red] (1,3) circle (1pt) node[left]{\small{$P(\mu)$}};
%     \draw[ thick, red] (1,3) to [out=290,in=175] (5,0);
%     \draw[ thick, red] (5,0) to (8,1);
%     % \draw[very thick, red] (0,0) -- (1.5,0);
%     % \draw[very thick, red] (1.5,0.5) to [out=70,in=225] (4,3.9);
%     % \draw[loosely dashed] (1.5,0) node[below]{$\rho^\hy_c$} -- (1.5,4);
%     \draw[-] (6.5,-0.07) -- (6.5,0.07) node[below]{$\rho_{\mathrm{c}}$} ;
%     \draw[-] (5,-0.07) -- (5,0.07) node[below]{\small{$\brho(\mu)$}} ;
%     \end{tikzpicture}
\begin{tikzpicture}[scale=0.85]
    \draw[->] (1,0) -- (8.2,0) node[below]{$x$};
    \draw[->] (1,0) -- (1,5)node[above]{$I_\mu(x)$};
    \filldraw[blue] (1,4) circle (2pt) node[left]{\small{$P(0)$}};
    \draw[very thick, blue, dotted] (1,4) to [out=270,in=180] (6.5,0);
    \draw[very thick, blue, dotted] (6.5,0) to (7.9,0);
    \filldraw[] (1,2) circle (2pt) node[left]{\small{$P(\mu)$}};
    \draw[thick] (1,2) to [out=270,in=180] (3.5,0) to [out=0,in=220] (6.5,1) -- (8,2.26) ;
    \draw[-] (6.5,-0.07) -- (6.5,0.07) node[below]{$\rho_{\mathrm{c}}$} ;
    \draw[-] (3.5,-0.07) -- (3.5,0.07) node[below]{\small{$\brho(\mu)$}} ;
    \end{tikzpicture}
    \caption{Two sketches of $I$, once for $\mu<0$ (black, solid line) and once for $\mu=0$ ({\color{blue}blue}, dotted).}
    \label{LDP function}
\end{figure}

\begin{proofsect}{\textbf{Proof of Lemma \ref{lemmafreeratefimctopm}}}
The claim for $x<0$ is immediate. For $x\ge 0$, the result follows from Cram\'{e}r's theorem. For $x=0$, note that
\begin{equation}
    I(0)=\sup_{t\in\R}\{-\phi(t)\}\, .
\end{equation}
As $P(\mu+t/\beta)$ is monotonous, positive and converges to zero as $t\to \infty$, the claim follows for $x=0$. For $x\in (0,\rho_{\mathrm{c}})$, notice that the equation
\begin{equation}
    x=\phi'(t)\quad\Longleftrightarrow\quad x=\brho(\mu+t/\beta)\, ,
\end{equation}
has the unique solution $\tilde t=[\bmu(x)-\mu]\beta$. We then use that $I(x)=x\tilde t-\phi(\tilde t)$. For $x\ge \rho_{\mathrm{c}}$, notice that due to Equation \eqref{MGF},
\begin{equation}
    \R\ni t\mapsto xt-\phi(t)\, ,
\end{equation}
is maximised at $t=-\beta \mu$. This proves the claim for $x\ge \rho_{\mathrm{c}}$. For $\mu=0$, $I_0(x)$ is constant for $x\ge \rho_{\mathrm{c}}$ and hence does not have compact level-sets (and is therefore not good).\qed
\end{proofsect}

\begin{lemma}\label{LemDerivRat}
For every $\mu<0$, $\brho(\mu)$ is smooth at $\mu$, same for $\bmu(\rho)$ for $\rho<\rho_{\mathrm{c}}$. Furthermore, for every $\e\in (0,1/2)$ as $\mu\uparrow 0$
\begin{equation}
  \brho(\mu)=\begin{cases}
    \rho_{\mathrm{c}}+\mu\brho'(0)+\Ocal\left(\mu^{3/2}\right)&\text{ if }d\ge 5\, ,\\
    \rho_{\mathrm{c}}+ \c_4\mu\beta^{-1}\log(-\mu^{-1})\left(1+o(1)\right)&\text{ if }d=4\, ,\\
    \rho_{\mathrm{c}}-(-2\mu)^{1/2}(\pi\beta)^{-1}\left(1+o(1)\right)&\text{ if }d=3\, ,
    \end{cases}
\end{equation}
which for $d=4$ means $\brho(\mu)=\rho_{\mathrm{c}}+ \mu^{1-\e}o(1)$.\\
Let $W_{-1}$ be the $-1$ branch of the Lambert W function, see \cite{corless1996lambertw} for a definition. As $\rho\uparrow\rho_{\mathrm{c}}$
\begin{equation}
 \bmu(\rho)=\begin{cases}
    (\rho-\rho_{\mathrm{c}}) \brho'(0)^{-1}+\Ocal([\rho_{\mathrm{c}}-\rho]^{3/2})&\text{ if }d\ge 5\, ,\\
    (\rho-\rho_{\mathrm{c}}){\beta\c_4^{-1}}{ W_{-1}\left((\rho-\rho_{\mathrm{c}})\beta/\c_4\right)^{-1}(1+o(1)) }&\text{ if }d=4\, ,\\
   -(\rho-\rho_{\mathrm{c}})^2 2\beta^2\pi^2(1+o(1)) &\text{ if }d=3\, ,
    \end{cases}
\end{equation}
which for $\d=4$ means $\bmu(\rho)=\Ocal\left(\left(\rho-\rho_{\mathrm{c}}\right)^{1+\e}\right)$.
\end{lemma}
%%%%%%%%%%%%%%%%%%.,%%%%%%%%%%%%%%%%%%%%%%%%%%%%%%%%%%%%%%
\begin{proofsect}{\textbf{Proof of Lemma \ref{LemDerivRat}}}
For $\mu<0$, $\brho(\mu)$ is smooth due to the exponential factor in the sum.\\
We expand in general
\begin{equation}\label{EquationExpnasion-General}
    \brho(\mu)-\brho(0)=\beta^{-1}(-\mu)^{d/2-1}\c_d\sum_{j\ge 1}(-\mu\beta) \left(\ex^{\beta \mu j}-1\right)(-\beta \mu j)^{-d/2}\, .
\end{equation}
For $d=3$, this implies that
\begin{equation}\label{EquationExpnasion-dimthree}
    \brho(\mu)-\brho(0)\sim\c_3 \beta^{-1}(-\mu)^{1/2}\int_{0}^\infty (\ex^{-x}-1)x^{-3/2}\d x=-(\beta\pi)^{-1}(-2\mu)^{1/2}\, .
\end{equation}
Indeed, by the convergence of the Riemann integral:
\begin{equation}
    \sum_{j\ge 1}(-\mu\beta) \left(\ex^{\beta \mu j}-1\right)(-\beta \mu j)^{-d/2}=\int_{0}^\infty (\ex^{-x}-1)x^{-3/2}\d x(1+o(1))\, ,\quad\text{as }(-\beta\mu)\to 0\, .
\end{equation}
Therefore, $\brho(\mu)-\brho(0)\sim-(\beta\pi)^{-1}(-2\mu)^{1/2}$. By relabelling the variables, we can see that
\begin{equation}
     \rho-\rho_{\mathrm{c}}\sim -[-2\bmu(\rho)]^{1/2}(\beta\pi)^{-1}\quad\Longleftrightarrow\quad\bmu(\rho)\sim-2\beta^2\pi^2(\rho-\rho_{\mathrm{c}})^2\, .
\end{equation}
For $d=4$, we use the same argument but with an additional truncation. Fix $\e\in(0,1/2)$ and observe that
\begin{equation}
    \c_4\mu^2\sum_{j=1}^{-\mu^{\e-1}}\left(\ex^{\beta\mu j}-1\right)(\beta\mu j)^{-2}\sim \c_4(1-\e)\mu\beta^{-1}\log(-\mu^{-1})\, .
\end{equation}
We also have that for some universal $C>0$
\begin{equation}
    \mu^2\sum_{j=-\mu^{\e-1}}^{-\mu^{-1}}\left(\ex^{\beta\mu j}-1\right)(\beta\mu j)^{-2}\le C \mu\e\log(-\mu^{-1})\, ,
\end{equation}
and
\begin{equation}
    \mu^2\sum_{j=-\mu^{-1}}^{\infty }\left(\ex^{\beta\mu j}-1\right)(\beta\mu j)^{-2}\le C\mu \, .
\end{equation}
Combining the last three formulae with the expansion from Equation \eqref{EquationExpnasion-General} and letting $\e\downarrow 0$, we can conclude that
\begin{equation}
    \brho(\mu)-\brho(0)\sim \c_4\mu\beta^{-1}\log(-\mu^{-1})\, .
\end{equation}
Recall that $y\ex^y=x$ for $\ex^{-1}\le x<0$ if and only if $y=W_{-1}(x)$. From this, it follows that
\begin{equation}
    \bmu(\rho)\sim (\rho-\rho_{\mathrm{c}})\frac{\beta}{\c_4 W_{-1}\left((\rho-\rho_{\mathrm{c}})\beta/\c_4\right) }\, .
\end{equation}
For $d\ge 5$, the result follows from the implicit function theorem and the fact that the density $\brho$ is differentiable. 
This concludes the proof.\qed
\end{proofsect}
\begin{lemma}\label{lemm_deriv_rf}
As $h\downarrow 0$
\begin{equation}\label{EquationApproximationI}
    I_0(\rho_{\mathrm{c}}-h)=I(\rho_{\mathrm{c}}-h)\sim  \begin{cases}
   \frac{h^{2}\beta }{2\brho'(0)}&\text{ if }d\ge 5\, ,\\
    \frac{-2h^2\beta^2}{\c_4W_{-1}\left(-\frac{h\beta}{\c_4}\right)}&\text{ if }d=4\, ,\\
  4{h^3 \beta^3\pi^2}&\text{ if }d=3\, .
    \end{cases}
\end{equation}
\end{lemma}
Note the faster-than-quadratic decay for $d=3,4$. This can be expected as $N_\L$ does not have a second moment in these dimensions under $\P_{\L,\beta,0}$. Indeed, this is because $\brho'(0)$ does not exist for these dimensions.
\begin{proofsect}{\textbf{Proof of Lemma \ref{lemm_deriv_rf}}}
The slightly awkward proof is warranted by the fact that we cannot simply apply the chain rule twice to observe the cancellations. Indeed, this is only possible for $d\ge 7$. Therefore, we expand our function into leading-order term plus remainder.\\
We begin with the case $d\ge 5$: we write $\bmu(\rho_{\mathrm{c}}-x)=-x\brho'(0)^{-1}+\e(x)$, where $\e(x)=\bmu(\rho_{\mathrm{c}}-x)+x\brho'(0)^{-1}$. We then expand
\begin{equation}
    P(h)=P(0)+h\beta\rho_{\mathrm{c}}+h^2\beta\brho'(0)/2+o(h^2)\, .
\end{equation}
This implies that
\begin{equation}\label{Eq4221}
    P(\bmu(\rho_{\mathrm{c}}-x))-P(0)=-x\brho'(0)^{-1}\beta\rho_{\mathrm{c}}+\beta\rho_{\mathrm{c}}\e(x)+x^2\brho'(0)^{-1}\beta/2-\beta^2x\e(x)+o(\e(x)x)\, .
\end{equation}
We furthermore expand
\begin{equation}\label{Eq4222}
    \beta (\rho_{\mathrm{c}}-x)\bmu(\rho_{\mathrm{c}}-x)=-\beta \rho_{\mathrm{c}} x\brho'(0)^{-1}+\beta\rho_{\mathrm{c}}\e(x)+x^2\beta\brho'(0)^{-1}+o(x\e(x))\, .
\end{equation}
Recall that $I(\rho_{\mathrm{c}}-h)$
\begin{equation}
    I(\rho_{\mathrm{c}}-h)=\beta (\rho_{\mathrm{c}}-x)\bmu(\rho_{\mathrm{c}}-x)-P(\bmu(\rho_{\mathrm{c}}-x))+P(0)\, ,
\end{equation}
for $h\downarrow 0$. Note that substituting Equation \ref{Eq4221} and Equation \eqref{Eq4222} into the above leads to the term $\beta\rho_{\mathrm{c}}\e(x)$ appearing with opposite sign. This implies that Equation \eqref{EquationApproximationI} is equal to
\begin{equation}
    \beta (\rho_{\mathrm{c}}-x)\bmu(\rho_x-x)- P(\bmu(\rho_{\mathrm{c}}-x))+P(0)=\frac{x^2\beta}{2\brho'(0)}+\Ocal\left(x^{5/2}\right)\, .
\end{equation}
This gives the result for $d\ge 5$. For $d=3$, we can expand
\begin{equation}
    P(h)=P(0)+h\beta\rho_{\mathrm{c}}-h^{3/2}\sqrt{2}\pi+o(h^{3/2})\, ,
\end{equation}
and thus, following the same route as for $d\ge 5$,
\begin{equation}
    \beta (\rho_{\mathrm{c}}-x)\bmu(\rho_{\mathrm{c}}-x)- P(\bmu(\rho_{\mathrm{c}}-x))+P(0)=4x^3\beta^3\pi^2+o(x^3)\, .
\end{equation}
For $d=4$,
\begin{equation}
    \beta (\rho_{\mathrm{c}}-x)\bmu(\rho_{\mathrm{c}}-x)=-\frac{\beta^2\rho_{\mathrm{c}} x}{\c_4W_{-1}\left(-\frac{x\beta}{\c_4}\right)}+\frac{\beta^2x^2}{\c_4W_{-1}\left(-\frac{x\beta}{\c_4}\right)}+\beta\rho_{\mathrm{c}}\e(x)+o(\e(x)x)\, .
\end{equation}
We also expand
\begin{equation}
    P(h)=P(0)+h\beta\rho_{\mathrm{c}}+h^2\c_4\log(-h^{-1})(1+o(1))\, .
\end{equation}
Note that 
\begin{equation}
\begin{split}
    \c_4\left[\frac{-x\beta}{\c_4W_{-1}\left(-\frac{x\beta}{\c_4}\right)}+\e(x)\right]^2\log\left(\frac{-x\beta}{\c_4W_{-1}\left(-\frac{x\beta}{\c_4}\right)}+\e(x)\right)\sim -\frac{x^2\beta^2}{\c_4W_{-1}\left(-\frac{x\beta}{\c_4}\right)}\, .
\end{split}
\end{equation}
Observing the same cancellations as before concludes the proof.\qed 
\end{proofsect}
We are now introducing the \textit{truncated} large deviation rate function:
We define $\P_\Delta^\mathrm{q}$ the PPP with intensity measure $M_\Delta^\mathrm{q}$ where
\begin{equation}
    M_\Delta^\mathrm{q}=\sum_{x\in\Delta}\sum_{j=1}^{q_n}\frac{1}{j}\P_{x,x}^{\beta j}\, .
\end{equation}
Here, $q_n$ is the same scale mentioned in Theorem \ref{TheoremMain}.\\
The log moment generating function of $N_o$ under $\P_\L^\mathrm{q}$ is given by
\begin{equation}
    \phi^\mathrm{q}(t)=
    P^\mathrm{q}(\mu+t/\beta)-P^\mathrm{q}(\mu)\, 
\end{equation}
where
\begin{equation}
     P^\mathrm{q}(\mu)=\sum_{j=1}^{q_n}\frac{\ex^{\beta \mu j}}{j}\p_{\beta j}(0)\, .
\end{equation}
We also define
\begin{equation}
    \brho^\mathrm{q}(\mu)=\sum_{j=1}^{q_n}{\ex^{\beta \mu j}}\p_{\beta j}(0)\, .
\end{equation}
Set $\bmu^\mathrm{q}(x)$ the unique $\mu$, such that $\brho^\mathrm{q}(\mu)=x$.
The Legendre transform of $\phi^\mathrm{q}$ is given by $I^q$, with
\begin{equation}
    I^\mathrm{q}(x)=\beta x\left(\bmu^\mathrm{q}(x)-\mu\right)-P^\mathrm{q}\left(\bmu^\mathrm{q}(x)\right)+P^\mathrm{q}(\mu)\, .
\end{equation}
% Note the general bound: for any $\delta_n\to\zero$, we have that
% \begin{equation}
%     \rho^Q(\delta_n)-\rho^Q(0)=\Ocal\left(\delta_n\int_1^{q_n}x^{-d/2+1}\d x\right)\, .
% \end{equation}
% Note that the sequence with suffix $n$ is separate from the $N$-sequence. In general, we have a finer control: 
\begin{lemma}\label{finecontrolapproxrate}
We have that
    \begin{equation}
        \forall \e>0 \,\exists C>0\, \forall \rho\in [0,\rho_{\mathrm{c}}-\e]\colon \,\abs{I^\mathrm{q}(\rho)-I(\rho)}\le \Ocal\left(\ex^{-\beta C q_n}\right)
    \end{equation}
    as $n\to\infty$.
\end{lemma}
\begin{proofsect}{\textbf{Proof of Lemma \ref{finecontrolapproxrate}}}
For $\mu<0$,
\begin{equation}
    \brho^\mathrm{q}(\mu)=\brho(\mu)+\Ocal\left(\ex^{\beta \mu q_n}\right)\, .
\end{equation}
Note that for $\rho>0$ bounded away from $\rho_{\mathrm{c}}$, $\bmu$ is differentiable at $\rho$ with uniformly bounded derivative. Furthermore, taking $n$ sufficiently large, we may assume without loss of generality that $y<0$, where $y$ is the value such that $\bmu^\mathrm{q}(y)=\rho$. Note that
\begin{equation}
    \bmu^\mathrm{q}(\rho)-\bmu(\rho)=y-\bmu(\brho^\mathrm{q}(y))=\Ocal\left(\ex^{\beta yq_n}\right)\, .
\end{equation}
From there on, the result follows in the manner of Lemma \ref{LemDerivRat}.\qed
\end{proofsect}

\subsection{Calculation of the partition function, supercritical}\label{SubsectionCalculationPartitionFUnction}
The goal in this subsection will be to calculate the value of the partition function
\begin{equation}
   Z^{(\mathrm{Can, HYL})}_{\Lambda,\beta,\rho}=\E_{\L,\beta,0}\left[\ex^{-\beta\Ham}, N_\L= \rho\abs{\L}\right]\quad\text{for }\rho>\rho_{\mathrm{c}}\, .
\end{equation}
Here, and henceforth, we use the notation $\E\left[F,A\right]$ to abbreviate $\E\left[F\1_A\right]$, for $F$ a function and $A$ a set. We introduce the parameter $\rho_\mathrm{e}=\rho-\rho_{\mathrm{c}}>0$. We assume for ease of reading that $\rho\abs{\L}\in\N$.\\
To calculate the partition function, we introduce the quadratic
\begin{equation}
    Q(t)=Q_{\rho_\mathrm{e}}(t)=\beta b\left[(t+\rho_\mathrm{e})^2-\beta\rho_\mathrm{e}^2\right]/2\, .
\end{equation}
This polynomial gives the gain of the function $b\beta x^2/2$ at $x=t+\rho_\mathrm{e}$ with respect to $b\beta\rho_\mathrm{e}^2/2$. This quantity is crucial: it represents the maximal potential energy gain if we increase the density of the interlacements from $\rho_\mathrm{e}$ to $t+\rho_\mathrm{e}$.

Later we see that the term $Q(t)- I(\rho_{\mathrm{c}}-t)$ is the total ``cost" of adding $t$ density to the interlacements, as the large deviation cost is given by $I$. Motivated by that, we define
\begin{equation}\label{EquationH0}
    S_o=\sup_{t\in [0,\rho_{\mathrm{c}}]}\left\{Q(t)- I(\rho_{\mathrm{c}}-t)\right\}\, .
\end{equation}
Crucial is the following parameter: set
\begin{equation}\label{EquationrhoS}
    \rho_{\mathrm{S}}=\rho_{\mathrm{S}}(b,\rho_\mathrm{e},d) \quad\text{ such that}\colon\quad Q(\rho_{\mathrm{S}})- I(\rho_{\mathrm{c}}-\rho_{\mathrm{S}})=S_o\, .
\end{equation}
$\rho_{\mathrm{S}}+\rho_\mathrm{e}$ will be the density of the random interlacements and therefore the $\rho_{\mathrm{S}}$ represents the extra density of interlacement arising from the non-mean-field part of the interaction.
\begin{lemma}\label{LemPropertiesOfJump}
The parameter $\rho_{\mathrm{S}}\in (0,\rho_{\mathrm{c}})$ is well defined. Furthermore, $\rho_{\mathrm{S}}(b,\rho_\mathrm{e},d)=o(1)$ as $\rho_\mathrm{e}\to 0$, given $b\le\tfrac{1}{\brho'(0)}$ and $d\ge 5$. Otherwise, 
\begin{equation}\label{EquationJumpRhoB}
    \lim_{\rho_\mathrm{e}\downarrow 0}\rho_{\mathrm{S}}(b,\rho_\mathrm{e},d)>0\, .
\end{equation}
\end{lemma}
\begin{figure}[ht]
    \centering
    \begin{subfigure}[t]{\textwidth}
    \centering
    \begin{tikzpicture}[scale=1.2]
        \draw[->] (0,-2) -- (0,2)node[above]{$Q(\rho)-I(\rho_{\mathrm{c}}-\rho)$};
        \draw[->] (0,0) -- (4.2,0)node[below]{$\rho$};
        \draw[dashed] (0,0) -- (4,1);
        \draw[dashed] (4,-2) -- (4,2)node[right]{$\rho_{\mathrm{c}}$};
        \draw[thick] (0,0) to [out=14,in=180] (2,1) to [out=0,in=90] (4,-1);
        \draw[dashed] (2,1) -- (2,0)node[below]{$\rho_{\mathrm{S}}$};
        \draw[->,very thick] (4.5,0) -- (5.5,0)node[above left]{$\rho_\mathrm{e}\to0$};
        \draw[->] (6,-2) -- (6,2)node[above]{$Q(\rho)-I(\rho_{\mathrm{c}}-\rho)$};
        \draw[->] (6,0) -- (10.2,0)node[below]{$\rho$};
        \draw[dashed] (10,-2) -- (10,2)node[right]{$\rho_{\mathrm{c}}$};
        \draw[thick] (6,0) to [out=0,in=180] (8,0.5) to [out=0,in=90] (10,-2);
        \draw[dashed] (8,0.5) -- (8,0)node[below]{$\rho_{\mathrm{S}}$};
    \end{tikzpicture}
    \caption{$d=3,4$, or $d\geq 5$ and $b>\tfrac{1}{\brho'(0)}$}
    \end{subfigure}
    \hfill
    \begin{subfigure}[t]{\textwidth}
    \centering
    \begin{tikzpicture}[scale=1.2]
        \draw[->] (0,-2) -- (0,2)node[above]{$Q(\rho)-I(\rho_{\mathrm{c}}-\rho)$};
        \draw[->] (0,0) -- (4.2,0)node[below]{$\rho$};
        \draw[dashed] (0,0) -- (4,1);
        \draw[dashed] (4,-2) -- (4,2)node[right]{$\rho_{\mathrm{c}}$};
        \draw[thick] (0,0) to [out=14,in=180] (2,0.25) to [out=0,in=90] (4,-1);
        \draw[dashed] (2,0.25) -- (2,0)node[below]{$\rho_{\mathrm{S}}$};
        \draw[->,very thick] (4.5,0) -- (5.5,0)node[above left]{$\rho_\mathrm{e}\to0$};
        \draw[->] (6,-2) -- (6,2)node[above]{$Q(\rho)-I(\rho_{\mathrm{c}}-\rho)$};
        \draw[->] (6,0) -- (10.2,0)node[below]{$\rho$};
        \draw[dashed] (10,-2) -- (10,2)node[right]{$\rho_{\mathrm{c}}$};
        \draw[thick] (6,0) to [out=0,in=90] (10,-2);
        \draw (6,0) node[above right]{$\rho_{\mathrm{S}}\to0$};
    \end{tikzpicture}
    \caption{$d\geq 5$ and $b\le\tfrac{1}{\brho'(0)}$}
    \end{subfigure}
    \caption{Sketch of the behaviour of $\rho_{\mathrm{S}}$ as $\rho_\mathrm{e}\to0$. The diagonal dashed lines follow $b\rho_\mathrm{e} \rho$. In the first case the function $Q(\rho)-I(\rho_{\mathrm{c}}-\rho)$ initially goes above this line and as $\rho_e\to0$ the maximising argument $\rho_{\mathrm{S}}$ stays away from $0$. In the second case the function stays below the line and $\rho_{\mathrm{S}}\to0$.}
    \label{figure3}
\end{figure}
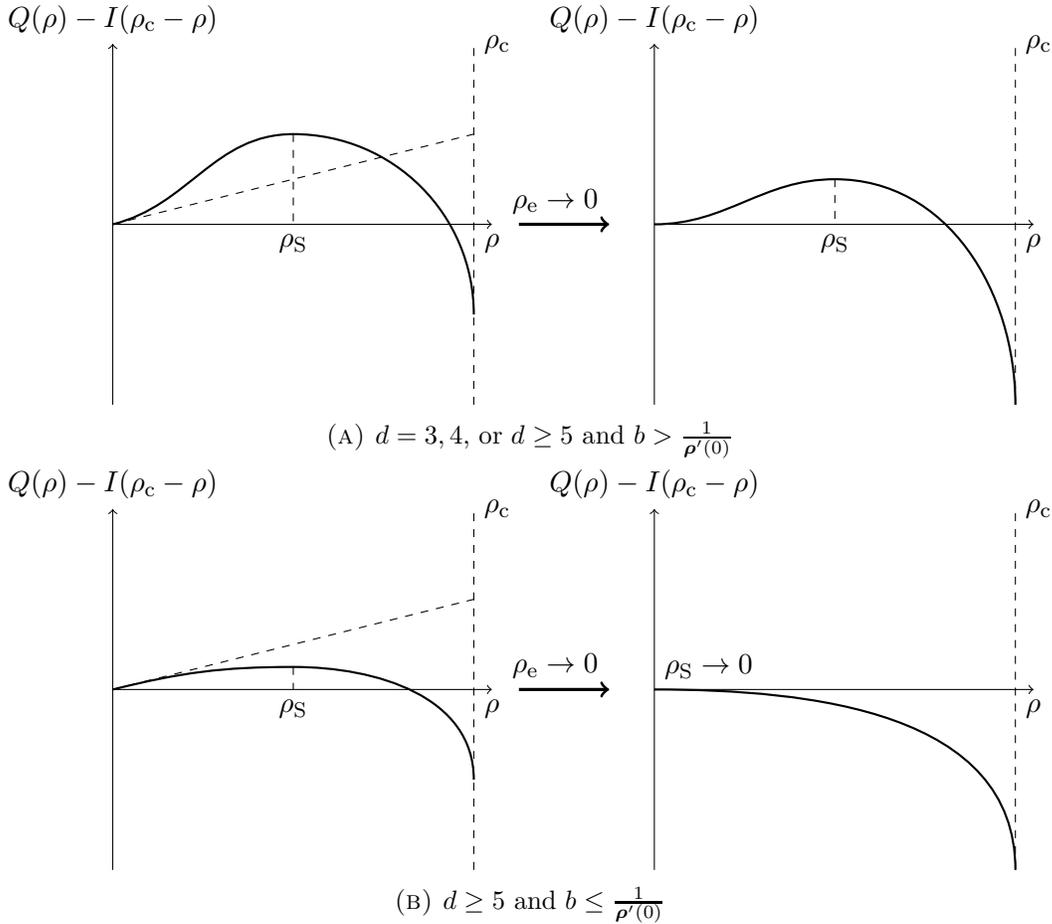
\begin{proofsect}{\textbf{Proof of Lemma \ref{LemPropertiesOfJump}}}
We supply Figure \ref{figure3} to illustrate the proof. Also recall that sketches of $\bmu,\brho$ and their derivatives are given in Figure \ref{figurethermodynamic functions}. We first show that $\rho_{\mathrm{S}}$ is well defined, i.e. that there is exactly one $\rho\in (0,\rho_{\mathrm{c}})$ such that $Q(\rho_{\mathrm{S}})-I(\rho_{\mathrm{c}}-\rho_{\mathrm{S}})=S_0$. For this, fix $\rho_\mathrm{e}>0$ and note that for $x\in (0,\rho_{\mathrm{c}})$
\begin{equation}\label{EquationDerivativeRateFunction}
    I'(x)=\beta \bmu(x)\quad\text{and}\quad I''(x)=\beta \bmu'(x)\, .
\end{equation}
This implies that
\begin{equation}
    \frac{\d}{\d x}\left[Q(x)- I(\rho_{\mathrm{c}}-x)\right]=\beta b(x+\rho_\mathrm{e})+\beta\bmu(\rho_{\mathrm{c}}-x)\text{ and }\frac{\d^2}{\d x^2}\left[Q(x)- I(\rho_{\mathrm{c}}-x)\right]=\beta b-\beta\bmu'(\rho_{\mathrm{c}}-x)\, .
\end{equation}
Note $Q(x)-I(\rho_{\mathrm{c}}-x)$ is increasing for $x\in[0,\e)$ (for some $\e>0$ small enough). Indeed, from Lemma \ref{lemm_deriv_rf} we know that for small $x$, the function $x\mapsto I(\rho_{\mathrm{c}}-x)$ grows at most quadratically in $x$ but $Q(x)$ is of linear growth. On the other hand, as $\bmu(\rho_{\mathrm{c}}-x)\to -\infty$, as $x\to\rho_{\mathrm{c}}$. This means that $Q(x)-I(\rho_{\mathrm{c}}-x)$ decreases for large enough $x$. Note that $x\mapsto \beta b-\beta\bmu'(\rho_{\mathrm{c}}-x)$ is decreasing monotonically. By the continuity of the involved functions, it follows that $Q(x)-I(\rho_{\mathrm{c}}-x)$ must attain its supremum inside the interval $(0,\rho_{\mathrm{c}})$. As the second derivative is decreasing, this supremum is indeed a maximum.\\
Next, we examine the situation as $\rho_\mathrm{e}\to 0$. We expand $Q(x)-I(\rho_{\mathrm{c}}-x)$ for small $x$ as
\begin{equation}
    \begin{cases}
    \frac{\beta bx^2}{2}+\beta b\rho_\mathrm{e} x-\frac{x^{2}\beta }{2\brho'(0)}(1+o(1))&\text{ if }d\ge 5\, ,\\
    \frac{\beta bx^2}{2}+\beta b\rho_\mathrm{e} x-\frac{2h^2\beta^2}{\c_4W_{-1}\left(-\frac{x\beta}{\c_4}\right)}(1+o(1))&\text{ if }d=4\, ,\\
  \frac{\beta bx^2}{2}+\beta b\rho_\mathrm{e} x-4{x^3 \beta^3\pi^2}(1+o(1))&\text{ if }d=3\, .
    \end{cases}
\end{equation}
Therefore, as long as $d=3,4$ or $b>\brho'(0)^{-1}$, $Q(x)-I(\rho_{\mathrm{c}}-x)$ increases in a small neighbourhood around zero, disregarding of the value of $\rho_\mathrm{e}$.\\
On the other hand, for $d\ge 5$, by the inverse function theorem
\begin{equation}
    \frac{1}{\beta}\frac{\d^2}{\d x^2}\left[Q(x)- I(\rho_{\mathrm{c}}-x)\right]=b-\bmu'(\rho_{\mathrm{c}}-x)=b-\frac{1}{\brho'(\bmu(\rho_{\mathrm{c}}-x))}\, .
\end{equation}
Thus, if $b<\brho'(0)^{-1}$, the second derivative is negative. Following from that, with $m=(b-1/\brho'(0))/2<0$, we have the expansion
\begin{equation}
   Q(x)-I(\rho_{\mathrm{c}}-x)=\beta mx^2+\beta b \rho_\mathrm{e} x+o(x^2)\, , 
\end{equation}
where the small-o term is independent of $\rho_\mathrm{e}$. This is (asymptotically) a parabola with zeros at the origin and at $x\sim -{\rho_\mathrm{e}}{m}$. Thus, $\rho_{\mathrm{S}}=\Ocal\left(\frac{\rho_\mathrm{e}}{-2m}\right)=o(1)$ and the result follows. For $b=\brho'(0)^{-1}$, the argument is similar, $\beta b(x+\rho_\mathrm{e})+\beta\bmu(\rho_{\mathrm{c}}-x)= b\beta\rho_\mathrm{e}+\Ocal(x^{3/2})$ and one can show that the zero of that function is of order $\Ocal(\rho_\mathrm{e})$. \qed
% and $b-\beta\bmu'(\rho_{\mathrm{c}}-x)$ is strictly decreasing on $(0,\rho_{\mathrm{c}})$, the claim follows.\\
% Recall that by Lemma \ref{lemm_deriv_rf}, $I(\rho_{\mathrm{c}}-x)$ decays faster $Q(x)$ \textit{for any $\rho_\mathrm{e}$} unless $d\ge 5$ and $b\le \tfrac{\beta}{\brho'(0)}$. Indeed, for $d=3,4$, $I(\rho_{\mathrm{c}}-x)$ decays faster than order square and for $d\ge 5$ and $b>\tfrac{\beta}{\brho'(0)}$ the square term is bigger than that of $Q(\rho)$. From that, one can deduce Equation \ref{EquationJumpRhoB}. This concludes the proof.\qed
\end{proofsect}

Abbreviate
\begin{equation}\label{Equationrhob}
   \rhob=\rho_{\mathrm{S}}+\rho_\mathrm{e}\, ,\quad\text{ which will turn out to be the total density of the condensate.} 
\end{equation}
 Indeed, $\rho_\mathrm{e}$ will be the contribution from the free gas, while $\rho_{\mathrm{S}}$ comes from the Hamiltonian.\\ 
Our goal is to calculate the partition function by expanding around $\rhob$. For this, we split our loop soup into long and short loops:
\begin{equation}
    \Ns=\Ns_\L=\sum_{\omega\in\eta_\L}\ell(\omega)\1\{\ell(\omega)<q_n\}\quad\text{and}\quad \Nl=\Nl_\L=N_\L-\Ns_\L\, .
\end{equation}
Using the independence of $\Ns$ and $\Nl$ from the Poisson property, we expand
\begin{equation}
    \E_{\L,\beta,0}\left[\ex^{-\Ham},N_\L= \rho\abs{\L}\right]=\int_0^{\rho\abs{\L}} \E_{\L,\beta,0}\left[\ex^{-\Ham},\Nl= \rho\abs{\L}-x\right]\d\P_{\L,\beta,0}\left(\Ns=x\right)\, .
\end{equation}
Indeed, $\Ham$ is measurable with respect to $\Nl$. Let us analyse the behaviour of $\P_{\L,\beta,0}\left(\Nl=x\right)$:
\begin{lemma}\label{LemDistOfQN}
For $c>0$ fixed,
\begin{equation}
    \P_{\L,\beta,0}\left(\Nl= x\right)\sim \frac{\beta\abs{\L}\c_d}{(\beta x)^{d/2+1}} \, ,
\end{equation}
uniformly in $c\abs{\L}>x>c^{-1}\abs{\L}$.
\end{lemma}
\begin{proofsect}{\textbf{Proof of Lemma \ref{LemDistOfQN}}}
% Note that by \cite{berger2019notes}, we have that
% \begin{equation}
%     \P_{\L,\beta,0}\left(N=x\right)\sim \P_{\L,\beta,0}\left(\exists \omega\colon \ell(\omega)=x\right)
% \end{equation}
Define the auxiliary sequence $r_N$ as
\begin{equation}
    r_n=M_{\L,\beta,0}[\ell(\omega)\ge q_n]=\Ocal\left(\abs{\L}q_n^{-d/2}\right)=o\left(1\right)\, .
\end{equation}
Note here that $a_n=o(q_n)$ is needed for some $d\ge 3$. By the fundamental properties of Poisson point processes,
\begin{equation}\label{EquationHelp123}
    \P_{\L,\beta,0}\left(\exists \omega_1,\ldots,\omega_k\colon \ell(\omega_i)\ge q_n\text{ for all }i=1,\ldots,k\right)=\Ocal\left(r_n^k\right)\, .
\end{equation}
Let $\eta_q$ denote the number of loops which are longer than $q_n$. We then split
\begin{equation}
    \P_{\L,\beta,0}\left(\Nl= x\right)=\P_{\L,\beta,0}\left(\Nl= x,\eta_q\le d\right)+\P_{\L,\beta,0}\left(\Nl= x,\eta_q>d\right)\, .
\end{equation}
By Equation \eqref{EquationHelp123}, the second term is negligible. Furthermore, note that
\begin{equation}
    \P_{\L,\beta,0}\left(\Nl= x,\eta_q=1\right)\sim M_{\L,\beta,0}[\ell(\omega)\ge q_n]\sim \frac{\beta\abs{\L}\c_d}{(\beta x)^{d/2+1}}\, .
\end{equation}
Now for $k\in\{2,\ldots,d\}$ fixed
\begin{equation}
    \P_{\L,\beta,0}\left(\Nl= x,\eta_q=k\right)\le \P_{\L,\beta,0}\left(\exists\omega\colon \ell(\omega)\ge x/k\text{, and }\exists \omega_1,\ldots,\omega_{k-1}\colon \omega_i\ge q_n,\,\forall i+1,\ldots,k-1\right),
\end{equation}
which implies that
\begin{equation}
    \P_{\L,\beta,0}\left(\Nl= x,\eta_q=k\right)\le C(k) \P_{\L,\beta,0}\left(\Nl= x,\eta_q=1\right)r_n^{k-1}\, ,
\end{equation}
where $C(k)$ is some $k$-dependent constant. Hence, the event $\{\eta_q=1\}$ is the only relevant one in the limit.
% We expand
% \begin{equation}\label{EquationTem11}
%      \P_{\L,\beta,0}\left(\Nl= x\right)=\sum_{k\ge 1}\sum_{x_1,\ldots,x_k\in\Pfrak_x^q} \P_{\L,\beta,0}\left(\exists (\omega_i)_{i=1}^k\colon \ell(\omega_i)=x_i\right)\, ,
% \end{equation}
% where $\Pfrak_x^q$ is the collection of tuples $x_1,\ldots,x_k$ such that $x_i>q_n$ and $x_1+\ldots+x_k=x$, for any $k\ge 1$. Using the properties of the Poisson process, we can estimate
% \begin{equation}
%     \P_{\L,\beta,0}\left(\exists (\omega_i)_{i=1}^k\colon \ell(\omega_i)=x_i\right)=\Ocal\left(\abs{\L}^k\ex^{-(d/2+1)\sum_{i=1}^k\log x_i}\right)\, .
% \end{equation}
% Indeed, the $\abs{\L}^k$ is the gain we get from having $\abs{\L}$ choices to start each loop (at a square of unit volume) and the exponential term is the cost of sampling a loop that the prescribed length. Taking the sum over all tuples in $\Pfrak_x^q$ gives
% \begin{equation}
%     \sum_{x_1,\ldots,x_k\in\Pfrak_x^q} \P_{\L,\beta,0}\left(\exists (\omega_i)_{i=1}^k\colon \ell(\omega_i)=x_i\right)=\Ocal\left(\abs{\L}^{-k(d/2-1)}\right)\, .
% \end{equation}
% This implies that the leading order term in Equation \eqref{EquationTem11} is the term for $k=1$ and therefore
% \begin{equation}
%   \P_{\L,\beta,0}\left(\Nl= x\right)\sim \P_{\L,\beta,0}\left(\exists \omega\colon\ell(\omega)=x\right)\sim \frac{\beta\abs{\L}\c_d}{(\beta x)^{d/2+1}} \, .
% \end{equation}
This concludes the proof.\qed
\end{proofsect}
Given that $\Nl=x$, the Hamiltonian becomes predictable:
\begin{cor}\label{CorConcentrationHamil}
For any $c>0$
\begin{equation}
    \E_{\L,\beta,0}\left[\ex^{-\beta \Ham(\eta)}|\Nl= x\right]\sim \E_{\L,\beta,0}\left[\ex^{-\beta\Ham(\eta)}|\exists \omega\colon \ell(\omega)=x\right]=\ex^{\beta bx^2/(2\abs{\L})}\, ,
\end{equation}
for any $x$ with $c\abs{\L}>x>c^{-1}\abs{\L}$, uniformly.
\end{cor}
\begin{proofsect}{\textbf{Proof of Corollary \ref{CorConcentrationHamil}}}
Let the decreasing sequence $(s_i)_{i=1}^{M_N}$ be the different values $-\Ham$ can attain, restricted to the set $\Nl= x$. One has that
\begin{equation}
    s_1=bx^2/(2\abs{\L})\quad\text{ and }s_2=s_1-\frac{bq_nx}{2\abs{\L}}(1+o(1))\, .
\end{equation}
Indeed, a simple calculation using Lagrange multipliers shows that the maximal value of $\Ham$ is achieved by placing all the particles in the same cycle, thereby proving the first equality above. Given the requirement $x>c^{-1}\abs{\L}$ and the growth of $(q_n)_n$, it follows that $s_2$ is negligible compared to $s_1$. 

For the second equality, we observe that this argument can be used inductively, i.e. the next best strategy is to place the particles in two cycles. As the cycle-lengths are bounded from below by $q_n$, the second equality follows. 
Lemma \ref{LemDistOfQN} asserts that the events $\{\Nl= x\}$ and $\{\exists \omega\colon \ell(\omega)=x\}$ have asymptotically the same mass. This concludes the proof.\qed
\end{proofsect}
Next, we give an asymptotic relation for the distribution of the short loops.
\begin{lemma}\label{NoLongLoops}
For any $\e>0$ and for $y\in \left[0,(\rho_{\mathrm{c}}-\e)\abs{\L}\right]$, it holds
\begin{equation}
    \P_{\L,\beta,0}\left(\Ns=y\right)\sim \P_{\L,\beta,0}\left(N_\L=y\right)\, .
\end{equation}
\end{lemma}
\begin{proofsect}{\textbf{Proof of Lemma \ref{NoLongLoops}}}
Note that we can find $c_o>0$ such that
\begin{equation}
      {c_o^{-1}}\ex^{-\abs{\L}I^{\mathrm{q}}(y)}\ge\P_{\L,\beta,0}\left(\Ns=y\right)\ge {c_o}\ex^{-\abs{\L}I^{\mathrm{q}}(y)}\, .
\end{equation}
Indeed, this expansion can be found in \cite{bahadur1960deviations}.

Set $p_n=q_n/\abs{\L}$ and expand
\begin{equation}
   \P_{\L,\beta,0}\left(N_\L=y,\exists\omega\colon \ell(\omega)>q_n\right)\le \int_{q_n}^{\rho\abs{\L}}\P_{\L,\beta,0}\left(\Ns=y-s\right)\d s=\Ocal\left(\ex^{-\abs{\L}I^{\mathrm{q}}(y-p_n/2)}\right)\, .
\end{equation}
By Lemma \ref{finecontrolapproxrate}, we may replace $\abs{\L}I^{\mathrm{q}}(y-p_n/2)$ by $\abs{\L}I(y-p_n/2)(1+o(1))$, where the $o(1)$ term is bounded from above by $\Ocal\left(\ex^{-Cq_n}\right)$. This gives
\begin{multline}
    \P_{\L,\beta,0}\left(N_\L=y,\exists\omega\colon \ell(\omega)>q_n\right)\le\Ocal\left(\ex^{-\abs{\L}I(y-p_n/2)(1+o(1))}\right)\\
    \le \Ocal\left(\ex^{-\abs{\L}I(y)}\right)\Ocal\left(\ex^{-\abs{\L}[I(y-p_n/2)-I(y)](1+o(1))}\right)
\end{multline}
However, invoking the same lemma again, we see that 
\begin{equation}
    I(y-p_n/2)-I(y)=\Ocal\left(p_n\right)\quad\Rightarrow\quad  \P_{\L,\beta,0}\left(N_\L=y,\exists\omega\colon \ell(\omega)>q_n\right)= \Ocal\left(\ex^{-\abs{\L}I(y)}\right)\Ocal\left(\ex^{-cq_n}\right)\, ,
\end{equation}
for some $c>0$. Here, we use the assumptions on $(q_n)_n$. On the other hand by \cite{bahadur1960deviations}, for some $d_o>0$
\begin{equation}
     \P_{\L,\beta,0}\left(N_\L=y\right)\sim \frac{d_o}{\sqrt{y}}\ex^{-\abs{\L}I(y)}\, .
\end{equation}
Therefore, by combining the two previous equations with the independence from the Poisson point process, we get
\begin{equation}
    \P_{\L,\beta,0}\left(N_\L=y\right)\sim\P_{\L,\beta,0}\left(N_\L=y,\forall\omega\colon \ell(\omega)<q_n\right)=\P_{\L,\beta,0}\left(\Ns=y\right)\, .
\end{equation}
This concludes the proof.\qed 
\end{proofsect}
We now split the partition function
\begin{multline}
%\begin{split}
    \E_{\L,\beta,0}\left[\ex^{-\beta\Ham},\bN_\L= \rho\right]=\E_{\L,\beta,0}\left[\ex^{-\beta\Ham},\bN_\L= \rho,\Nl\in\bB_T^{+}(\rhob)\right]\\
    +\E_{\L,\beta,0}\left[\ex^{-\beta\Ham},\bN_\L= \rho,\Nl\notin\bB_T^{+}(\rhob)\right]\, ,
%\end{split}
\end{multline}
where $\bB_T^{+}(\rhob)=\abs{\L}\rhob+\abs{\L}^{1/2}[-T,T]$, for $T>0$ which we will let diverge to $+\infty$ later. 
\begin{lemma}\label{LemmaTreash}
Set $S_1=S_0+b\beta\rho_\mathrm{e}^2/2$. There exists a positive, increasing function $T\mapsto\gamma_T$ diverging to $+\infty$ as $T$ diverges to $+\infty$, such that
\begin{equation}
    \E_{\L}\left[\ex^{-\beta\Ham},\bN_\L= \rho,\Nl\notin\bB_T^{+}(\rhob)\right]=\Ocal\left(\frac{\ex^{\abs{\L}S_1-\gamma_T}}{\abs{\L}^{d/2+1}}\right)\, .
\end{equation}
\end{lemma}
\begin{proofsect}{\textbf{Proof of Lemma \ref{LemmaTreash}}}
Recall that $Q(x)=\beta b[(x+\rho_\mathrm{e})^2-\rho_\mathrm{e}]/2$. On the event that $\Nl=x$, we bound
\begin{equation}
    -\Ham\le \abs{\L}\left[Q(x/\abs{\L}-\rho_\mathrm{e})+\rho_\mathrm{e}^2/2\right]\, .
\end{equation}
Thus 
\begin{multline}
      \E_{\L}\left[\ex^{-\beta\Ham},\bN_\L= \rho,\Nl\notin\bB_T^{+}(\rhob)\right]\\
      \le\ex^{\abs{\L}\frac{b\beta\rho_\mathrm{e}^2}{2}} \int_{0}^{\rho\abs{\L}}\ex^{\abs{\L}Q(x/\abs{\L}-\rho_\mathrm{e})}\1\{x\notin\bB_T^{+}(\rhob)\} \P_\L(\Nl=\floor{x})\d\P_\L\left(\Ns=\rho\abs{\L}-x\right)\, .
\end{multline}
Now, choose $\e,\delta>0$ such that $Q(\e)+b\beta\rho_\mathrm{e}^2/2\le S_1-\delta$. This is possible as $\rhob$ is strictly between 0 and $\rho_{\mathrm{c}}$. For such choice, we get that
\begin{equation}
    {\ex^{\abs{\L}\frac{b\beta\rho_\mathrm{e}^2}{2}}}\int_{0}^{(\rho_{\mathrm{c}}+\e)\abs{\L}}\ex^{\abs{\L}Q(x/\abs{\L}-\rho_\mathrm{e})}\1\{x\notin\bB_T^{+}(\rhob)\}\d\P_\L\left(\Ns=\rho\abs{\L}-x\right)\le \ex^{\abs{\L}({S_1-\delta})}\, .
\end{equation}
Thus, it remains to estimate (after a change of variables)
\begin{equation}
     \frac{\ex^{\abs{\L}\frac{b\beta\rho_\mathrm{e}^2}{2}}}{\abs{\L}^{d/2+1}}\int_{\e\abs{\L}}^{\rho_{\mathrm{c}}\abs{\L}}\ex^{\abs{\L}Q(x/\abs{\L})}\1\{x\notin\bB_T^+(\rho_{\mathrm{S}})\}\d\P_\L\left(\Ns=\rho_{\mathrm{c}}\abs{\L}-x\right)\,  .
\end{equation}
We recall that the $\abs{\L}^{d/2+1}$ factor comes from contribution of $\Nl$. Applying Lemma \ref{NoLongLoops}, we can bound the above integral by
\begin{equation}
    C\int_{\e\abs{\L}}^{\rho_{\mathrm{c}}\abs{\L}}\ex^{\abs{\L}Q(x/\abs{\L})}\1\{x\notin\bB_T^+(\rho_{\mathrm{S}})\}\d\P_\L\left(N=\rho_{\mathrm{c}}\abs{\L}-x\right)\,  ,
\end{equation}
for some $C>0$. Using \cite{bahadur1960deviations}, the above is bounded by
\begin{equation}
    C\int_{\e\abs{\L}}^{\rho_{\mathrm{c}}\abs{\L}}\ex^{\abs{\L}Q(x/\abs{\L})}\1\{x\notin\bB_T^+(\rho_{\mathrm{S}})\}\frac{\ex^{-\abs{\L}I(\rho_{\mathrm{c}}-x/\abs{\L})}}{\sqrt{x}}\d x\,  .
\end{equation}
Note that $x\mapsto Q(x)-I(\rho_{\mathrm{c}}-x)$ is differentiable at its minimum $\rho_{\mathrm{S}}$, and thus as $x\to \rhob$ 
\begin{equation}
Q(x/\abs{\L})-I(\rho_{\mathrm{c}}-x/\abs{\L})=S_0-C_q\left(\frac{x}{\abs{\L}}-\rhob\right)^2(1+o(1))\, ,
\end{equation}
where $C_q=b\beta-\beta\bmu(\rho_{\mathrm{c}}-\rhob)>0$.\\
This implies that
\begin{equation}
     \ex^{\abs{\L}\frac{b\beta\rho_\mathrm{e}^2}{2}}\int_{\e\abs{\L}}^{\rho_{\mathrm{c}}\abs{\L}}\ex^{\abs{\L}Q(x/\abs{\L})}\1\{x\notin\bB_T^+(\rho_{\mathrm{S}})\}\d\P_\L\left(\Ns=\rho_{\mathrm{c}}\abs{\L}-x\right)\le \Ocal\left(\frac{\ex^{\abs{\L}S_1-\gamma_T}}{\abs{\L}^{d/2+1}}\right)\, .
\end{equation}
This completes the proof.\qed
\end{proofsect}

% Note that
% \begin{equation}\label{EquationExponentialOutsideIntegral}
%     \begin{split}
%         \E_{\L}\left[\ex^{-\Ham},\bN_\L= \rho,\Nl\notin\bB_T^{+}(\rhob)\right]\le C \int_{0}^{\rho_\mathrm{e}\abs{\L}}\!\ex^{\left[Q(x/\abs{\L})+\rho_\mathrm{e}^2/2\right]\abs{\L}}\d \P_\L\left(\Ns\!=-x+\rho_{\mathrm{c}}\abs{\L}\notin\bB_\delta(\rho-\rhob)\right).
%     \end{split}
% \end{equation}
% {\color{red}MORE HERE}
% Indeed, the event that the short loops carry more than $\abs{\L}\rho_{\mathrm{c}}$ particles can be disregarded, since the the free system gives vanishing probability to it and the Hamiltonian discourages it. By approximation, we can bound the integral by
% \begin{equation}\label{EquationSupremum}
%     \exp\left(\abs{\L}\sup_{x\in [0,\rho_{\mathrm{c}}]\setminus\bB_T^{+}(\rhob)/\abs{\L}}\left(Q(x)-I(\rho_{\mathrm{c}}-x)\right)+\abs{\L}b\rho_\mathrm{e}^2/2\right)\le \ex^{\abs{\L}S_1-c_N}\, ,
% \end{equation}
% where $S_1=S_0+b\rho_\mathrm{e}^2/2$ and $c_N$ a positive, increasing sequence, diverging to infinity at some power-law speed.
Recall Varadhan's theorem (see \cite{dembo2009large}): for $(\P_n)_n$ satisfying a large deviation principle with rate function $i$
\begin{equation}
    \E_n\left[\ex^{nF(x)}\right]=\exp\left(n \sup_{x}\left\{F(x)-i(x)\right\}\big(1+o(1)\big)\right)\, ,
\end{equation}
for $F$ continuous and bounded above.

In the spirit of the standard Laplace approximation, this result was refined in \cite{martin1982laplace}: assume the supremum in the above equation attained at $0$. Then:
\begin{equation}
     \E_n\left[\ex^{nF(x)}\right]=\sqrt{1+\frac{F''(0)}{i''(0)}}\ex^{n[ F(0)-i(0)]}\big(1+o(1)\big)\, .
\end{equation}
We are now ready to compute the partition function to the required accuracy. 
\begin{lemma}\label{LemmaPartitionFunction}
For $\rho>\rho_{\mathrm{c}}$
\begin{equation}
    \E_{\L,\beta,0}\left[\ex^{-\beta\Ham},\bN_\L= \rho\right]\sim \sqrt{1+\frac{b}{ \bmu'(\rho_{\mathrm{c}}-\rho_{\mathrm{S}})}}\exp\left(S_1\abs{\L}\right)\frac{\beta\abs{\L}\c_d}{(\beta \rhob\abs{\L})^{d/2+1}}\, .
\end{equation}
\end{lemma}
\begin{proofsect}{\textbf{Proof of Lemma \ref{LemmaPartitionFunction}}}
By Lemma \ref{LemmaTreash}, we can reduce the question to calculating
\begin{equation}
    \E_{\L,\beta,0}\left[\ex^{-\beta\Ham},\bN_\L= \rho,\, \Nl\in\bB_T^{+}(\rhob)\right]\, .
\end{equation}
By Lemma \ref{NoLongLoops}, we can expand
\begin{equation}
   \E_{\L,\beta,0}\left[\ex^{-\beta\Ham},\bN_\L= \rho,\, \Nl\in\bB_T^{+}(\rhob)\right]\sim\int_{\bB_T^{+}(\rhob)} \E_{\L,\beta,0}\left[\ex^{-\beta\Ham},\Nl= x\right]\d\P_{\L,\beta,0}\left(N_\L=\rho\abs{\L}-x\right)\, .
\end{equation}
By Corollary \ref{CorConcentrationHamil}, we may replace $ \E_{\L,\beta,0}\left[\ex^{-\beta\Ham}
|\Nl= x\right]$ by $\ex^{b\beta x^2/(2\abs{\L})}$. We are now in position to apply \cite[Theorem 3]{martin1982laplace}, to conclude that
\begin{equation}
    \int_{\bB_T^{+}(\rhob)} \ex^{b\beta x^2/(2\abs{\L})}\frac{\beta\abs{\L}\c_d}{(\beta x)^{d/2+1}}\d\P_{\L,\beta,0}\left(N_\L=\rho\abs{\L}-x\right)\sim \frac{\beta\abs{\L}\c_d\sqrt{1+\frac{b}{ \bmu'(\rho_{\mathrm{c}}-\rho_{\mathrm{S}})}}\ex^{S_1\abs{\L}}}{(\beta \rhob\abs{\L})^{d/2+1}}\, .
\end{equation}
Indeed, the polynomial term $(\beta x)^{d/2+1}$ varies sufficiently slowly. This concludes the proof. \qed
\end{proofsect}
\begin{proofsect}{\textbf{Proof of Corollary \ref{freeenergycorolary}, free energy}}
By Lemma \ref{LemmaPartitionFunction},
\begin{equation}
    f^\hy(\beta,\rho)=-\lim_{N\to\infty}\frac{1}{\beta\abs{\L}}\log Z^{(\mathrm{Can, HYL})}_{\Lambda,\beta,\rho}=-\frac{S_1}{\beta}=-b\left(\rho_\mathrm{S}+\rho_\mathrm{e}\right)^2/2+\beta^{-1}I\left(\rho_\mathrm{c}-\rho_\mathrm{S}\right)\, .
\end{equation}
Using the definition of $\rhob$, we recognise the above as $-b\rhob^2/2+\beta^{-1}I\left(\rho-\rhob\right)$. The result now follows from the explicit form of $I$ given in Equation \eqref{EquationEcplicitI}.\qed
\end{proofsect}
\subsection{Computation of the limiting measure, supercritical}\label{subsectionlimiting} As usual, the computation of the partition function already reveals the limiting structure of the ensemble. Hence, using the results from the previous section together with the approximation techniques from \cite{vogel2021emergence}, the result emerges quickly.

There are two steps to the proof:
\begin{enumerate}
    \item The measure governing the long loops is converging to the intensity measure of the random interlacements.
    \item The remaining loops are governed by the loop soup with density $\bmu(\rho-\rhob)$. We employ a change of measure trick here, which is new compared to \cite{vogel2021emergence}.
\end{enumerate}
\textbf{Step 1:} let $f\ge 0$, $F(\eta)=\ex^{-\eta[f]}$ be a test function, as in Definition \ref{DefinitionLocConv}. Similar to the previous section, we can approximate
\begin{multline}\label{EquationEins}
     \E_{\L,\beta,0}\left[F(\eta)\ex^{-\beta\Ham},\bN_\L= \rho\right]=\E_{\L,\beta,0}\left[F(\eta)\ex^{-\beta\Ham},\bN_\L= \rho,\, \exists \omega\colon\ell(\omega)=\Nl\in\bB_T^{+}(\rhob)\right]\\
     +o\left(Z_{\L,\beta,\rho}^{(\mathrm{Can, HYL})}\right)\, .
\end{multline}
Here, recall that $\bB_T^{+}(\rhob)=\abs{\L}\rhob+\abs{\L}^{1/2}[-T,T]$, for $T>0$. On the event $\{\bN_\L= \rho\}$, we can rewrite
\begin{equation}
    \Ham(\eta)=-\frac{b\beta}{2\abs{\L}}\left(\rho\abs{\L}-\Ns\right)^2\, .
\end{equation}
Using the Mecke equation (see \cite{last2017lectures}) we expand leading order term in Equation \eqref{EquationEins} as
\begin{equation}\label{EquationZwei}
    \int \d \P_{\L,\beta,0}(\eta)\ex^{\frac{b\beta}{2\abs{\L}}\left(\rho\abs{\L}-\Ns_\L\right)^2}\int\d M_{\L,\beta,\rho}(\omega)F(\eta+\delta_\omega)\1_{A(\eta,\omega)}\, ,
\end{equation}
where
\begin{equation}
    A(\eta,\omega)=\{\ell(\omega)+N_\L(\eta)=\rho\abs{\L}\text{ and }N_\L(\eta)=\Ns(\eta)\text{ and }\ell(\omega)\in\bB_T^{+}(\rhob) \}\, .
\end{equation}
Write $a=\rho-\rhob$ and set $a^-=a\abs{\L}-T\abs{\L}^{1/2}$ and $a^+=a\abs{\L}+T\abs{\L}^{1/2}$. We rewrite Equation \eqref{EquationZwei} as
\begin{equation}
    \int \d \P_{\L,\beta,0}(\eta){\1_{\{N_\L(\eta)=\Ns(\eta)\in [a^-,a^+]\}}}\ex^{\frac{b\beta}{2\abs{\L}}\left(\rho\abs{\L}-N_\L\right)^2-\eta[f]}\int\d M_{\L,\beta,\rho}(\omega)\1_{\{\ell(\omega)=\rho\abs{\L}-N_\L(\eta)\}}\ex^{-f(\omega)}\, .
\end{equation}
Now by \cite{vogel2021emergence}, uniformly on the  event $\{N_\L(\eta)\in [a^-,a^+]\}$,
\begin{equation}
    \frac{M_{\L,\beta,\rho}\left[\1_{\{\ell(\omega)=\rho\abs{\L}-N_\L(\eta)\}}\ex^{-f}\right]}{\left(\frac{\beta\abs{\L}\c_d}{(\beta\rhob\abs{\L})^{1+d/2}}\right)}\sim \nu_{\rhob}[\ex^{-f}]\, .
\end{equation}
Here, $\nu_{\rhob}$ is the intensity measure of the (Brownian) random interlacements with density $\rhob>0$. While \cite{vogel2021emergence} was written for the case of the random walk and not the Brownian motion, this does not change the proof for the convergence to the interlacements. Indeed, the only time \cite{vogel2021emergence} uses properties of the random walk is in its Lemma 5.9. However, the continuum version of that lemma exists as \cite[Theorem 3]{uchiyama2018brownian}.

% Now note that I can expand on the event $A(\eta,\omega)$ the measure $\P_{\L,\beta,0}$ as
% \begin{equation}
%     \sum_{k=a^-}^{a^+}
% \end{equation}
% Write $a=\rho-\rhob$ and let $\tau=\beta\bmu(a)$ and set $r=\ex^{-I(a)}$. We can change the measure such that
% \begin{equation}
%     \E_{\L,\beta,0}[G]=\ex^{-\abs{\L}I(a)}{ \E_{\L,\beta,\bmu(a)}\left[G\ex^{-(\tau N_\L-a\tau)}\right]}\, ,
% \end{equation}
% for any function $G$.
% Now by \cite{vogel2021emergence}, we can do the approximation on the event $\{N_\L(\eta)\in [a^-,a^+]\}$
% \begin{equation}
%     \frac{M_{\L,\beta,\rho}\1\{\ell(\omega)=\rho\abs{\L}-N_\L(\eta)\}}{\left(\frac{\beta\abs{\L}\c_d}{(\beta\rhob\abs{\L})^{1+d/2}}\right)}\sim \nu_{\rhob}\, .
% \end{equation}
% .\\
To summarise the previous steps, we have now shown that
\begin{equation}\label{EquationSumIntlSetp}
    \frac{\E_{\L,\beta,0}\left[F(\eta)\ex^{-\beta\Ham},\bN_\L= \rho\right]}{Z_{\L,\beta,\rho}}\sim \frac{\E_{\L,\beta,0}\otimes\nu_{\rhob}\left[\1\{N_\L(\eta)=\Ns(\eta)\in [a^-,a^+]\}\ex^{\frac{b\beta}{2\abs{\L}}\left(\rho\abs{\L}-N_\L\right)^2}F(\eta+\delta_\omega)\right]}{\sqrt{1+\frac{b}{ \bmu'(\rho_{\mathrm{c}}-\rho_{\mathrm{S}})}}\ex^{S_1\abs{\L}}}\, .
\end{equation}
This concludes the first step.

\textbf{Step 2:} as done in the proof of Lemma \ref{LemmaTreash}, we can simplify the above to
\begin{equation}\label{EquationDrei}
    \frac{\E_{\L,\beta,0}\otimes\nu_{\rhob}\left[\1\{N_\L(\eta)\in [a^-,a^+]\}\ex^{\frac{b\beta}{2\abs{\L}}\left(\rho\abs{\L}-N_\L\right)^2}F(\eta+\delta_\omega)\right]}{\sqrt{1+\frac{b}{ \bmu'(\rho_{\mathrm{c}}-\rho_{\mathrm{S}})}}\ex^{S_1\abs{\L}}}\, ,
\end{equation}
using Lemma \ref{NoLongLoops}.
As $F(\eta+\delta_\omega)=F(\eta)F(\delta_\omega)$, we will omit the $\nu_{\rhob}$ part of the limiting process, to aid legibility. This means, we now examine
\begin{equation}
    \frac{\E_{\L,\beta,0}\left[\1\{N_\L(\eta)\in [a^-,a^+]\}\ex^{\frac{b\beta}{2\abs{\L}}\left(\rho\abs{\L}-N_\L\right)^2}F(\eta)\right]}{\sqrt{1+\frac{b}{ \bmu'(\rho_{\mathrm{c}}-\rho_{\mathrm{S}})}}\ex^{S_1\abs{\L}}}\, .
\end{equation}
Let $\tau=\beta\bmu(a)$. We can change the measure such that for any measurable $G$
\begin{equation}
    \E_{\L,\beta,0}[G]=\ex^{\abs{\L}\left[P(\bmu(a))-P(0)\right]}\E_{\L,\beta,\mu(a)}[G\ex^{-\beta\bmu(a)N_\L}]=\ex^{-\abs{\L}I(a)}{ \E_{\L,\beta,\bmu(a)}\left[G\ex^{-(\tau N_\L-a\tau\abs{\L})}\right]}\, . 
\end{equation}
Applying this to Equation \eqref{EquationDrei} leads to
\begin{equation}\label{Equation22220222}
    \frac{\E_{\L,\beta,\bmu(a)}\left[\1\{N_\L(\eta)\in [a^-,a^+]\}\ex^{\frac{b\beta}{2\abs{\L}}\left(\rho\abs{\L}-N_\L\right)^2-(\tau N_\L-a\tau\abs{\L})}F(\eta)\right]}{\ex^{\abs{\L}I(a)}\sqrt{1+\frac{b}{ \bmu'(\rho_{\mathrm{c}}-\rho_{\mathrm{S}})}}\ex^{S_1\abs{\L}}}\, .
\end{equation}
Note that the denominator can be simplified to
\begin{equation}
  \ex^{\abs{\L}I(a)}\sqrt{1+\frac{b}{ \bmu'(\rho_{\mathrm{c}}-\rho_{\mathrm{S}})}}\ex^{S_1\abs{\L}}=\ex^{\frac{b\rhob^2}{2}\abs{\L}}\sqrt{1+\frac{b}{\bmu'(\rho_{\mathrm{c}}-\rho_{\mathrm{S}})}}\, ,
\end{equation}
which we abbreviate by $\widetilde{Z}_\L$. 
We now expand the numerator in Equation \eqref{Equation22220222} as
\begin{equation}\label{EquationNewStage}
    \sum_{j=a^-}^{a^+}\ex^{\frac{b\beta}{2\abs{\L}}\left(\rho\abs{\L}-j\right)^2-(\tau j-a\tau\abs{\L})}\P\left(N_\L=j\right)\E_{\L,\beta,\bmu(a)}\left[F(\eta)\big|N_\L=j\right]\, .
\end{equation}
We have the following lemma:
\begin{lemma}\label{Lemmarefapprox}
Uniformly for $j\in [a^-,a^+]$, 
\begin{equation}
    \E_{\L,\beta,\bmu(a)}\left[F(\eta)\big|N_\L=j\right]\sim \E_{\L,\beta,\bmu(a)}\left[F(\eta)\right]\, .
\end{equation}
\end{lemma}
We give the proof of this lemma at the end of this subsection.

Using Lemma \ref{Lemmarefapprox}, we rewrite Equation \eqref{EquationNewStage} as
\begin{equation}
    \E_{\L,\beta,\bmu(a)}\left[F(\eta)\right]\sum_{j=a^-}^{a^+}\ex^{\frac{b\beta}{2\abs{\L}}\left(\rho\abs{\L}-j\right)^2-(\tau j-a\tau\abs{\L})}\P\left(N_\L=j\right)\sim \E_{\L,\beta,\bmu(a)}\left[F(\eta)\right]\ex^{\frac{b\rhob^2}{2}\abs{\L}}\sqrt{1+\frac{b}{\bmu'(\rho_{\mathrm{c}}-\rho_{\mathrm{S}})}}\, ,
\end{equation}
using \cite[Theorem 3]{martin1982laplace} again. This shows that
\begin{equation}
      \frac{\E_{\L,\beta,0}\left[\1\{N_\L(\eta)\in [a^-,a^+]\}\ex^{\frac{b}{2\abs{\L}}\left(\rho\abs{\L}-N_\L\right)^2}F(\eta)\right]}{\sqrt{1+\frac{b}{ \bmu'(\rho_{\mathrm{c}}-\rho_{\mathrm{S}})}}\ex^{S_1\abs{\L}}}\sim \E_{\L,\beta,\bmu(a)}\left[F(\eta)\right]\, ,
\end{equation}
and hence, using Equation \eqref{EquationSumIntlSetp}
\begin{equation}
    \frac{\E_{\L,\beta,0}\left[F(\eta)\ex^{-\beta\Ham},\bN_\L= \rho\right]}{Z_{\L,\beta,\rho}}\sim \E_{\L,\beta,\bmu(a)}\otimes \nu_{\rhob}\left[F(\eta+\delta_\omega)\right]\, ,
\end{equation}
in each box. Using the independence of the boxes and the superposition of Poisson processes, like in \cite{vogel2021emergence}, we get that
\begin{equation}
    \E_{n,\beta,0}\left[F\ex^{-\beta\Ham},\bN_\L= \rho\right]\sim \E_{\R^d,\beta,\bmu(\rho-\rhob)}\otimes\E_{\rhob}^\iota[F]\, .
\end{equation}
This concludes the proof of Theorem \ref{TheoremMain} for $\rho>\rho_{\mathrm{c}}$.
\begin{proofsect}{\textbf{Proof of Lemma \ref{Lemmarefapprox}}}
Choose $M=\abs{\L}^{1/10}$. Set $\bB_M=\bB_M(0)$ the ball of radius $M$, centred at the origin. We abbreviate $N_M=N_{\bB_M}$, $N_{\L\setminus M}=N_{\L\setminus \bB_M}$, $\E_{M,\beta,\bmu(a)}=\E_{\bB_M,\beta,\bmu(a)}$, and $\E_{\L\setminus M,\beta,\bmu(a)}=\E_{\L\setminus \bB_M,\beta,\bmu(a)}$. Set $R=M^2$. Using the superposition of Poisson point processes, we expand
\begin{multline}
    \E_{\L,\beta,\bmu(a)}\left[F,N_\L=j\right]=\E_{\L\setminus M,\beta,\bmu(a)}\otimes \E_{M,\beta,\bmu(a)}\left[F,N_{M}\le R,\,N_\L=j \right]\\
    +\E_{\L\setminus M,\beta,\bmu(a)}\otimes \E_{M,\beta,\bmu(a)}\left[F,N_{M}>R,\,N_\L=j \right]\, .
\end{multline}
As $N_{M}$ has exponential tails under $\P_{M,\beta,\bmu(a)}$ and $j-R$ is inside the range of Gnedenko's Local Limit Theorem (see \cite[Theorem 8.4.1]{bingham1989regular}), we have that
\begin{align}
    \E_{\L\setminus M,\beta,\bmu(a)}\otimes\E_{M,\beta,\bmu(a)}\left[F,N_{M}>R,\,N_\L=j \right]&\le C \P_{\L\setminus M,\beta,\bmu(a)}\otimes\P_{M,\beta,\bmu(a)}\left[N_{M}>R,\,N_\L=j \right]\nonumber\\
    &=o(1)\P_{\L\setminus M,\beta,\bmu(a)}\left[N_{\L\setminus M}=j \right]\nonumber\\
    &=o(1)\P_{\L,\beta,\bmu(a)}\left[N_{\L}=j \right],
\end{align}
as $M\to\infty$. We then expand
\begin{align}
        \E_{\L\setminus M,\beta,\bmu(a)}\otimes \E_{M,\beta,\bmu(a)}&\left[F,N_{M}\le R,\,N_\L=j \right]\nonumber\\
        &=\sum_{k=0}^R\E_{\L\setminus M,\beta,\bmu(a)}\otimes \E_{M,\beta,\bmu(a)}\left[F,N_{M}=k,\,N_{\L\setminus M}=j-k \right]\nonumber\\
        &=\sum_{k=0}^R\E_{M,\beta,\bmu(a)}\left[F,N_{M}=k\right]\E_{\L\setminus M,\beta,\bmu(a)}\left[F,\,N_{\L\setminus M}=j-k \right]\, ,
\end{align}
where we used the multiplication property of the test function $F(\eta)=\ex^{-\eta[f]}$. We want to replace $\E_{\L\setminus M,\beta,\bmu(a)}\left[F,\,N_{\L\setminus M}=j-k \right]$ by $\P_{\L\setminus M,\beta,\bmu(a)}\left[N_{\L\setminus M}=j-k \right]$. Note that 
\begin{equation}
    \{F\neq 1\}\subset\{\exists \omega\colon\supp(\omega)\cap\supp(f)\neq \emptyset \}\, .
\end{equation}
As the support of $f$ is compact, we may set $\supp(f)=[0,1]^d$ without loss of generality. We estimate
\begin{equation}
    \P_{\L\setminus M,\beta,\bmu(a)}\left(\exists \omega\in\eta\text{ such that }[0,1]^d\cap\omega \neq \emptyset\right)\le \sum_{j\ge 1}\frac{\ex^{\beta \bmu(a)j}}{j}\int_{\L\setminus \bB_M}\d x\, \P_{x,x}^{\beta j}\left(H_{[0,1]^d}<\beta j\right)\, ,
\end{equation}
where $H_{[0,1]^d}$ is the first hitting time of the unit cube, centred at the origin. Using that the distance between $x$ and $[0,1]^d$ is at least $R$, we bound (see for example \cite{bm})
\begin{equation}\label{hitting time estimateseq}
    \P_{x,x}^{\beta j}\left(H_{[0,1]^d}<\beta j\right)=\Ocal\left(\ex^{-cR/j^2}\right)\, . 
\end{equation}
This shows that the weight on the event $\{F\neq 1\}$ is of stretch-exponential order or less. However, $\P_{\L\setminus M,\beta,\bmu(a)}\left[N_{\L\setminus M}=j-k \right]$ is of polynomial order (again, by Gnedenko's Local Limit Theorem) and thus we can replace $F$ by $1$. This leads to
\begin{equation}
    \E_{\L\setminus M,\beta,\bmu(a)}\otimes \E_{M,\beta,\bmu(a)}\left[F,N_{M}\le R,N_\L=j \right]\sim\sum_{k=1}^R\E_{M,\beta,\bmu(a)}\left[F,N_{M}=k\right]\P_{\L\setminus M,\beta,\bmu(a)}\left[N_{\L\setminus M}=j-k \right],
\end{equation}
However, as $j-k$ is still in the CLT regime and $j-R=j(1+o(1))$, we get that
\begin{equation}
    \P_{\L\setminus M,\beta,\bmu(a)}\left[N_{\L\setminus M}=j-k \right]\sim \P_{\L\setminus M,\beta,\bmu(a)}\left[N_{\L\setminus M}=j \right]\sim \P_{\L,\beta,\bmu(a)}\left[N_{\L}=j \right]\, .
\end{equation}
This implies that 
\begin{equation}
    \E_{\L,\beta,\bmu(a)}\left[F(\eta)\big|N_\L=j\right]\sim\sum_{k=1}^R\E_{M,\beta,\bmu(a)}\left[F,N_{M}=k\right]\, .
\end{equation}
However, using the same reasoning as Equation \eqref{hitting time estimateseq}, we can see that
\begin{equation}
\begin{split}
    \E_{\L,\beta,\bmu(a)}\left[F(\eta)\big|N_\L=j\right]&\sim\sum_{k=1}^R\E_{M,\beta,\bmu(a)}\left[F,N_{M}=k\right]\\
    &=\E_{M,\beta,\bmu(a)}\left[F,N_{M}\le R\right]\sim \E_{\L,\beta,\bmu(a)}\left[F(\eta)\right]\, .
\end{split}
\end{equation}
This concludes the proof.\qed
\end{proofsect}
\subsection{The subcritical case}\label{subsectionbelowctricial} The case for $
\rho\le
\rho_{\mathrm{c}}$ is easier. Indeed, as the reference measure $\P_{\L,\beta,\rho}$ does not generate random interlacements by default, only the large deviation contribution from $\rhob$ will matter. Fix $\rho_o\le \rho_{\mathrm{c}}$ and define $S_1$ by
\begin{equation}
    S_1=\sup_{\rho
    \in(0,\rho_o)}\left\{\frac{b\beta \rho^2}{2}-I(\rho_o-\rho)\right\}\, .
\end{equation}
Furthermore, set $\rhob=\rhob(b,\rho_o,d)$ to be the value at which the supremum is achieved.
\begin{lemma}\label{LemmaMinimizerSubcritical} For every $b>0$
\begin{enumerate}
\item $\rhob=\rhob(\rho_o)\in [0,\rho_o)$ is well defined except for at most one $\rho_o$.
\item If $d\ge 5$ and $b\le 1/\brho'(0)$, then $\rhob=0$ for every $\rho_o\le \rho_{\mathrm{c}}$.
\item When $d=3,4$ or $b> 1/\brho'(0)$, $\rhob>0$ for $\rho_o$ sufficiently close to $\rho_{\mathrm{c}}$.
\item When $d=3,4$ or $b> 1/\brho'(0)$, $\rho_o\mapsto\rhob$ has a jump discontinuity. 
%     \item 
% $\rhob\in [0,\rho_o)$ is well defined. 
% \item 
% Given $b\ge\tfrac{\beta}{\brho'(0)}$ and $d\ge 5$, we have that $\rhob(b,\rho_o,d)=0$. Otherwise, 
% \begin{equation}\label{EquationJumpRhoBsub}
%     \lim_{\rho_o\uparrow \rho_{\mathrm{c}}}\rhob(b,\rho_o,d)>0\, .
% \end{equation}
% \item
% Whenever $ \lim_{\rho_o\uparrow \rho_{\mathrm{c}}}\rhob(b,\rho_o,d)>0$, there exists a $\rho_{\mathrm{c}}^\hy<\rho_{\mathrm{c}}$ such that
% \begin{equation}\label{EquationDichtonomy}
%     \rhob=0\text{ for }\rho_o<\rho_{\mathrm{c}}^\hy
%     \quad\text{and}\quad\rhob>0\text{ for }\rho_o>\rho_{\mathrm{c}}^\hy\, .
% \end{equation}
% \item $\lim_{\rho\downarrow\rho_{\mathrm{c}}^\hy}\rhob>0$.
\end{enumerate}
\end{lemma}
% \begin{figure}[h]
%     \centering
%     \begin{tikzpicture}[scale=1.2]
% \draw[->] (0,-2) -- (0,2)node[above]{$R(x)$};
%         \draw[->] (0,1) -- (4.2,1)node[below]{$x$};
%         %\draw[dashed] (4,-2) -- (4,2)node[right]{$\rho_o$};
%         \draw[thick] (0,0.5) to [out=0,in=90] (4,-1);
%         \draw[thick, dashed] (0,0.5) to [out=0,in=180] (0.75,-0.9) to [out=0,in=180] (2,0.25) to [out=0,in=90] (4,-1);
%         \draw[thick, dotted] (0,0.5) to [out=0,in=180] (0.75,0) to [out=0,in=180] (2.5,1.25) to [out=0,in=90] (4,-1);
%         \filldraw[black] (0,-1)  circle (1pt)node[left]{$R(\rho_o)$};
%         \draw[dashed] (0,-1) -- (4,-1)node[right]{$\rho_o$};
%         %\filldraw[black] (2,1)  circle (1pt)node[below]{$\rhob$};

% \draw[->] (6,-2) -- (6,2)node[above]{$R'(x)$};
%         \draw[->] (6,1) -- (10.2,1)node[below]{$x$};
%         \filldraw[black] (9.5,1)  circle (1pt)node[below]{$\rho_o$};
%         \draw[thick] (6,0.5) to [out=0,in=90] (9.5,-2);
%         \draw[thick, dashed] (6,0.5) to [out=0,in=180] (6.95,0.75) to [out=0,in=180] (8,0.25) to [out=0,in=90] (9.5,-2);
%         \draw[thick, dotted] (6,0.5) to [out=0,in=180] (7.5,1.75) to [out=0,in=90](9.5,-2);
%       \draw[dashed] (9.5,2) -- (9.5,-2);

% \end{tikzpicture}
% \end{figure}

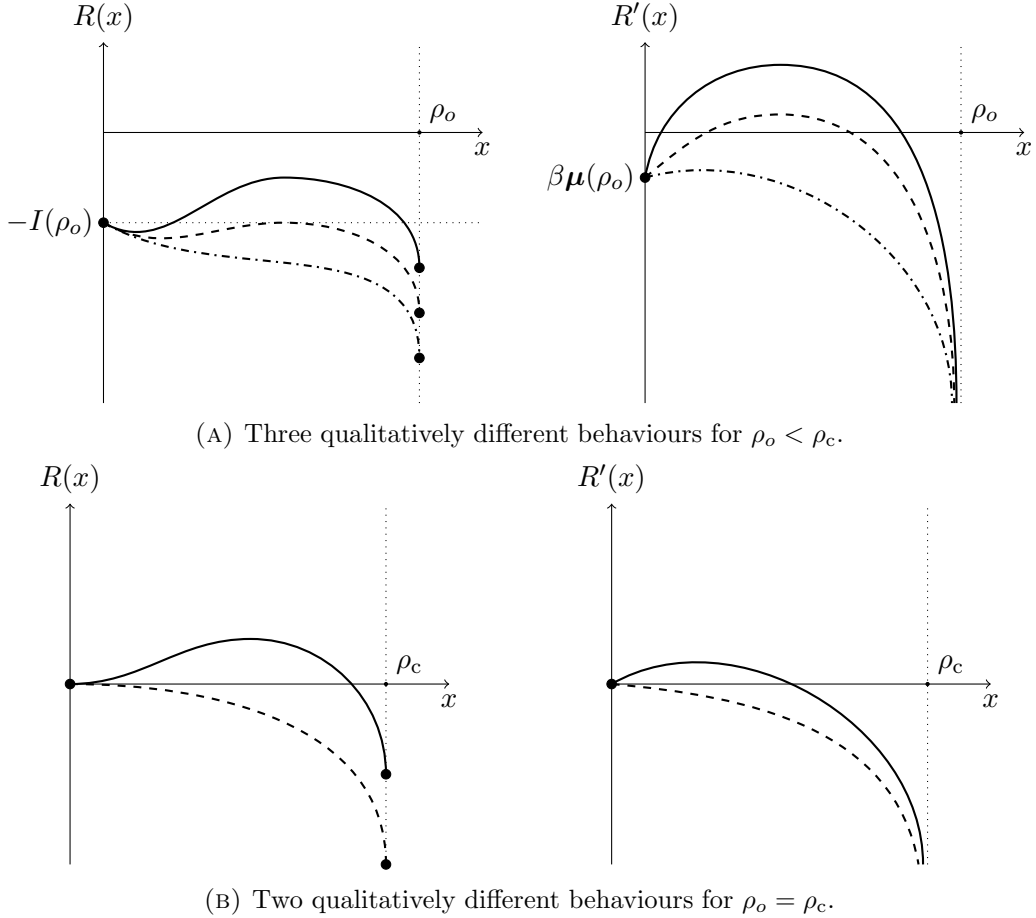
\begin{figure}[ht]
    \centering
    \begin{subfigure}[t]{\textwidth}
    \centering
    \begin{tikzpicture}[scale=1.2]
\draw[->] (0,-2) -- (0,2)node[above]{$R(x)$};
        \draw[->] (0,1) -- (4.2,1)node[below]{$x$};
        \draw[thick,dashdotted] (0,0) to [out=330,in=90] (3.5,-1.5);
        \draw[thick,dashed] (0,0) to [out=330,in=180] (2,0) to [out=0,in=90] (3.5,-1);
        \draw[thick] (0,0) to [out=330,in=180] (2,0.5) to [out=0,in=90] (3.5,-0.5);
        \draw[dotted] (0,0) -- (4.2,0);
        \draw[dotted] (3.5,-2) -- (3.5,2);
        \filldraw (3.5,-1.5) circle (1.5pt);
        \filldraw (3.5,-1) circle (1.5pt);
        \filldraw (3.5,-0.5) circle (1.5pt);
        \filldraw (0,0) circle (1.5pt)node[left]{$-I(\rho_o)$};
        \filldraw (3.5,1) circle (0.5pt)node[above right]{$\rho_o$};

\draw[->] (6,-2) -- (6,2)node[above]{$R'(x)$};
        \draw[->] (6,1) -- (10.2,1)node[below]{$x$};
        \draw[thick,dashdotted] (6,0.5) to [out=15,in=90] (9.4,-2);
        \draw[thick,dashed] (6,0.5) to [out=40,in=180] (7.5,1.2) to [out=0,in=90] (9.425,-2);
        \draw[thick] (6,0.5) to [out=80,in=180] (7.5,1.75) to [out=0,in=90] (9.45,-2);
        \draw[dotted] (9.5,-2) -- (9.5,2);
        \filldraw (6,0.5) circle (1.5pt)node[left]{$\beta\bmu(\rho_o)$};
        \filldraw (9.5,1) circle (0.5pt)node[above right]{$\rho_o$};

\end{tikzpicture}
\caption{Three qualitatively different behaviours for $\rho_o<\rho_{\mathrm{c}}$.}
    \end{subfigure}
    \begin{subfigure}[t]{\textwidth}
    \centering
    \begin{tikzpicture}[scale=1.2]
\draw[->] (0,-2) -- (0,2)node[above]{$R(x)$};
        \draw[->] (0,0) -- (4.2,0)node[below]{$x$};
        \draw[thick,dashed] (0,0) to [out=0,in=90] (3.5,-2);
        \draw[thick] (0,0) to [out=0,in=180] (2,0.5) to [out=0,in=90] (3.5,-1);
        \draw[dotted] (3.5,-2) -- (3.5,2);
        \filldraw (3.5,-2) circle (1.5pt);
        \filldraw (3.5,-1) circle (1.5pt);
        \filldraw (0,0) circle (1.5pt);
        \filldraw (3.5,0) circle (0.5pt)node[above right]{$\rho_{\mathrm{c}}$};

\draw[->] (6,-2) -- (6,2)node[above]{$R'(x)$};
        \draw[->] (6,0) -- (10.2,0)node[below]{$x$};
        \draw[thick,dashed] (6,0) to [out=355,in=100] (9.4,-2);
        \draw[thick] (6,0) to [out=30,in=90] (9.45,-2);
        \draw[dotted] (9.5,-2) -- (9.5,2);
        \filldraw (6,0) circle (1.5pt);
        \filldraw (9.5,0) circle (0.5pt)node[above right]{$\rho_{\mathrm{c}}$};

\end{tikzpicture}
\caption{Two qualitatively different behaviours for $\rho_o=\rho_{\mathrm{c}}$.}
    \end{subfigure}
    \caption{The function $x\mapsto R(x)$ and its derivative. Different possibilities (depending on the values of $\rho_o,\beta,b$) are drawn in different styles.}
    \label{fig:my_label}
\end{figure}

 \begin{proofsect}{\textbf{Proof of Lemma \ref{LemmaMinimizerSubcritical}}}
 Set $R(x)=R_{\rho_o}(x)=b\beta x^2/2-I(\rho_o-x)$. Note, as before $R'(x)=b\beta x+\beta\bmu(\rho_o-x)$ and $R''(x)=\beta b-\beta \bmu'(\rho_o-x)$. Recall that Figure \ref{figurethermodynamic functions} sketches the functions $\brho,\,\bmu$ as well as their derivatives.
 \begin{enumerate}
     \item Now in contrast to the supercritical case, whilst $\rho_o<\rho_{\mathrm{c}}$ we have that $R(x)$ is decreasing for small $x$. Indeed, $R'(x)$ is strictly negative in a neighbourhood around the origin. By looking at the second derivative, which is monotonously decreasing, we distinguish three cases: $x\mapsto R'(x)$ has either none, one or two zeros. If $x\mapsto R'(x)$ has none or one zero, the maximum value of $R(x)$ is attained at the origin. If $x\mapsto R'(x)$ has two zeros, the maximum value of $x\mapsto R(x)$ will be attained at the origin or at the rightmost zero of $x\mapsto R'(x)$.
     
     We first show that for $\rho_o>0$ sufficiently small, $\rhob=0$. Note that for sufficiently small $\rho_o>0$
     \begin{equation}
         \sup_{x\in (0,\rho_o)}\{R''(x)\}=\sup_{x\in (0,\rho_o)}\{\beta b-\beta\bmu'(\rho_o-x)\}<0\quad\text{as}\quad \bmu'(y)\to \infty,\,\text{ when }y\downarrow 0\, .
     \end{equation}
     This shows that for $\rho_o$ small enough, $R'(x)$ is decreasing. As $R'(0)<0$, the supremum of $R(x)$ is attained at the origin. Thus, for small $\rho_o$, $\rhob$ is well defined.\\
     For larger $\rho_o$ we have to prove that the supremum is attained at only one point. For this, note that the change of $x\mapsto R(x)$ with respect to $\rho_o$ is positive, i.e., $\frac{\d}{\d \rho_o}R_{\rho_o}(x)>0$ for every $x$ in the domain. Thus, $\rhob$ is either increasing or jumps to zero, as $\rho_o$ increases. Write $x_o=x_o(\rho_o)$ for the location of the second zero, if it exists. However,
     \begin{equation}\label{euqationdsd}
         \frac{\d}{\d \rho_o}R_{\rho_o}(0)=-\beta\bmu(\rho_o)<-\beta\bmu(\rho_o-\e)=\frac{\d}{\d \rho_o}R_{\rho_o}(\e)\, ,
     \end{equation}
     for any $\e>0$. This implies that $R(x)$ grows faster away from the origin. Hence, as soon as $\rhob(\rho_o)>0$, it has to be greater than zero for any larger $\tilde\rho_o>\rho_o$. Equation \eqref{euqationdsd} also shows that there can be at most one $\rho_o$ for which $R(x_o)=R(0)$. At this point, $\rhob$ is not well defined.
    %  Notice that
    %  \begin{equation}
    %      r_{\rho_o+\delta}(x)=r_{\rho_o}(x)+I(\rho_o-x)-I(\rho_o-x+\delta)=r_{\rho_o}(x)-\delta\beta\bmu(\rho_o-x)(1+o(1))\, .
    %  \end{equation}
    %  Thus, $\rhob$ is increasing in $\rho_o$.
     \item Here, we argue similar as before: if $b<1/\brho'(0)$ and $d\ge 5$, the second derivative is always negative. Furthermore, the first derivative is negative. Thus, in that case, $R(x)$ is always decreasing and the minimum is attained at the origin.
     \item If $b>1/\brho'(0)$ or $d=3,4$, note that 
     \begin{equation}
         \exists \epsilon>0 \exists\delta>0\exists\gamma>0\forall \rho_o\in(\rho_{\mathrm{c}}-\e,\rho_{\mathrm{c}})\forall x\in [0,\delta\wedge \rho_o)\colon R''_{\rho_o}(x)>\gamma\, .
     \end{equation}
     This is because $\bmu'(x)$ decays either faster than linearly or with a linear coefficient larger than $b$. This gives us the bound
     \begin{equation}
         R'(x)\ge \gamma x+\beta\bmu(\rho_o)\, .
     \end{equation}
     Hence, the interval around the origin on which $R'$ is negative shrinks to zero, as $\rho_o$ increases to $\rho_\mathrm{c}$. Thus, by letting $\rho_o$ tend to $\rho_{\mathrm{c}}$, we see that in this limit, $\rhob$ has to be positive for all $\rho_o$ with $\rho_{\mathrm{c}}-\rho_o>\e_\mathrm{c}$ for some $\e_\mathrm{c}>0$.
     \item Define
     \begin{equation}\label{EqdefrhocHY}
         \rho_{\mathrm{c}}^\hy=\sup\{\rho_o\colon \rhob(\rho_o)=0\}\, .
     \end{equation}
     Now for $b>1/\brho'(0)$ or $d=3,4$
     \begin{equation}
         0<\rho_{\mathrm{c}}^\hy<\rho_{\mathrm{c}}\, .
     \end{equation}
     Suppose that $\rho_o\mapsto\rhob$ was continuous at $\rho_{\mathrm{c}}^\hy$. Recall that for $\rho_o>\rho_{\mathrm{c}}^\hy$, the supremum is attained at the second zero of $x\mapsto R'(x)$. However, as $\rho_{\mathrm{c}}^\hy<\rho_{\mathrm{c}}$, $R'(x)$ is bounded away uniformly (in both $x,\rho_o$) from zero. This implies that the second zero of $R'(x)$ is bounded away from the origin uniformly in $\rho_o$. This is a contradiction, as $\rhob$ is either zero or equal to that second zero. This concludes the proof.\qed
 \end{enumerate}

 \end{proofsect}
%  Note that we can determine further properties of $\rhob$. The implicit function theorem gives $\partial_{\rho_o}\rhob=\tfrac{\beta \bmu'(\rho_o-\rhob)}{\beta \bmu'(\rho_o-\rhob)-b}$.\\
 Given the previous lemma,  we can follow the same steps as in the supercritical case. The partition function is asymptotically equal to
\begin{equation}
    Z^{(\mathrm{Can, HYL})}_{\Lambda,\beta,\rho}\sim \sqrt{1+\frac{b}{ \bmu'(\rho_o-\rhob)}}\exp\left(S_1\abs{\L}\right)\frac{\beta\abs{\L}\c_d}{(\beta \rhob\abs{\L})^{d/2+1}}\, 
\end{equation}
Similarly, one can then show the convergence of the conditional measure. This concludes the proof of Theorem \ref{TheoremMain} for the case $\rho<\rho_{\mathrm{c}}$. The case $\rho=\rho_{\mathrm{c}}$ is only relevant for $d=3,4$ or $d\ge 5$ with $b>\beta \brho'(0)^{-1}$. However, in that case the proof of Lemma \ref{LemmaMinimizerSubcritical} does not change as $\rho_{\mathrm{c}}^\hy<\rho_{\mathrm{c}}$.\qed
\subsection{The grand-canonical ensemble}\label{subsectionGC}
Fix $\mu\in\R$. We now examine the properties of the Bose gas in terms of the parameters $\beta>0$ and $\mu$. We expand the partition function
\begin{equation}
    \E_{\L,\beta,\mu}\left[\ex^{-\beta\Ham-\beta\Hpmf}\right]=\E_{\L,\beta,\mu}\left[\ex^{-\beta\Hpmf}\right]\sum_{j=0}^\infty \E_{\L,\beta,\mu}\left[\ex^{-\beta\Ham}|N_\L=j\right]\P_{\L,\beta,\mu}^\pmf(N_\L=j)\, .
\end{equation}
Here, we used that $\E_{\L,\beta,\mu}\left[\ex^{-\beta\Ham}|N_\L=j\right]=\E_{\L,\beta,\mu}^\pmf\left[\ex^{-\beta\Ham}|N_\L=j\right]$. Note that by Varadhan's Lemma, $\P_{\L,\beta,\mu}^\pmf$ satisfies a large deviation principle with rate function
\begin{equation}\label{EquationrhoGC}
    I_\mu(x)+\frac{\beta a x^2}{2}-\inf_x \left\{I_\mu(x)+\frac{\beta a x^2}{2}\right\}=I(x)+\frac{\beta a x^2}{2}-\beta \mu x-\inf_x\left\{ I(x)+\frac{\beta a x^2}{2}-\beta \mu x\right\}\, .
\end{equation}
Define $\rho^{\mathrm{GC}}=\rho^{\mathrm{GC}}(\mu)$ to be the maximiser of the function $J$ with
\begin{equation}
    J(\rho)=\frac{\beta b\rhob^2}{2}-I(\rho-\rhob)-\frac{\beta a\rho^2}{2}+\beta \mu \rho-I(\rho)\, .
\end{equation}
Here, $\rhob=\rhob(b,\rho,d)$ maximises the function $x\mapsto \beta bx^2/2-I(\rho-x)$ on $(0,\rho)$. For $\rho=\rho_{\mathrm{c}}^\hy$, this may not be well defined as there could exist two maximizers. However, $b\rhob^2/2-I(\rho-\rhob)$ is well defined at this point. The next lemma gives the properties of $\rho^{\mathrm{GC}}$.
\begin{lemma}\label{LemmaGCMinimizer} The following statements hold true for any $b>0$.
\begin{enumerate}
    \item 
For any $\mu\in\R$, $\rho^{\mathrm{GC}}$ is well defined and $\rho^{\mathrm{GC}}\in (0,\infty)$. 
\item
The map $\mu\mapsto\rho^{\mathrm{GC}}(\mu)$ is strictly increasing.
\item 
$\lim_{\mu\to\infty}\rho^{\mathrm{GC}}=\infty$ and $\lim_{\mu\to-\infty}\rho^{\mathrm{GC}}=0$.
\end{enumerate}
\end{lemma}
\begin{proofsect}{\textbf{Proof of Lemma \ref{LemmaGCMinimizer}}}
Recall that the thermodynamic functions are sketched in Figure \ref{figurethermodynamic functions}.
\begin{enumerate}
    \item 
    It is clear that the maximizer(s) have to be bounded away from zero and $+\infty$: as $\rhob<\rho$ and $a>b$, $J(\rho)$ diverges to $-\infty$ as $\rho$ gets large. For $\rho\neq \rho_{\mathrm{c}}^\hy$
\begin{equation}
    \frac{J'(\rho)}{\beta}=b\rhob \rhob'-\bmu(\rho-\rhob)(1-\rhob')-a\rho+\mu-\bmu(\rho)\, .
\end{equation}
Set $\beta=1$, to shorten notation. Recall that for $\rho<\rho_{\mathrm{c}}^\hy$, $\rhob=0$. From there, it follows
\begin{equation}\label{EquationIncrease}
    J'(\rho)=\begin{cases}
    \mu-a\rho-2\bmu(\rho)&\text{ if }\rho<\rho_{\mathrm{c}}^\hy\, ,\\
    \mu-a\rho-\bmu(\rho-\rhob)-\bmu(\rho)&\text{ if }\rho>\rho_{\mathrm{c}}^\hy\, .
    \end{cases}
\end{equation}
Indeed, if $\rho>\rho_{\mathrm{c}}^\hy$, this means that $\rhob>0$. As $\rhob$ is defined as the maximizer of the differentiable function $x\mapsto bx^2/2-I(\rho-x)$, it holds
\begin{equation}
    b\rhob+\bmu(\rho-\rhob)=0\, .
\end{equation}
From Equation \eqref{EquationIncrease}, we can see that $J$ increases in a small neighbourhood of the origin. This shows that the maximiser(s) have to be contained in $(0,\infty)$. Now we show that there can only be one maximiser. It holds
{ \renewcommand{\arraystretch}{1.5}
 \begin{align}
   J''(\rho)&=\left\{\begin{array}{l}
    -a-2\bmu'(\rho)\\
    -a-\bmu'(\rho-\rhob)(1-\rhob')-\bmu'(\rho)
   \end{array}\right.
  &\begin{array}{l}
   \text{ if }\rho<\rho_\mathrm{c}^\hy\, ,\\
   \text{ if }\rho>\rho_\mathrm{c}^\hy\, ,\
  \end{array}
  \nonumber\\
   &=\left\{\begin{array}{l}
    -a-2\bmu'(\rho) \\
    -(a-b)-\bmu'(\rho) 
   \end{array}\right.
  &\begin{array}{l}
   \text{ if }\rho<\rho_\mathrm{c}^\hy\, ,\\
   \text{ if }\rho>\rho_\mathrm{c}^\hy\, .\
  \end{array}
 \end{align}
}
Indeed, by the implicit function theorem
\begin{equation}
    \rhob'=\frac{\d }{\d \rho}\rhob=\frac{\bmu'(\rho-\rhob)+b}{\bmu'(\rho-\rhob)}>0\quad\text{and thus}\quad 1-\rhob'=\frac{-b}{\bmu'(\rho-\rhob)}<0\, .
\end{equation}
Recall that $a>b$. As $J$ is continuous, with first derivative positive in a neighbourhood around the origin and the second derivative strictly negative, it follows that it attains its maximum at a single point.
% We first show that $\rho^{\mathrm{GC}}<\infty$. For this, note that $\rhob<\rho$ and that $b<a$. Both $I(\rho-\rhob)$ and $I(\rho)$ are bounded while $\mu \rho$ is of linear speed. It follows that $J(\rho)$ diverges to $-\infty$ at quadratic speed and therefore has a well defined maximiser. To show that $\rho^{\mathrm{GC}}>0$, observe that $J'(\rho)=-2\beta\bmu(\rho)+\mu-a\rho$ on the interval $(0,\rho_{\mathrm{c}}^\hy)$. This implies that $J'(\rho)$ is positive in a neighbourhood of the origin. This concludes the proof of the first claim.
\item
To prove the monotonicity of $\rho^{\mathrm{GC}}$, we calculate its derivative with respect to $\mu$. As $J'$ is increasing in a neighbourhood around the origin, we can assume that $\rho^{\mathrm{GC}}$ is bigger than zero. Let us assume that $\rho^{\mathrm{GC}}\neq \rho_{\mathrm{c}}^\hy$. In that case, $\rho^{\mathrm{GC}}$ satisfies $J'(\rho^{\mathrm{GC}})=0$. If $\rho^{\mathrm{GC}}<\rho_{\mathrm{c}}^\hy$, this implies by the implicit function theorem that $\tfrac{\d}{\d \mu}\rho^{\mathrm{GC}}=\left[2\bmu'(\rho)+a\right]^{-1}>0$. If $\rho^{\mathrm{GC}}>\rho_{\mathrm{c}}^\hy$,
\begin{equation}
   \frac{\d}{\d \mu}\rho^{\mathrm{GC}}=\frac{1}{(a-b)+\bmu'(\rho)} >0\, .
\end{equation}
Thus, the positivity of $\tfrac{\d}{\d \mu}\rho^{\mathrm{GC}}$ follows. 
Finally, if $\rho^{\mathrm{GC}}=\rho_{\mathrm{c}}^\hy$, it is straight forward to show the claim.
\item
For this claim, note that if $\mu\to-\infty$, this implies that any $\rho$ solving the equation $\mu=2\beta\bmu(\rho)+a\rho$ goes to zero. Furthermore, we can see that for the case $\mu\to\infty$, $\rho$ diverges to $+\infty$. This can be seen by looking at the area where the first derivative is positive.
\qed
\end{enumerate}
\end{proofsect}
Next, we calculate the partition function.
\begin{cor}\label{LemmaGCPartitionFunction}
As $\L\uparrow\R^d$, it holds that
\begin{equation}
    \E_{\L,\beta,\mu}\left[\ex^{-\beta \Hpmf}\right]\sum_{j=0}^\infty \E_{\L,\beta,\mu}\left[\ex^{-\beta \Ham}|N_\L=j\right]\P_{\L,\beta,\mu}^\pmf(N_\L=j)\sim Z_{\L,\beta,\rho^{\mathrm{GC}}}^{(\mathrm{Can, HYL})}\, .
\end{equation}
\end{cor}
\begin{proofsect}{\textbf{Proof of Corollary \ref{LemmaGCPartitionFunction}}}
Note that analogous to Lemma \ref{LemmaTreash}
\begin{equation}
   \sum_{j=0}^\infty \E_{\L,\beta,\mu}\left[\ex^{-\beta \Ham}|N_\L=j\right]\P_{\L,\beta,\mu}^\pmf(N_\L=j)\sim \sum_{j=\rho^{\mathrm{GC}}\abs{\L}-T\abs{\L}^{1/2}}^{\rho^{\mathrm{GC}}\abs{\L}+T\abs{\L}^{1/2}} \E_{\L,\beta,\mu}\left[\ex^{-\beta \Ham}|N_\L=j\right]\P_{\L,\beta,\mu}^\pmf(N_\L=j)\, .
\end{equation}
For $j\in\abs{\L}\rho^{\mathrm{GC}}+\abs{\L}^{1/2}[-T,T]$, we may employ the expansion from Lemma \ref{LemmaPartitionFunction}. This concludes the proof.\qed
\end{proofsect}
Having calculated the partition function, we can compute the limiting measure with no difficulties. The steps are the same as in Section \ref{subsectionlimiting} and we leave the details to the reader.
\begin{proofsect}{\textbf{Proof of Theorem \ref{TheoremGrandCanonical}}}
Using the tools from Section \ref{subsectionlimiting} and Corollary \ref{LemmaGCPartitionFunction}, we can show that
\begin{equation}
    \frac{1}{Z_{\L,\beta,\rho^{\mathrm{GC}}}}\sum_{j=0}^\infty \E_{\L,\beta,\mu}\left[F(\eta)\ex^{-\beta \Ham}|N_\L=j\right]\P_{\L,\beta,\mu}^\pmf(N_\L=j)\sim \E_{\L,\beta,\bmu(a)}\otimes\nu_{\rhob}\left[F(\eta+\delta_\omega)\right]\, ,
\end{equation}
where $\rhob=\rhob(b,\rho^{\mathrm{GC}},d)$ and $a=\rho^{\mathrm{GC}}-\rhob$. Proceeding analogous to Section \ref{subsectionlimiting} concludes the proof.\qed
\end{proofsect}
\begin{proofsect}{\textbf{Proof of Corollary \ref{freeenergycorolary}, pressure}}
The proof is an immediate consequence of the above result.\qed
\end{proofsect}
%%%%%%%%%%%%%%%%%%%%%%%%%%%%%%%%%%%%%%%%%%%%%%%%%%%%%%%%%%%%%%%%%%%%%%%%%%%%%%%%%%%%%%%%%%%%%%%%%%%%%%%%%%%%%%%%%%%%%%%%%%%%%%%%%%

\subsection{Proof of Proposition \ref{propositionfullvspartial}, pressure comparison}\label{SubsectionPressureDifference}

Recall that Proposition~\ref{propositionfullvspartial} stated that the thermodynamic pressures produced by the partial HYL Hamiltonian $\Ham_\Lambda$ and the full HYL Hamiltonian $\widetilde{\Ham}_\Lambda$ (where we fix $q_\Lambda=1$) are indeed different.

We first introduce some notation. Recall the thermodynamic pressures $P^{\hy}\left(\beta,\mu\right)$ and $\widetilde{P}^{\hy}\left(\beta,\mu\right)$ given in \eqref{eqn:thermodynamicPressures}. We now use $P^{\hy}_\Lambda\left(\beta,\mu\right)$ and $\widetilde{P}^{\hy}_\Lambda\left(\beta,\mu\right)$ to denote the respective finite-volume pressures:
\begin{equation}
\label{eqn:finitevolumePressures}
    P^{\hy}_\Lambda\left(\beta,\mu\right) = \frac{1}{\beta \abs*{\Lambda}}\E_{\Lambda,\beta,0}\left[\ex^{\beta\mu N_\Lambda-\beta \Ham_\Lambda}\right],\qquad \widetilde{P}^{\hy}_\Lambda\left(\beta,\mu\right) = \frac{1}{\beta \abs*{\Lambda}}\E_{\Lambda,\beta,0}\left[\ex^{\beta\mu N_\Lambda-\beta \widetilde{\Ham}_\Lambda}\right],
\end{equation}
so $P^{\hy}_\Lambda\left(\beta,\mu\right)\to P^{\hy}\left(\beta,\mu\right)$ and $\widetilde{P}^{\hy}_\Lambda\left(\beta,\mu\right)\to \widetilde{P}^{\hy}\left(\beta,\mu\right)$. Also recall how the grand-canonical measure with HYL interaction, $\GCanonHY{\Lambda}{\mu}$, is defined using the Hamiltonian $\Ham_\Lambda$ (see \eqref{eqn:GCHYLradonnikodymDerivative}). We now use $\E^\hy_{\Lambda,\beta,\mu}$ to denote the expectation with respect to this measure, and respectively define $\widetilde{\P}^\hy_{\Lambda,\beta,\mu}$ and $\widetilde{\E}^\hy_{\Lambda,\beta,\mu}$ with $\widetilde{\Ham}_\Lambda$ replacing $\Ham_\Lambda$. It will also be convenient to define
\begin{equation}
    \lambda_\L(\eta) = \left\{\lambda_{\L,j}(\eta)\right\}_{j\geq 1} := \left\{\frac{1}{\abs*{\L}}\#\left\{\omega\in\eta_\Lambda:\ell(\omega) = j\right\}\right\}_{j\geq 1}
\end{equation}
in $\ell^1\left(\R_{\geq0}\right)$ for each locally finite counting measure $\eta$ on $E$. Each entry $\lambda_{\L,j}(\eta)$ then gives the density of loops of length $j$ rooted in $\L$.

\begin{proof}[Proof of Proposition~\ref{propositionfullvspartial}]
Let $\Ham^*_\L = \Ham_\Lambda - \widetilde{\Ham}_\Lambda$ be the difference in the Hamiltonian densities. So given a locally finite counting measure $\eta$ on $E$, we have
\begin{equation}
    \Ham^*_\L\left(\eta\right) = \frac{b}{2\abs*{\L}}\sum^{q_\L-1}_{j=1}j^2\#\left\{\omega\in\eta_\Lambda:\ell(\omega) = j\right\}^2 = \abs*{\L}\frac{b}{2}\sum^{q_\L-1}_{j=1}j^2\lambda_{\L,j}(\eta)^2.
\end{equation}
By using a change of measure from the non-interacting distribution to the full loop HYL distribution, $\E_{\L,\beta,0}\left[\ex^{\beta\mu N_\Lambda-\beta \Ham_\L}\right] = \E_{\L,\beta,0}\left[\ex^{\beta\mu N_\Lambda-\beta \widetilde{\Ham}_\L}\right]\widetilde{\E}^\hy_{\L,\beta,0}\left[\ex^{-\beta \Ham^{*}_\L}\right]$. Therefore
\begin{equation}
    P^{\hy}_\Lambda\left(\beta,\mu\right) = \widetilde{P}^{\hy}_\Lambda\left(\beta,\mu\right) + \frac{1}{\beta \abs{\L}}\log \widetilde{\E}^\hy_{\L,\beta,0}\left[\ex^{-\beta \Ham^{*}_\L}\right].
\end{equation}
Since $\Ham^*_\L$ is non-negative, we are reassured that the logarithm term in non-positive. We can further bound $\Ham^*_\L\left(\eta\right) \geq \abs*{\L}\frac{b}{2}\lambda^2_{\L,1}(\eta)$, and so
\begin{equation}
\label{eqn:PressureDifference}
    P^{\hy}_\Lambda\left(\beta,\mu\right) \leq \widetilde{P}^{\hy}_\Lambda\left(\beta,\mu\right) + \frac{1}{\beta \abs{\L}}\log \widetilde{\E}^\hy_{\L,\beta,0}\left[\ex^{-\abs{\L}\beta \frac{b}{2}\lambda_{\L,1}^2}\right].
\end{equation}

A large deviation principle for $\lambda_\Lambda$ under the full HYL model was derived in \cite[Theorem~1.6]{adams2021large}. This principle holds with respect to the $\ell^1$-topology, has rate $\abs{\L}$, and has rate function given by
\begin{multline}
\label{eqn:fullloopHYLratefunction}
    \mathcal{I}\left(x\right) = \sum^\infty_{j=1}x_j\left(\log\frac{jx_j}{\p_{\beta j}(0)}-1\right) -\mu\beta D\left(x\right) + \frac{a\beta}{2}D\left(x\right)^2 - \frac{b\beta}{2}\sum^\infty_{j=1}j^2x^2_j - \frac{\beta}{2\left(a-b\right)}\left(\mu-aD\left(x\right)\right)^2_+ \\- P\left(\beta,0\right) + \widetilde{P}^\hy\left(\beta,\mu\right), 
\end{multline}
where $D\left(x\right) = \sum^\infty_{j=1}jx_j\in\left[0,+\infty\right]$ and $P\left(\beta,0\right)$ is the thermodynamic pressure of the non-interacting model with $\mu=0$. If $D\left(x\right)=+\infty$ then we set $\mathcal{I}\left(x\right)=+\infty$. Since $x\mapsto \beta\frac{b}{2}x^2_1$ is continuous with respect to the $\ell^1$-topology and is bounded from below by $0$, Varadhan's Lemma gives the existence and equality of the following limit:
\begin{equation}
\label{eqn:UsedVaradhan}
    \lim_{\abs*{\L}\to\infty}\frac{1}{\beta \abs{\L}}\log \widetilde{\E}^\hy_{\L,\beta,0}\left[\ex^{-\abs{\L}\beta \frac{b}{2}\lambda_{\L,1}^2}\right] = -\frac{1}{\beta}\inf_{x\in\ell^1\left(\R_+\right)}\left\{\mathcal{I}\left(x\right) + \frac{b\beta}{2}x^2_1\right\}.
\end{equation}

Let us first remark on the existence of global minimisers of $\mathcal{I}$. From the expression \eqref{eqn:fullloopHYLratefunction}, note that there exist $C_1 = C_1\left(\mu,\beta,a,b\right)<0$ and $C_2 = C_2\left(\mu,\beta,a,b\right) >0$ such that $\mathcal{I}\left(x\right) \geq C_1 + C_2D\left(x\right)^2 \geq C_1 + C_2\norm*{x}_1^2$. Therefore the level sets of $\mathcal{I}$ are $\ell^1$-bounded. Along with the lower semicontinuity of $\mathcal{I}$, this implies that there exists at least one global minimiser of $\mathcal{I}$.

We now show that any such global minimiser is in the set $\left\{x_1 \geq \varepsilon \right\}$ for some $\varepsilon>0$. Taking the $x_1$-partial derivative of $\mathcal{I}$, we find
\begin{equation}
    \frac{\partial \mathcal{I}}{\partial x_1}\left(x\right) = \log \frac{x_1}{\p_{\beta }(0)} - b\beta x_1 - \beta \left(\mu - a D\left(x\right)\right)
	\begin{Bmatrix}
	1&:aD\left(x\right) \geq \mu\\
	-\frac{b}{a-b}&:aD\left(x\right) \leq \mu
	\end{Bmatrix}.
\end{equation}
The bound $\mathcal{I}\left(x\right) \geq C_1 + C_2D\left(x\right)^2$ implies that the value of $D\left(x\right)$ at any global minimiser of $\mathcal{I}$ cannot be arbitrarily large for given parameters $\mu$, $\beta$, $a$, and $b$. Since $\log x_1\to -\infty$ as $x_1\downarrow0$, this tells us that there exists $\varepsilon = \varepsilon\left(\mu,\beta,a,b\right)>0$ such that $\frac{\partial I}{\partial x_1}\left(x\right)<0$ for $x$ satisfying $x_1 < \varepsilon$. Hence no such $x$ can be a global minimiser.

Restricting the optimisation to $x_1\geq \varepsilon$ (for the $\varepsilon$ used above) we find
\begin{equation}
    \inf_{x : x_1\geq \varepsilon}\left\{\mathcal{I}\left(x\right) + \frac{b\beta}{2}x^2_1\right\} \geq \inf_{x : x_1\geq \varepsilon}\mathcal{I}\left(x\right) + \frac{b\beta}{2}\varepsilon^2 = \frac{b\beta}{2}\varepsilon^2.
\end{equation}
Furthermore, since $\frac{b\beta}{2}x^2_1\geq0$ and any global minimiser of $\mathcal{I}$ has $x_1\geq \varepsilon$,
\begin{equation}
    \inf_{x : x_1 < \varepsilon}\left\{\mathcal{I}\left(x\right) + \frac{b\beta}{2}x^2_1\right\} \geq \inf_{x : x_1< \varepsilon}\mathcal{I}\left(x\right) >0.
\end{equation}
These two bounds then imply that the limit \eqref{eqn:UsedVaradhan} is strictly negative and with the bound \eqref{eqn:PressureDifference} the result is therefore proven.
\end{proof}
%%%%%%%%%%%%%%%%%%%%%%%%%%%%%%%%%%%%%%%%%%%%%%%%%%%%%%%%%%%%%%%%%%%%%%%%%%%%%%%%%%%%%%%%%%%%%%%%%%%%%%%%%%%%%%%%%%%%%%%%%%%%%%%%%%%%%%%%%%
\subsection{Proof of GMF results}\label{SubsectionGMF}
Fix $G\colon [0,\infty)\to\R\cup\{+\infty\}$. Set $L=\inf_y \{G(y)+I(y)\}$ and $M=\{x\colon I(x)+G(x)=L\}\neq \emptyset$.

To aid readability, we restrict ourselves to the case $\abs{M}=1$ and write $\{x_\min\}=M$. We comment on how to generalise the result to multiple minima at the end of the proof. Recall the conditions we require on $G$ in Assumption~\ref{conditionG}.

% We ask for the following condition:
% \begin{enumerate}
%     \item $x\neq \rho_{\mathrm{c}}$ (explain later).
%     \item If $x<\rho_{\mathrm{c}}$, we want $G$ to be twice differentiable at $x$ and continuous in a neighbourhood of $x$. We furthermore require that for any $\e>0$, there exists a $\delta>0$ such that $G^{-1}\left[[K,K+\delta)\right]\subset\bB_\e(x)$.
%     \item If $x>\rho_{\mathrm{c}}$, we require that for any $\e>0$, there exists a $\delta>0$ such that $G^{-1}\left[[K,K+\delta)\right]\subset\bB_\e(x)$. We also require that
%     \begin{enumerate}
%         \item either $G$ has a jump-discontinuity from the left (or the right), is differentiable in a right (resp. left) neighbourhood of $x$ and the first derivative\footnote{at x, we take the right (resp. left) derivative} is bounded uniformly from below.
%         \item or $G$ is twice differentiable at $x$ and there exists $\delta_1>0$ such that for all $y$ in a neighbourhood around $x$
%         \begin{equation}
%             G(x+y)-G(x)\ge \delta_1 y^2\, .
%         \end{equation}
%     \end{enumerate}
% \end{enumerate}
To ease the reading, we abbreviate $\beta G$ by $G$, this multiplicative factor does not affect our calculations. This doesn't affect the proof in any way.\\ We split the proof into two parts. Recall that we shorten $G(x_\min)=K$.

\textbf{Case $x_\min<\rho_{\mathrm{c}}$:} as usual, we start with the partition function. We expand 
\begin{equation}
    \E_{\L,\beta,0}\left[\ex^{-\abs{\L}G\left(\bN\right)}\right]=\E_{\L,\beta,0}\left[\ex^{-\abs{\L}G\left(\bN\right)},\, \bN\in \bB_\e(x_\min)\right]+\E_{\L,\beta,0}\left[\ex^{-\abs{\L}G\left(\bN\right)},\, \bN\notin \bB_\e(x_\min)\right]\, .
\end{equation}
Using Assumption \ref{conditionG}, the second term can be bounded by for some $\e>0$
\begin{equation}
    \E_{\L,\beta,0}\left[\ex^{-\abs{\L}G\left(\bN\right)},\bN\notin \bB_\e(x_\min)\right]=\ex^{-(L+\e)\abs{\L}\left(1+o(1)\right)}\, ,
\end{equation}
and will prove to be negligible. For the first term, we can apply \cite[Theorem 3]{martin1982laplace} to compute
\begin{align}
    \E_{\L,\beta,0}\left[\ex^{-\abs{\L}G\left(\bN\right)},\, \bN\in \bB_\e(x_\min)\right]&\sim\sqrt{1+\frac{G''(x)}{\beta\bmu'(x_\min)}}\ex^{-\abs{\L}\left(G(x_\min)+I(x_\min)\right)}\nonumber\\
    &=\sqrt{1+\frac{G''(x_\min)}{\beta\bmu'(x_\min)}}\ex^{-\abs{\L}L}\, .
\end{align}
Now we can follow the argument made in Section \ref{subsectionlimiting}, to approximate for any test function $F$
\begin{equation}
    \frac{\E_{\L,\beta,0}\left[F(\eta)\ex^{-\abs{\L}G\left(\bN\right)}\right]}{\E_{\L,\beta,0}\left[\ex^{-\abs{\L}G\left(\bN\right)}\right]}\sim \E_{\L,\beta,\bmu(x_\min)}\left[F(\eta)\right]\, .
\end{equation}
This concludes the proof for the case $x_\min<\rho_{\mathrm{c}}$.\qed

\textbf{Case $x_\min>\rho_{\mathrm{c}}$}: for the case $x_\min>\rho_{\mathrm{c}}$, we first treat the case that $G$ has a jump-discontinuity from the left and is twice differentiable from the right. Note that for $x_\min>\rho_{\mathrm{c}}$, it holds that $I(x_\min)=0$ and hence $G(x_\min)=K=L$. Denote $\e_1=\lim_{y\uparrow x_\min}G(y)-K>0$. Fix $\delta>0$ such that
\begin{equation}
    G^{-1}\left[\R\setminus[x_\min+\delta,\infty)\right]\subset (K+\e_1/2,\infty)\, .
\end{equation}
We then expand
\begin{equation}
    \E_{\L,\beta,0}\left[\ex^{-\abs{\L}G\left(\bN\right)}\right]=\E_{\L,\beta,0}\left[\ex^{-\abs{\L}G\left(\bN\right)},\, \bN\in x_\min+[0,\delta)\right]+\Ocal\left(\ex^{-\abs{\L}(L+\e_1/2)}\right)\, .
\end{equation}
The second term will turn out to be negligible. We expand
\begin{equation}
   \E_{\L,\beta,0}\left[\ex^{-\abs{\L}G\left(\bN\right)},\, \bN\in x_\min+[0,\delta)\right]=\sum_{j=0}^{\delta\abs{\L}}\ex^{-\abs{\L}G(x_\min+j/\abs{\L})}\P_{\L,\beta,0}(N_\L=x_\min\abs{\L}+j)\, .
\end{equation}
As $x_\min>\rho_{\mathrm{c}}$, it holds $\P_{\L,\beta,0}(N_\L=x_\min\abs{\L}+j)=\abs{\L}\beta\c_d[\beta(x_\min\abs{\L}+j)]^{-d/2-1}(1+o(1))$. We factor out the dominant terms
\begin{equation}
    \ex^{-\abs{\L}G(x_\min)}\P_{\L,\beta,0}(N_\L=x_\min\abs{\L})\sum_{j=0}^{\delta\abs{\L}}\ex^{-\abs{\L}\left[G(x_\min+j/\abs{\L})-G(x_\min)\right]}\frac{\P_{\L,\beta,0}(N_\L=x_\min\abs{\L}+j)}{\P_{\L,\beta,0}(N_\L=x_\min\abs{\L})}\, .
\end{equation}
As the first derivative is uniformly positive, we can bound $\left[G(x_\min+j/\abs{\L})-G(x_\min)\right]\ge \delta_1j/\abs{\L}$, for some $\delta_1>0$. As the ratio of probabilities in the above equation is uniformly bounded (see \cite{berger2019notes} for this again), this means that the for any $\e_2$, we can find a $J>0$ such that
\begin{equation}
    \sum_{j=J+1}^{\delta\abs{\L}}\ex^{-\abs{\L}\left[G(x_\min+j/\abs{\L})-G(x_\min)\right]}\frac{\P_{\L,\beta,0}(N_\L=x_\min\abs{\L}+j)}{\P_{\L,\beta,0}(N_\L=x_\min\abs{\L})}<\e_2\, .
\end{equation}
For $j\in \{0,\ldots,J\}$,
\begin{equation}
    \frac{\P_{\L,\beta,0}(N_\L=x_\min\abs{\L}+j)}{\P_{\L,\beta,0}(N_\L=x_\min\abs{\L})}=1+o(1)\, .
\end{equation}
Expanding $\left[G(x_\min+j/\abs{\L})-G(x_\min)\right]=G'(x_\min)j/\abs{\L}+o(\abs{\L}^{-1})$, we find that
\begin{equation}
    \sum_{j=0}^{J}\ex^{-\abs{\L}\left[G(x_\min+j/\abs{\L})-G(x_\min)\right]}\frac{\P_{\L,\beta,0}(N_\L=x_\min\abs{\L}+j)}{\P_{\L,\beta,0}(N_\L=x_\min\abs{\L})}\sim  \sum_{j=0}^{\infty}\ex^{-G'(x_\min)j}
    =\frac{1}{1-\ex^{-G'(x_\min)}}\, .
\end{equation}
By letting $\e_2\to 0$, we find that
\begin{equation}
\begin{split}
     \E_{\L,\beta,0}\left[\ex^{-\abs{\L}G\left(\bN\right)}\right]&\sim \frac{\ex^{-\abs{\L}G(x_\min)}\abs{\L}\beta\c_d}{[\beta(x\abs{\L})]^{d/2+1}\left[1-\ex^{-G'(x_\min)}\right]}\\
     &=\frac{\c_d\ex^{-\abs{\L}G(x_\min)}}{(\abs{\L}\beta^{d/2})x_\min^{d/2+1}\left[1-\ex^{-G'(x_\min)}\right]}\, .
\end{split}
\end{equation}
From here on, our proof of Theorem~\ref{thm:GMFlimit} will be very similar to that of \cite[Theorem~2.3]{vogel2021emergence}. Therefore we will focus our attention on the steps that are actually novel for the general mean-field Hamiltonian. Note that we are assuming that $\abs*{\L}$ are such that certain values are integers. For other $\abs*{\L}$, the same argument follows with the introduction of floor or ceiling functions.

We begin with an auxiliary lemma.
\begin{lemma}\label{LemmaPartionFunctionGMF}
If $T=\left[\rho-x_\min\right]\abs{\L}+\Ocal\left(\abs{\L}^{5/6}\right)$, then
\begin{equation}
      \E_{\L,\beta,0}\left[\ex^{-\abs{\L}G\left(\bN\right)}\right]\sim \ex^{-\abs{\L}G(x_\min)} M_{\L}\left[\ex^{-G'\left(\ell(\omega)-T\right)}\1\{\ell(\omega)\ge T\}\right]\, .
\end{equation}
\end{lemma}
The proof of Lemma \ref{LemmaPartionFunctionGMF} is analogous to the above computations and is therefore omitted.

Let $\Theta(N_\L)=x_\min\abs{\L}-N_\L$ and define the probability measure for $\Delta$ a translate of $\L$
\begin{equation}
    \d\Pfrak_\Delta=\frac{1}{\Zfrak_\Delta}\exp\left\{-G'\left[\ell(\omega)-{\Theta(N_\L)}\right]\right\}\1\{\ell(\omega)\ge \Theta(N_\L)\}\d\P_\Delta\otimes M_\Delta\, .
\end{equation}
% where we set $\bB_{\delta}^+(x,\rho)=\bB_{\delta\abs{\L}^{5/6}}(\left[\rho-x_\min\right]\abs{\L})$. 
It follows from the previous lemma that
\begin{cor}
Under the above conditions,
\begin{equation}
     \E_{\L,\beta,0}\left[\ex^{-\abs{\L}G\left(\bN\right)}\right]\sim\ex^{-\abs{\L}G(x_\min)} \Zfrak_\Delta\, .
\end{equation}
\end{cor}
In the next lemma, we remove the influence of the Hamiltonian,
\begin{lemma}\label{RemoveHamil}
It holds that
\begin{equation}
    \sum_{y\in C_N}\Pfrak_{yN+\L}\left[\Big|\ex^{-\abs{\L}G(\bN)+\abs{\L}G(x)+G'\left[\ell(\omega)-{\Theta(N_\L)}\right]}-1\Big|\1\{\eta\cap\{0\}\neq \emptyset\}\right]=o(1)\, .
\end{equation}
\end{lemma}
The proof of Lemma \ref{RemoveHamil} is almost analogous to the proof of \cite[Lemma 5.14]{vogel2021emergence}. The only difference is that in that case, we had an explicit remainder of a square, whereas in our case we can only bound $G(x_\min+\delta)-G(x_\min)-\delta G'(x_\min)$ by a $o(\delta)$ term. This only gives a decay of $o(1)$ (as opposed to $o(\abs{\L}^{-1}$)) but this is enough for our purposes.
\begin{lemma}\label{LemmaConvRandInt}
For $\Theta\in\R$, we define the probability measure $\Mfrak^\Theta_\Delta$
\begin{equation}
    \d\Mfrak^\Theta_\Delta(\omega)=\frac{1}{\Zfrak(\Theta)}\ex^{-G'\left[\ell(\omega)-\Theta\right]}\1\{\ell(\omega)\ge \Theta\}\d M_\Delta(\omega)\, .
\end{equation}
We furthermore set for $\Theta=(\Theta_y)_{y\in C_N}$
\begin{equation}
    \d\Mfrak^\Theta_K=\sum_{y\in C_N}\1\{\omega\cap K\neq\emptyset\}\d\Mfrak^{\Theta_y}_\Delta\, ,
\end{equation}
for $K\subset \R^d$ compact.

Define ${\Mfrak}_K^{\Theta,*}={\Mfrak}_K^\Theta\circ\amalg$, as in \cite{vogel2021emergence}. We then have that there for every sequence $\Theta_y\in \rho_{\mathrm{c}}\abs{\L}+\Ocal\left(\abs{\L}^{5/6}\right)$ that
\begin{equation}
    {\Mfrak}_K^{\Theta,*}[E]=\nu[E]\left(1+o(1)\right)\, ,
\end{equation}
where the $o(1)$ can be chosen uniform in $\Theta$ and $E$ is an element of the dense approximating class defined in \cite[Definition 5.7]{vogel2021emergence}.
\end{lemma}
These lemmas fill in the sections of the proof of \cite[Theorem~2.3]{vogel2021emergence} that extend that result for $x_{\min} > \rho_{\mathrm c}$ to general mean-field interaction with jump-discontinuity from the left. The case with the jump-discontinuity for the right is analogous to the one treated above. \qed

Suppose now that $G$ does not have a jump-discontinuity and is twice differentiable around $x_\min$. We expand
\begin{equation}
    \E_{\L,\beta,0}\left[\ex^{-\abs{\L}G\left(\bN\right)}\right]=\E_{\L,\beta,0}\left[\ex^{-\abs{\L}G\left(\bN\right)},\, \bN\in \bB_\delta(x_\min)\right]+\Ocal\left(\ex^{-\abs{\L}(K+\e_1/2)}\right)\, .
\end{equation}
The second term will turn out to be negligible. By the virtue of Assumption \ref{conditionG},
\begin{equation}
    \left[G(x_\min+j/\abs{\L})-G(x_\min)\right]\ge \delta_1\frac{j^2}{\abs{\L}^2}\, .
\end{equation}
Therefore, for any $\e_2>0$, we can find and $R>0$ such that
\begin{equation}
    \left|\sum_{j=-\delta\abs{\L}}^{\delta\abs{\L}}\ex^{-\abs{\L}\left[G(x_\min+j/\abs{\L})-G(x_\min)\right]}-\sum_{j=-R\sqrt{\abs{\L}}}^{R\sqrt{\abs{\L}}}\ex^{-\abs{\L}\left[G(x_\min+j/\abs{\L})-G(x_\min)\right]}\right|<\e_2\, .
\end{equation}
Using the Riemann approximation on the scale $\sqrt{\abs{\L}}$, we find that
\begin{equation}
    \sum_{j=-R\sqrt{\abs{\L}}}^{R\sqrt{\abs{\L}}}\ex^{-\abs{\L}\left[G(x_\min+j/\abs{\L})-G(x_\min)\right]}\sim \sum_{j=-R\sqrt{\abs{\L}}}^{R\sqrt{\abs{\L}}}\ex^{-\frac{G''(x_\min)j^2}{2\abs{\L}}}\sim \sqrt{\abs{\L}}\int_{-R}^R\ex^{-\frac{G''(x_\min)t^2}{2}}\d t\, .
\end{equation}
By letting $R\to\infty$, we conclude that
\begin{equation}
    \E_{\L,\beta,0}\left[\ex^{-\abs{\L}G\left(\bN\right)}\right]\sim \frac{\c_d \sqrt{\pi\abs{\L}}\ex^{-\abs{\L}G(x_\min)}}{(\abs{\L}\beta^{d/2})x_\min^{d/2+1}\sqrt{G''(x_\min)}}\, .
\end{equation}
From here one, the proof works in the same way as the case for the jump-discontinuity. Having covered all the cases of Theorem \ref{thm:GMFlimit}, we conclude the proof.\qed
\begin{remark}\label{remarkmin}
What happens if $I+G$ has more than one minimizer, i.e. $\abs{M}>1$? The condition $G^{-1}\left[[K,K+\delta)\right]\subset\bB_\e(x)$ has to be replaced by $G^{-1}\left[[K,K+\delta)\right]\subset\bigcup_{x\in M}\bB_\e(x)$. We then follow the same procedure, expanding the partition function around each neighbourhood for each point. This requires $M$ to not have accumulation points. The final result will be a weighted mixture of different loop soups with individual intensities corresponding to the values induced by $x$ for $x\in M$. We leave the details to the reader.
\end{remark}
\section{Appendix}\label{SectionAppendix}
{
\begin{table}[H]
  \renewcommand{\arraystretch}{1.4}
 \caption{List of frequently used notation}
\begin{tabular}{ |p{1.2cm}||p{5.6cm}|p{4.9cm}|p{3.65cm}|  }
 \hline
% \multicolumn{4}{|c|}{} \\
% \hline
 Symbol & Definition & Explanation &Class\\
 \hline \hline
 $\c_d$   &  $(2\pi)^{-d/2}$   & &   Constant\\
 \hline
 $\rho$   &     &Density&   Model parameter\\
 \hline
  $\mu$   &     &Chemical potential&   Model parameter\\
 \hline
 $\beta$   &     &Inverse temperature&   Model parameter\\
 \hline
 $a,b$   &   See Eq. \eqref{eqn:partialHYLinteraction}  &Interaction strength&   Model parameter\\
 \hline
 $\rho_\mathrm{c}$   &  $\sum_{j\ge 1}(\beta j)^{-d/2}=\zeta(d/2)\beta^{-d/2}$   &Critical density&   Parameter\\
 \hline
 $\rho_\mathrm{c}^\hy$   &   See Eq. \eqref{EqdefrhocHY}  &HYL-critical density&   Parameter\\
 \hline
 $\rho_\mathrm{e}$   &  $\max\{\rho-\rho_\mathrm{c},0\}$   &Excess density&   Parameter\\
 \hline
 $\rho_\mathrm{S}$   &  See Eq. \eqref{EquationrhoS}   &Density from interaction&   Parameter\\
 \hline
 $\rhob$   &  $\rho_\mathrm{e}+\rho_\mathrm{S}$   &Total condensate density&   Parameter\\
 \hline 
 $\rho^{\mathrm{GC}}$   &  See Eq. \eqref{EquationrhoGC}   &Grand-can. condensate dens.&   Parameter\\
 \hline
 $P$   & $P(x)\!=\!\beta^{-d/2}\c_d\sum_{j\ge 1}{\ex^{\beta x j}}{j^{-1-d/2}}$    &Pressure {\emph{(non-interacting)}}&   Thermodyn. function\\
 \hline
 $\brho$   & $\brho(x)\!=\!\beta^{-d/2}\c_d\sum_{j\ge 1}{\ex^{\beta x j}}{j^{-d/2}}$    &Density {\emph{(non-interacting)}}&   Thermodyn. function\\
 \hline
 $\bmu$   & $\bmu(x)\!=\!\brho^{-1}(x)$    &Chemical potential&   Thermodyn. function\\
 \hline
 $I$ & See Eq. \eqref{EquationEcplicitI} & Rate function &  Thermodyn. function\\
 \hline
\end{tabular}\label{table1} 
\end{table}}
\subsection*{Acknowledgements} The authors would like to thank the anonymous referees for their many suggestions. Quirin Vogel would like to thank Julius Damarackas for his help with improving the presentation of the article. He would further like to thank Roberto Fernandez, Vedran Sohinger and Daniel Ueltschi for the discussions on this topic.

\bibliography{thoughts}{}
\bibliographystyle{alpha}
\end{document}